\documentclass[12pt]{amsart}
\usepackage{amscd,verbatim}
\usepackage{amssymb}
\usepackage[all]{xy}
\usepackage[colorlinks,linkcolor=blue,citecolor=blue,urlcolor=red]{hyperref}

\newcounter{spec}
{\end{list}}%

\newtheorem{lemma}{Lemma}[section]
\newtheorem{thm}[lemma]{Theorem}

\newtheorem{prop}[lemma]{Proposition}

\newtheorem{cor}[lemma]{Corollary}

\theoremstyle{definition}
\newtheorem{defn}[lemma]{Definition}

\newtheorem{definition}[lemma]{Definition}

\theoremstyle{remark}

\newtheorem{rk}[lemma]{Remark}
\newtheorem{remark}[lemma]{Remark}

\newtheorem{claim}[lemma]{Claim}

\numberwithin{equation}{section}

\newcounter{zaehler} 
\numberwithin{equation}{section}

\newcommand{\A}{\mathbf{A}}

\renewcommand{\P}{\mathbf{P}}

\newcommand{\Z}{\mathbb{Z}}

\newcommand{\sC}{\mathcal{C}}
\newcommand{\sD}{\mathcal{D}}

\newcommand{\sH}{\mathcal{H}}

\newcommand{\sL}{\mathcal{L}}
\newcommand{\sO}{\mathcal{O}}

\newcommand{\sS}{\mathcal{S}}
\newcommand{\sT}{\mathcal{T}}
\newcommand{\sK}{\mathcal{K}}

\newcommand{\sU}{\mathcal{U}}
\newcommand{\sV}{\mathcal{V}}
\newcommand{\fX}{\mathfrak{X}}
\newcommand{\sX}{\mathcal{X}}

\newcommand{\sY}{\mathcal{Y}}
\newcommand{\fY}{\mathfrak{Y}}

\newcommand{\sZ}{\mathcal{Z}}
\newcommand{\fZ}{\mathfrak{Z}}

\newcommand{\bZ}{\mathbb{Z}}
\newcommand{\fm}{\mathfrak{m}}

\newcommand{\Xb}{{\overline{X}}}
\newcommand{\Yb}{{\overline{Y}}}

\newcommand{\Lb}{{\overline{L}}}
\newcommand{\Zb}{{\overline{Z}}}
\newcommand{\Hb}{{\overline{H}}}

\newcommand{\Cor}{\operatorname{\mathbf{Cor}}}

\newcommand{\HI}{\operatorname{\mathbf{HI}}}

\newcommand{\Rec}{{\operatorname{\mathbf{Rec}}}}

\newcommand{\RSC}{{\operatorname{\mathbf{RSC}}}}
\newcommand{\RSCNis}{{\operatorname{\mathbf{RSC}}}_{\Nis}}

\newcommand{\ul}[1]{{\underline{#1}}}

\newcommand{\PST}{{\operatorname{\mathbf{PST}}}}

\newcommand{\NST}{\operatorname{\mathbf{NST}}}

\newcommand{\Hom}{\operatorname{Hom}}
\newcommand{\Nm}{\operatorname{Nm}}
\newcommand{\uHom}{\operatorname{\underline{Hom}}}

\newcommand{\Ker}{\operatorname{Ker}}

\newcommand{\Coker}{\operatorname{Coker}}

\newcommand{\Pic}{\operatorname{Pic}}

\newcommand{\Spec}{\operatorname{Spec}}
\newcommand{\Spf}{\operatorname{Spf}}

\newcommand{\Sm}{\operatorname{\mathbf{Sm}}}
\newcommand{\Sch}{\operatorname{\mathbf{Sch}}}
\newcommand{\Ab}{\operatorname{\mathbf{Ab}}}

\newcommand{\tr}{{\operatorname{tr}}}

\newcommand{\fin}{{\operatorname{fin}}}

\newcommand{\op}{{\operatorname{op}}}

\newcommand{\Nis}{{\operatorname{Nis}}}
\newcommand{\et}{{\operatorname{\acute{e}t}}}

\newcommand{\inj}{\hookrightarrow}

\newcommand{\id}{{\operatorname{id}}}

\newcommand{\codim}{{\operatorname{codim}}}
\newcommand{\ch}{{\operatorname{ch}}}

\newcommand{\CH}{{\operatorname{CH}}}

\renewcommand{\lim}{\operatornamewithlimits{\varprojlim}}
\newcommand{\colim}{\operatornamewithlimits{\varinjlim}}

\newcommand{\ol}{\overline}

\renewcommand{\phi}{\varphi}
\renewcommand{\epsilon}{\varepsilon}
\renewcommand{\div}{\operatorname{div}}

\newcommand{\MNST}{\operatorname{\mathbf{MNST}}}

\newcommand{\MCor}{\operatorname{\mathbf{MCor}}}

\newcommand{\MSm}{\operatorname{\mathbf{MSm}}}
\newcommand{\MP}{\operatorname{\mathbf{MSm}}}

\newcommand{\MPST}{\operatorname{\mathbf{MPST}}}
\newcommand{\CI}{\operatorname{\mathbf{CI}}}

\newcommand{\CIt}{{^{\tau}{\CI}}}
\newcommand{\CItsp}{\CIt^{sp}}

\newcommand{\lsCIt}{\CIt^{ls}}

\newcommand{\lsCItspNis}{\CIt^{ls,sp}_{\Nis}}
\newcommand{\CItspNis}{\CIt^{sp}_{\Nis}}

\newcommand{\lsCItsp}{\CIt^{ls,sp}}

\newcommand{\Bl}{{\mathbf{Bl}}}

\newcommand{\bcube}{{\ol{\square}}}
\newcommand{\cube}{\square}
\newcommand{\lscube}{\operatorname{ls-{\ol{\square}}}}

\def\fin{\operatorname{fin}}

\renewcommand{\Lb}{{\overline{L}}}
\renewcommand{\Zb}{{\overline{Z}}}

\def\indlim#1{\underset{{\underset{#1}{\longrightarrow}}}{\mathrm{lim}}\; }
\def\rmapo#1{\overset{#1}{\longrightarrow}}
\def\lmapo#1{\overset{#1}{\longleftarrow}}

\newcommand{\ulMPST}{\operatorname{\mathbf{\underline{M}PST}}}
\newcommand{\ulMPSTfin}{\ulMPST^{\fin}}
\newcommand{\ulMNST}{\operatorname{\mathbf{\underline{M}NST}}}
\newcommand{\ulMNSTfin}{\ulMNST^{\fin}}
\newcommand{\ulMNSTls}{\operatorname{\mathbf{\underline{M}NST}}_{ls}}
\newcommand{\ulMCor}{\operatorname{\mathbf{\underline{M}Cor}}}
\newcommand{\ulMCorfin}{\ulMCor^{\fin}}
\newcommand{\ulMCorls}{\ulMCor_{ls}}

\newcommand{\ulomega}{\underline{\omega}}

\newcommand{\Comp}{\operatorname{\mathbf{Comp}}}
\newcommand{\CompM}{\Comp(M)}

\def\bZ{\mathbb{Z}}
\def\Ztr{\bZ_\tr}

\def\Image{\mathrm{Image}}

\def\pb{\overline{p}}

\def\isom{\overset{\cong}{\longrightarrow}}

\def\tW{\widetilde{W}}

\def\Xinf{X_{\infty}}

\def\Cinf{C_{\infty}}
\def\Yinf{Y_{\infty}}
\def\Linf{L_{\infty}}
\def\Zinf{Z_{\infty}}
\def\Hinf{H_{\infty}}

\def\tXinf{\widetilde{X}_{\infty}}

\def\tX{\widetilde{X}}
\def\tW{\widetilde{W}}
\def\tC{\widetilde{C}}

\def\piinf{\pi_{\infty}}
\def\ttheta{\tilde{\theta}}
\def\thetarelf{\theta_{f,rel}}
\def\thetarelg{\theta_{g,rel}}
\def\tthetarelf{\ttheta_{f,m}}
\def\tthetarelg{\ttheta_{g,m}}

\newcommand{\ulMCorScube}{\operatorname{\mathbf{\underline{M}Cor}_S^{\square}}}
\newcommand{\ulMCorS}{\operatorname{\mathbf{\underline{M}Cor}_S}}

\def\qaq{\;\text{ and }\;}
\def\qwith{\;\text{ with }\;}

\def\qfor{\;\text{ for }\;}

\def\ep{\epsilon}
\def\qwith{\;\;\text{ with }}

\def\cd#1#2{{#1}^{(#2)}}

\def\Fc{F_{-1}}
\def\Fcc{F_{-2}}
\def\Fccc#1{F_{-#1}}
\def\tFccc#1{\tilde{F}_{-#1}}
\def\Fcont{F_{-1}}
\def\Fcontt#1{F^{(#1)}_{-1}}

\def\ttimes{\tilde{\times}}

\def\phib{\overline{\phi}}

\def\Cbb{{\overline{C}}}

\def\Txinf{T_{x,\infty}}

\def\sCFU#1{\sC^{#1}_{\sU}(F)}

\def\sDFUphii#1#2{\sD^{#1}_{\sU}(F)_{\phi_{#2}}}

\def\sDFUphih#1#2{\sD^{#1}_{\sU}(F)^{hor}_{\phi_{#2}}} 
\def\sDFUphiv#1#2{\sD^{#1}_{\sU}(F)^{ver}_{\phi_{#2}}} 

\def\ep{\epsilon}
\def\tul{\underline{t}}

\def\hep{\widehat{\epsilon}}
\def\ul#1{{\underline{#1}}}
\def\hen#1#2{#1^{h}_{|#2}}

\def\hF{\hat{F}}
\def\hgamma{\hat{\gamma}}

\def\dphi#1{\delta_{\phi_i}^{#1}}
\def\dphii#1{\delta_{\phi_{i-1}}^{#1}}

\def\parphiw#1{\partial^{#1}_{\phi_i,w}}

\def\parphi#1{\partial^{#1}_{\phi_i}}

\def\tauul{\underline{\tau}}

\def\tul{\underline{t}}

\def\hM#1{h_0^{\bcube}(#1)}
\def\hMM#1{h^0_{\bcube}(#1)}
\def\hMMnp{h^0_{\bcube}}

\def\hMw#1{h_0(#1)}
\def\ta{\tilde{a}}
\def\ulaNis{\underline{a}_{\Nis}}
\def\ulaNisfin{\underline{a}_{\Nis}^{\fin}}
\def\aVNis{{a}^V_{\Nis}}

\def\aNis{a_{\Nis}}

\def\omegaCI{\omega^{\CI}}

\def\tS{\tilde{S}}
\def\tnu{\tilde{\nu}}

\def\tSm{\widetilde{\Sm}}
\def\ulsO{\underline{\sO}}

\def\Sigmafin{\Sigma^{\fin}}

\def\lssemipure{\text{semipure }}

\begin{document}
\title{Purity of reciprocity sheaves}
\author{Shuji Saito}
\address{Interactive Research Center of Science\\
Graduate School of Science and Engineering\\
Tokyo Institute of Technology\\
2-12-1 Okayama, Meguro\\ Tokyo 152-8551\\ Japan}
\email{sshuji@msb.biglobe.ne.jp}
\date{}
\thanks{The author is supported by JSPS KAKENHI Grant (15H03606) and
the SFB grant ``Higher Invariants".}

\maketitle


\tableofcontents

\begin{abstract}
The purpose of this paper is to prove a conjecture of Kahn-Saito-Yamazaki 
\cite[Conj.1(1)]{rec}.
This is accomplished by extending Voevodsky's fundamental results on homotopy invariant (pre)sheaves with transfers \cite{voetri0} to its generalizations,
\emph{reciprocity sheaves} and \emph{cube-invariant sheaves} in the context of theory of \emph{modulus (pre)sheave with transfers} \cite{kmsy}.
The main results of this paper is expected to play a crucial role in deducing 
the main properties of the triangulated category of \emph{motives with modulus}, 
which is a new triangulated category enlarging Voevodsky's triangulated category of motives to encompass non homotopy invariant motivic phenomena. 
\end{abstract}


\section*{Introduction}\label{intro}

Let $\Sm$ be the category of separated smooth schemes of finite type over $k$.
Let $\Cor$ be the category of finite correspondences:
$\Cor$ has the same objects as $\Sm$ and morphisms in $\Cor$ are finite correspondences. 
Let $\PST$ be the category of additive presheaves of abelian groups on $\Cor$, called presheavs with transfers. 
In Voevodsky's theory of his triangulated category of motives, a fundamental role is played by homotopy invariant objects $F\in \PST$, namely such $F$ that $F(X)\to F(X\times\A^1)$ induced by the projection $X\times\A^1\to X$ are isomorphisms for all $X\in \Sm$.
The homotopy invariant objects form a full abelian subcategory $\HI\subset\PST$.
\medbreak

In order to extend Voevodsky's paradigm to a non-homotopy invariant framework, we use a new full abelian subcategory $\RSC\subset\PST$ of \emph{reiprocity presheaves} 
\footnote{The terminology ``reciprocity presheaves" was used in \cite{rec} for a slightly different notion, which is not used in this paper. 
In loc.cite. a full subcategory $\Rec$ of $\PST$ was introduced.
In \cite{recII} it is shown that $\RSC\subset \Rec$ and $\RSC\cap \NST =\Rec\cap \NST$.}.
It contains $\HI$ and many objects of $\PST$ which are not in $\HI$,
such as the presheaf associated to a commutative algebraic group $G$ (which may contains a unipotent part) and the presheaf $\Omega^i$ of K\"ahler differential forms and the de Rham-Witt presheaf $W_n\Omega^i$.

Let $\Ztr(X)\in \PST$ be the object which $X\in \Sm$ represents by the Yoneda functor. Recall (\cite[Lem. 2.16]{mvw}) that the inclusion $\HI\to \PST$ has a left adjoint $h_0$ such that $h_0(X):=h_0(\Ztr(X))$ is given by
\begin{equation}\label{eq;h0X}
\begin{aligned}
 h_0(X)(Y) &=\Coker\big(\Ztr(X)(Y\times\A^1)\rmapo{i_0^*-i_1^*} \Ztr(X)(Y)\big)
 \quad (Y\in \Sm)\\
&=\Coker\big(\Cor(Y\times\A^1,X)\rmapo{i_0^*-i_1^*} \Cor(Y,X)\big),\\
\end{aligned}
\end{equation} 
where $i_\ep^*$ for $\ep=0,1$ is the pullback by the section $i_\ep: \Spec(k) \to \A^1$. This implies that $F\in \PST$ is in $\HI$ if and only if for any $X\in \Sm$ and $a\in F(X)$, the map $a:\Ztr(X) \to F$ in $\PST$ associated to $a$ by the Yoneda functor, factors through the quotient $h_0(X)$ of $\Ztr(X)$. 
\medbreak

The key idea to define $\RSC$ is to introduce bigger quotients $h_0(\fX)$ of $\Ztr(X)$ associated to 
pairs $\fX=(\Xb,\Xinf)$ where $\Xb$ is a proper scheme over $k$
and $\Xinf$ is an effective Cartier divisor on $\Xb$ such that $X=\Xb - |\Xinf|$.
We have for $Y\in \Sm$
\begin{equation}\label{eq;h0sX;intro}
\hMw \fX(Y)=\Coker\big(\ulMCor(Y\otimes\bcube,\fX)\rmapo{i_0^*-i_1^*} \Cor(Y,X)\big),
\end{equation} 
where $\ulMCor(Y\otimes\bcube,\fX)$ is a subgroup of $\Cor(Y\times\A^1,X)$
generated by elementary correspondences\footnote{It means integral closed subschemes in $Y\times X$.} satisfying a certain admissibility condition
with respect to $\Xinf$.\footnote{
We have surjective maps $\Ztr(X)\to \hMw \fX \to h_0(X)$ in $\PST$.
Evaluated on $\Spec(k)$, it is identified with 
$Z_0(X)\to \CH_0(\Xb,\Xinf) \to H_0^S(X)$, where $Z_0(X)$ is the group of $0$-cycles on $X$, $\CH_0(\Xb,\Xinf)$ is the Chow group of $0$-cycles with modulus introduced in \cite{KSa} and $H_0^S(X)$ is the Suslin homology of $X$ introduced in \cite{sv-invent}. See Definition \ref{def;hM} for a more conceptual definition of $\hMw \fX$.}
Then $F\in \PST$ is defined to be a reciprocity presheaf, or to have reciprocity, if for any $X\in \Sm$ and $a\in F(X)$, the associated map $a:\Ztr(X) \to F$ factors through $h_0(\fX)$ for some $\fX$ as above. 
By \eqref{eq;h0X} and \eqref{eq;h0sX;intro}, $h_0(X)$ is a quotient of $h_0(\fX)$, and hence every $F\in \HI$ has reciprocity.\footnote{Heuristically 
$\RSC$ (resp. $\HI$) may be viewed as consisting of such $F\in \PST$ that 
for any $X\in \Sm$ and $a\in F(X)$, $a$ has ``bounded (resp. tame) ramification" along the boundary of a compactification of $X$. 
A manifestation of this viewpoint is given in \cite{rs}.}
By the definition those $F\in \PST$ which have reciprocity form a full abelian subcategory $\RSC\subset\PST$ which is closed under subobjects
and quotients in $\PST$. 
\medbreak

Let $\NST\subset \PST$ be the full subcategory of Nisnevich sheaves, i.e.
those objects $F\in \PST$ whose restrictions $F_X$ to the \'etale site $X_{\et}$ over $X$ are Nisnevich sheaves for all $X\in \Sm$. 
By \cite[Th. 3.1.4]{voetri} the inclusion $\NST\to \PST$ has an exact left adjoint $\aVNis$ such that $\aVNis F$ is the Nisnevich sheafication $F_{\Nis}$ of $F$
as a presheaf on $\Sm$.
We put $\RSCNis=\RSC\cap \NST$.
\medbreak

We now state our main results on $\RSC$.
Theorem \ref{thm;purityrec} gives an affirmative answer to \cite[Conjecture 1(1)]{rec}. 

\begin{thm}\label{thm0;rec}
The functor $\aVNis$ preserves $\RSC$.
\end{thm}

An analogous result is obtained in \cite{rec} by a different method. 
\medbreak

To state the second main result on $\RSC$, we introduce some notations.
For $F\in \PST$ and $n\in \Z_{>0}$ and $S\in \Sm$, define
\[ \tFccc n (S)=\Coker\big(\underset{1\leq i\leq n}{\bigoplus}\;
F((\A^1-0)^{i-1}\times \A^1 \times(\A^1-0)^{n-i}\times S)\to  F((\A^1-0)^n\times S)\big)
,\]
where $\A^1$ is the affine line over $k$ with $0\in \A^1$ the origin. 
\medbreak

For $X\in \Sm$ and $n\in \Z_{\geq 0}$, let $\cd X n$ be the set of points $x\in X$ 
such that the closure of $x$ in $X$ is of codimension $n$.
Assuming $k$ is perfect, there is an isomorphism (see Lemma \ref{lem;etalepoint})
\begin{equation}\label{eq;epsilon}
\ep :\hen{X}{x}\simeq  \Spec K\{t_1,\dots,t_n\},
\end{equation}
where $K=k(x)$ and $\hen{X}{x}$ is the henselization of $X$ at $x$ and $(t_1,\dots,t_n)$ is a regular system of parameters of $X$ at $x$, and
$K\{t_1,\dots,t_n\}$ is the henselization of $K[t_1,\dots,t_n]$ at 
$(t_1,\dots,t_n)$.
\medbreak


\begin{thm}\label{thm;purityrec}
Assume $k$ is perfect. Let $F\in \RSCNis$.
For $X\in \Sm$ and $x\in \cd X n$ with $n\in \Z_{>0}$, we have
\begin{equation}\label{eq1;thm;purityrec}
 H^i_x(X_{\Nis},F_X) =0\qfor i\not= n,
\end{equation}
and there exists an isomorphism depending on $\ep$ from \eqref{eq;epsilon}:
\begin{equation}\label{eq2;thm;purityrec}
\theta_\ep: \tFccc n (x) \simeq H^n_x(X_{\Nis},F_X).
\end{equation}
\end{thm}


As an immediate corollary we get the following.

\begin{cor}\label{cor;purityrec}
Let $F\in \RSCNis$. Let $\sX$ be the henselization of $X\in \Sm$ at a point of 
$ X$ and $\xi$ be its generic point. Then the Cousin complex
\begin{multline}\label{eq;Cousin}
0\to F(\sX)\to F(\xi)\to \bigoplus_{x\in \sX^{(1)}} H^1_{x}(\sX_{\Nis}, F_\sX)\to \\
\dots  \to \bigoplus_{x\in \sX^{(n)}} H^n_{x}(\sX_{\Nis}, F_\sX)\to \dots
\end{multline}
is exact.
\end{cor}

\bigskip

After replacing $\RSC$ by $\HI$, the above results were proved by
Voevodsky \cite{voetri0} and played a fundamental role in his theory of triangulated category of mixed motives in \cite{voetri} (see Theorem \ref{thm2-intro;purityM}).
\medbreak

\subsection{Modulus refinements}

We now explain refinements of the above results in the new categorical framework developed in \cite{kmsy}. A new category $\ulMCor$ (see Definition \ref{def;moduluspair}) is introduced: The objects are \emph{modulus pairs} $\fX=(\Xb,\Xinf)$ where $\Xb$ is a separated schemes of finite type over $k$ equipped with an effective Cartier divisor $\Xinf\subset \Xb$ such that $\Xb - |\Xinf|\in \Sm$.
The morphisms are finite correspondences satisfying some admissibility and properness conditions. 
Let $\MCor\subset \ulMCor$ be the full subcategory of such objects $(\Xb,\Xinf)$
that $\Xb$ is proper over $k$. 
We then define $\ulMPST$ (resp. $\MPST$) as the category of additive presheaves of abelian groups on $\ulMCor$ (resp. $\MCor$). 
For $F\in \ulMPST$ and $\fX=(\Xb,\Xinf)\in \ulMCor$ write $F_{\fX}$ for the presheaf
on the \'etale site $\Xb_{\et}$ over $\Xb$ given by $U\to F(\fX_U)$ for $U\to \Xb$ \'etale, where $\fX_U=(U,\Xinf\times_{\Xb}U)\in \ulMCor$.
We say $F$ is a Nisnevich sheaf if so is $F_{\fX}$ for all $\fX\in \ulMCor$.
We write $\ulMNST\subset \ulMPST$ for the full subcategory of Nisnevich sheaves and 
\begin{equation}\label{eq;HfXF}
 H^i(\fX_{\Nis},F) = H^i(\Xb_{\Nis},F_{\fX})\qfor F\in \ulMNST.
\end{equation}
It is shown in \cite{kmsy} (see Theorem \ref{thm;ulMNST}(2) in \S\ref{Mrec}) that the inclusion $\ulMNST\to \ulMPST$ has an exact left adjoint $\ulaNis$.\footnote{For $F\in \ulMPST$ and $\fX \in \ulMCor$, $(\ulaNis F)_{\fX}$ is NOT the Nisnevich sheafication of $F_{\fX}$ contrary to the case of $\aVNis$ (see Theorem \ref{thm;ulMNST}(2) for its description). }
We have a functor 
\[\omega:\ulMCor \to \Cor\;;\; (\Xb,\Xinf) \to \Xb - |\Xinf|,\]
and two pairs of adjunctions
\begin{equation}\label{eq;tauomegaadjunction-intro}
\MPST\begin{smallmatrix} \tau^*\\ \longleftarrow\\ \tau_!\\ \longrightarrow\\
\end{smallmatrix}\ulMPST, 
\quad
\MPST
\begin{smallmatrix} \omega^*\\ \longleftarrow\\ \omega_!\\ \longrightarrow\\
\end{smallmatrix}\PST,
\end{equation}
where $\tau^*$ is induced by the natural inclusion $\tau:\MCor\to \ulMCor$ and
$\tau_!$ is its left Kan extension, and $\omega^*$ is induced by $\omega$ and $\omega_!$ is its left Kan extension (see Lemma \ref{lem;MPST}(4) and (5) for descriptions of these left Kan extensions).
We now introduce a basic property which is an analogue of homotopy invariance exploited by Voevodsky: Put $\bcube=(\P^1,\infty)$. 
$F\in \MPST$ is called $\bcube$-invariant if $F(\fX)\simeq F(\fX\otimes\bcube)$ for all $\fX\in \MCor$ (see Definition \ref{def;moduluspair}(4) for the tensor product $\otimes$ in $\MCor$).
We write $\CI\subset \MPST$ for the full subcategory of $\bcube$-invariant objects.
The category $\CI$ and its essential image $\CIt$ under $\tau_!:\MPST\to\ulMPST$ will play a fundamental role in this paper.
Let $\CItsp\subset \CIt$ be the full subcategory of \emph{semipure} objects $F$,
namely such objects that the natural map
$F(\Xb,\Xinf)\to F(\Xb-\Xinf,\emptyset)$ is injective for all $\fX=(\Xb,\Xinf)\in \ulMCor$.
One can show (see Lemma \ref{lem2;CIRec}) $\omega_!(\CI)= \RSC$ and that the induced functor $\CI\to \RSC$ admits a fully faithful right adjoint $\omegaCI$ such that $\tau_!\omegaCI(\RSC)\subset \CItsp$.
Using this, Theorems \ref{thm0;rec} and \ref{thm;purityrec}
are deduced from Theorems \ref{thm-intro;localinjectivity} and \ref{thm-intro;purityM} 
below respectively. 

\begin{thm}\label{thm-intro;localinjectivity}
(Theorem \ref{thm;sheafication})
The functor $\ulaNis$ preserves $F\in \CItsp$. 
\end{thm}

Write $\CItspNis=\CItsp\cap \ulMNST$.

\begin{thm}\label{thm-intro;purityM}(Corollaries \ref{cor;vanishing})
Assume $k$ is perfect. Take $F\in \CItspNis$ and $\fX=(\Xb,\Xinf)\in \ulMCor$.
Assume $\Xb$ and $|\Xinf| \in \Sm$.
For $x\in \cd {\Xb} n$ with $K=k(x)$, we have
\begin{equation*}\label{eq1;thm-intro;purityM} 
H^i_x(\fX_{\Nis},F)=0 \qfor i\not=n,
\end{equation*}
where $H^i_x(\fX_{\Nis},F)$ is defined as \eqref{eq;HfXF} for cohomology with support. 
\end{thm}
\medbreak

We also show the following.

\begin{thm}\label{thm2-intro;purityM}(Theorem \ref{thm2;purityM})
Assume $k$ is perfect. 
Take $F\in \CItspNis$ and $\fX=(\Xb,\Xinf)\in \ulMCor$.
Assume that $\Xb\in \Sm$ and $|\Xinf|$ is a simple normal crossing divisor on $\Xb$.
Then we have an isomorphism
\begin{equation}\label{eq;thm2-intro;purityM}
 \pi^*: H^q(\fX_{\Nis},F) \isom H^q((\fX\otimes\bcube)_{\Nis},F) 
\end{equation}
induced by the projection $\pi:\fX\otimes\bcube \to \fX$.
\end{thm}

The following theorem of Voevodsky is a direct consequence of Theorems \ref{thm-intro;purityM} and \ref{thm2-intro;purityM} (see \S\ref{RSCproof}).

\begin{thm}\label{thm;Voev}(\cite[Th.5.6]{voetri0})
Assume $k$ is perfect. 
\[ H^i(X_{\Nis},F_X) \to H^i((X\times\A^1)_{\Nis},F_{X\times\A^1}) \]
induced by the projection $X\times\A^1\to X$ is an isomorphism.
\end{thm}

\bigskip

We give an overview of the content of the paper.

In \S\ref{Mrec} we review basic notions and facts on modulus sheaves with transfer and reciprocity sheaves. The whole content of \S1 is a joint work with B. Kahn, H. Miyazaki and T. Yamazaki. The content of \S1.1 is extracted from \cite{kmsy} except that of \S1.2 through \S1.4. The (enriched version of) content of those subsections will appear in a forthcoming joint paper while the proofs of all the statements used in this paper are given in \S1. 

In \S\ref{Vpair} we introduce \emph{$V$-pairs}, a technical key tool in this paper.
It is a generalization to the modulus world of \emph{standard triples} invented by Voevodsky (see \cite[Lecture 11]{mvw}).
We follow the formulation introduced by \cite{bv}.

In \S\ref{localinj} we prove using the result from \S\ref{Vpair}, 
\emph{the local injectivity result for $F\in \CItsp$},
which implies that for a semi-localization $X$ of an object of $\Sm$ and a dense open immersion $U\hookrightarrow X$, the restriction 
$F(X,\emptyset)\to F(U,\emptyset)$ is injective.

In \S\ref{CohP1} we prove a vanishing theorem of cohomology with coefficient in $F\in \CItspNis$ for a pair of an affine open $X\subset \P^1$ and an effective Cartier divisor $\Sigma\subset X $ (Theorem \ref{thm2-P1}).

In \S\ref{contraction} we introduce the \emph{contractions of $F\in \CItspNis$} as a modulus analogue of Voevodsky's contractions (see \cite[Lecture 23]{mvw}).
It describes cohomology groups with support of $F\in \CItspNis$.

In \S\ref{fibration} we introduce a fibration technique which is used in the proof of another vanishing theorem in \S\ref{vanishing}, which will implies Theorem \ref{thm-intro;purityM}. 

In \S\ref{Gysin} we introduce Gysin maps for cohomology of $F\in \CItspNis$,
which are used in the proof of the vanishing theorem in \S\ref{vanishing}.


In \S\ref{cubeinvcoh} we prove a vanishing theorem of cohomology of $\P^1$ and
Theorem \ref{thm2-intro;purityM} using Theorem \ref{thm-intro;purityM}. 
A basic idea is taken from \cite[Lecture 24]{mvw}.

In \S\ref{cubeNis} we prove that the sheafication preserves 
$\CItsp$ proving Theorem \ref{thm-intro;localinjectivity}.

In \S\ref{RSCproof} we deduce Theorems \ref{thm0;rec} and \ref{thm;purityrec}
from Theorems \ref{thm-intro;localinjectivity} and \ref{thm-intro;purityM}
respectively. We also deduce Theorem \ref{thm;Voev} from Theorems \ref{thm-intro;purityM} and \ref{thm2-intro;purityM}.
 
\subsection*{Acknowledgements} 
The author is deeply indebted to Bruno Kahn, H. Miyazaki and Takao Yamazaki for the collaborations \cite{kmsy} and \cite{kmsy2} from which basic ideas in this paper originate. He is grateful to the referee for very careful reading and for numerous useful comments and suggestions, and also to Ofer Gabber for his kindly answering questions related to this paper.
Finally he is very grateful to J. Ayoub for his pointing out a fundamental mistake in an earlier version of \cite{kmsy}. The author believes that the whole theory has become deepened by the effort to fixing it, 
Part of this work was done while the author stayed at the university of Regensburg supported by the SFB grant ``Higher Invariants". 
Another part was done in a Research in trio in CIRM, Luminy. 
He is grateful to the support and hospitality received in all places.

\subsection*{Notation and conventions}
In the whole paper we fix a base field $k$. 
Let $\Sm$ be the category of separated smooth schemes of finite type over $k$,
and let $\Sch$ be the category of separated schemes of finite type over $k$.

We call a morphism of schemes $X\to Y$ \emph{essentially \'etale}
(resp. \emph{essentially smooth }) if one can write $X$ as a limit 
$X= \lim_{i \in I} X_i$ over a filtered set $I$ 
where $X_i $ is \'etale (resp. smooth) over $Y$ and all transition maps are \'etale and affine. Let $\tSm$ be the category of $k$-schemes $X$ which are essentially smooth over $k$.
We frequently allow $F\in \PST$ to take values 
on objects of $\widetilde{\Sm}$ by
$F(X) := \colim_{i \in I} F(X_i)$.

For a closed immersion $Z\hookrightarrow X$ of affine schemes,
$X^{h}_{|Z}$ denote the henselization of $X$ along $Z$.

\medbreak


\section{$\bcube$-invariance and reciprocity}\label{Mrec}

In this section we review basic notions and facts from \cite{kmsy} and \cite{kmsy2}.

\subsection{Modulus pairs}

\begin{definition}\label{def;moduluspair}(see \cite[Def.1.1.1 and 1.3.1]{kmsy})
\begin{itemize}
\item[(1)]
A \emph{modulus pair} is a pair $\fX=(\Xb,\Xinf)$ where $\Xb \in \Sch$ and $\Xinf$ is an effective Cartier divisor on $\Xb$ such that  
$X=\Xb - |\Xinf|\in \Sm$ and is dense in $\Xb$ (The case $|\Xinf|=\emptyset$ is allowed). We call $X$ the interior of $\fX$.
We say that $\fX$ is \emph{proper} if $\Xb$ is proper over $k$.
\item[(2)]
Let $\fX=(\Xb,\Xinf)$ and $\fX'=(\Xb',\Xinf')$ be modulus pairs with
$X=\Xb - |\Xinf|$ and $X'=\Xb' - |\Xinf'|$.
Let $Z \in \Cor(X', X)$ be an elementary correspondence and $\bar Z^N$ be the normalization of the closure of $Z$ in $\Xb' \times \Xb$ with $p : \bar Z^N \to \Xb$ and $q : \bar Z^N \to \Xb'$ the induced morphisms.
We say $Z$ is \emph{admissible} for $(\fX', \fX)$ 
if $q^* \Xinf' \geq p^* \Xinf$ (an inequality of Cartier divisors).
We say $Z$ is \emph{left-proper} (resp. \emph{finite}) for $(\fX', \fX)$ 
if $\Zb$ is proper (resp. finite) over $\Xb'$.
An element of $\Cor(X', X)$ is called admissible (resp. left-proper, resp. finite)
if all of its irreducible components are admissible (resp. left-proper, resp finite).
\item[(3)]
Let $\ulMCor$ (resp. $\ulMCorfin$) be the additive category of modulus pairs and 
left-proper (resp. finite) admissible correspondences (see \cite[Pr.1.2.3]{kmsy}). 
By definition $\ulMCorfin$ is a full subcategory of $\ulMCor$.
Let $\MCor$ be the full subcategory of $\ulMCor$ whose objects are proper modulus pairs. We let 
\[ b :\ulMCorfin\to \ulMCor \qaq \tau: \MCor \to \ulMCor\]
denote the inclusion functors.
\item[(4)]
For $\fX=(\Xb,\Xinf),\;\fX'=(\Xb',\Xinf')\in \ulMCor$, we put
\[\fX\otimes\fX' = (\Xb\times,\Xb',\Xb\times\Xinf' + \Xinf\times\Xb')\in \ulMCor.\]
\item[(5)]
For $n\in \Z_{>0}$, put $\fX^{(n)}=(\Xb,n\Xinf)$, where $n\Xinf\hookrightarrow \Xb$ is the $n$-th thickening of $\Xinf\hookrightarrow \Xb$.
\item[(6)]
Put $\bcube=(\P^1,\infty)\in \MCor$.
\end{itemize}
\end{definition}

\begin{definition}\label{def;ulMPST}
Let $\ulMPST$ (resp. $\ulMPSTfin$, resp. $\MPST$) be the abelian category of additive presheaves of abelian groups on $\ulMCor$ (resp. $\ulMCorfin$, resp. $\MCor$).
For $\fX\in \ulMCor$ (resp. $\fX\in \MCor$) let $\Ztr(\fX)\in \ulMPST$ 
(resp. $\Ztr(\fX)\in \MPST$) be the object represented by $\fX$.
\end{definition}

\begin{rk}\label{rk;moduluspair}
We will use a limit $\fX=\lim_{i \in I} \fX_i$ over a filtered set $I$, where
$\fX_i=(\Xb_i,X_{i,\infty})\in \ulMCor$ and transition maps are \'etale maps 
$\Xb_i \to \Xb_j$ such that $X_{i,\infty}=X_{j,\infty}\times_{\Xb_i}\Xb_j$.
Those $\fX$ form a larger category $\widetilde{\ulMCor}$ containing $\ulMCor$ and  
the inclusion $\ulMCor\hookrightarrow \widetilde{\ulMCor}$ induces
an equivalence of $\ulMPST$ with the category of additive functors $\widetilde{\ulMCor}^{\op}\to \Ab$ which commute with colimits.
By abuse of notation we will write $\ulMCor$ for $\widetilde{\ulMCor}$.
\end{rk}

By \cite[Pr. 2.2.1, 2.3.1 and 2.4.1]{kmsy} there are pairs of adjoint functors 
\begin{equation}\label{eq;adjunction}
\MPST\begin{smallmatrix}\tau^*\\ \longleftarrow\\ \tau_!\\ \longrightarrow\\ 
\end{smallmatrix}\ulMPST, \;\; 
\ulMPST\begin{smallmatrix} \ulomega^*\\ \longleftarrow\\ \ulomega_!\\ \longrightarrow\\
\end{smallmatrix}\PST, \;\;
\MPST\begin{smallmatrix} \omega^*\\ \longleftarrow\\ \omega_!\\ \longrightarrow\\
\end{smallmatrix}\PST. \;\;
\end{equation}
Here $\tau^*$ is the restriction along $\tau$ from Definition \ref{def;moduluspair}(3)
and $\tau_!$ is its left Kan extension, and $\omega^*$ (resp. $\ulomega^*$) is induced by the functor $\omega:\MCor \to \PST$ (resp. $\ulomega:\ulMCor \to \PST$) given by $(\Xb,\Xinf) \to \Xb-|\Xinf|$, and $\omega_!$ (resp. $\ulomega_!$) is its left Kan extension.

For $X\in \Sm$ let $\MP(X)$ be the category of objects $\fX=(\Xb,\Xinf)\in \MCor$ such that $\Xb-|\Xinf|=X$. Given $\fX_1, \fX_2 \in \MP(X)$,
define $\MP(X)(\fX_1, \fX_2)$ to be $\{ 1_X \}$
if $1_X$ is admissible for $(\fX_1, \fX_2)$, and $\emptyset$ otherwise.


\begin{lemma}\label{lem;MPST}(\cite[\S2]{kmsy})
\begin{itemize}
\item[(1)]
All functors in \eqref{eq;adjunction} are exact, and $\ulomega^*$, $\omega^*$ and $\tau_!$ are fully faithful, and $\ulomega_!\tau_!=\omega_!$ and $\tau_!\omega^*=\ulomega^*$.
\item[(2)]
$\ulomega_!F(X)=F(X,\emptyset)$ for $F\in \ulMPST$ and $X\in \Sm$.
\item[(3)]
$\omega_!\Ztr(\fX)=\Ztr(X)$ for $F\in \MPST$, $X\in \Sm$, $\fX\in \MSm(X)$.

\item[(4)]
For $F\in \MPST$ and $X\in \Sm$, 
\begin{equation*}\label{rem:single-cpt}
\omega_!(F)(X)\simeq \colim_{\fX\in \MSm(X)} F(\fX) \simeq \colim_{n>0} F(\fY^{(n)}),
\end{equation*}
where $\fY$ is any fixed object of $\MSm(X)$.
\item[(5)]
For $F\in \MPST$ and $\fX=(\Xb,\Xinf)\in \ulMCor$, 
\[ \tau_!F(\fX) \simeq \underset{\fY\in \Comp(\fX)}{\colim} F(\fY),\]
where $\Comp(\fX)$ defined in \cite[Def. 1.8.1]{kmsy}, is the category of pairs 
$(\fY,j)$ of $\fY=(\Yb,\Yinf)\in\MCor$ and a dense open immersion $j:\Xb\hookrightarrow \Yb$ such that 
$\Yinf= \Xinf' +\Sigma$ for effective Cartier divisors $\Xinf',\Sigma$ on $\Yb$ such that $|\Sigma|=\Yb-\Xb$ and $j^*\Xinf'=\Xinf$.
For $\fY=(\Yb,\Yinf)\in\Comp(\fX)$ we have $\Xb-|\Xinf|=\Yb-|\Yinf|$ and $\fY$ is equipped with $j_\fY\in \ulMCor(\fX,\fY)$ which is the identity on $\Xb-|\Xinf|$.
For $\fY_1,\fY_2\in \Comp(\fX)$ put
\[\CompM(\fY_1,\fY_2)=\{\gamma\in \MCor(\fY_1,\fY_2)\;|\; \gamma\circ j_{\fY_1}  =j_{\fY_2}\}.\]
\item[(6)]
The restriction functor $b^*: \ulMPST \to \ulMPSTfin$ along $b$ from Definition \ref{def;moduluspair}(2) is fully faithful and admits a left adjoint $b_!$ given by 
\[ b_! F(\fX) =\colim_{\fY\in \Sigmafin\downarrow \sX} F(\fY)
\;\;(F\in \ulMPSTfin, \;\fX\in \ulMCor),\]
where $\Sigmafin$ is the subcategory of $\ulMCorfin$ which have the same objects as $\ulMCorfin$ such that the morphisms $f: \fX \to \fY$ with $\fX=(\Xb,\Xinf)$ and $\fY=(\Yb,\Yinf)$ are the graphs of proper morphisms $f:\Xb\to  \Yb$ which induce isomorphisms on the interiors and satisfy $\Xinf=f^*\Yinf$.
\end{itemize}
\end{lemma}

\begin{remark}\label{rem;m-rec}
Let $\fX=(\Xb,\Xinf)\in \ulMCor$. By \cite[Lem. 1.8.2]{kmsy}, $\Comp(\fX)$ is nonempty. For $\fY=(\Yb, \Xinf' +\Sigma)\in\Comp(\fX)$ as in Lemma \ref{lem;MPST}(5), we have
\[\fY_m:=(\Yb,\Xinf' + m\Sigma)\in \Comp(\fX)\;\text{ for all } m\in \Z_{>0}\]
and we have a natural isomorphism
\begin{equation}\label{eq;rem;m-rec}
\underset{m\in \Z_{>0}}{\colim} F(\fY_m) \isom
\underset{\fY\in \Comp(\fX)}{\colim} F(\fY).
\end{equation}
\end{remark}

\begin{definition}\label{def:Nissheaves}
For $F\in \ulMPST$ or $F\in \ulMPSTfin$ and for $\fX=(\Xb,\Xinf)\in \ulMCor$ write $F_{\fX}$ for the presheaf on $\Xb_{\et}$ given by $U\to F(\fX_U)$ for $U\to \Xb$ \'etale, where $\fX_U=(U,\Xinf\times_{\Xb}U)\in \ulMCor$. We say \emph{$F$ is a Nisnevich sheaf} if $F_{\fX}$ is a sheaf on $\Xb_{\Nis}$ for all $\fX\in \ulMCor$.
We write $\ulMNST\subset \ulMPST$ and $\ulMNSTfin\subset \ulMPSTfin$ for the full subcategories of Nisnevich sheaves.
Let $\MNST\subset \MPST$ be the full subcategory of those objects $F$ such that
$\tau_!F\in \ulMNST$. 
\end{definition}

We have the following (cf. \cite[Th.3.14 and Pr.3.1.8]{voetri}).

\begin{thm}\label{thm;ulMNST}(\cite[Th. 3.5.3, 4.5.5 and Th. 2]{kmsy} and \cite[Th. 4.2.4 and Th. 2]{kmsy2})
\begin{itemize}
\item[(1)]
The inclusion $\ulMNSTfin\to \ulMPSTfin$ has an exact left adjoint $\ulaNisfin$
such that $(\ulaNisfin F)_\fX$ is the Nisnevich sheafication of $F_\fX$ for 
every $F\in \ulMPSTfin$ and $\fX\in \ulMCorfin$. 
\item[(2)]
The inclusion $\ulMNST\to \ulMPST$ has an exact left adjoint $\ulaNis$
such that 
\begin{equation}\label{eq;ulaNisformular}
(\ulaNis F)(\fX) =\colim_{\fY\in \Sigmafin\downarrow \fX} \ulaNisfin b^*F(\fY)
\end{equation}
for every $F\in \ulMPST$ and $\fX\in \ulMCor$. 
\item[(3)]
The inclusion $\MNST\to \MPST$ has an exact left adjoint $\aNis$
such that $\tau_!\circ \aNis = \ulaNis\circ \tau_!$.
\end{itemize}
\end{thm}

\begin{definition}\label{def:Nissheafication}
For $F\in \ulMPST$ we write $F_{\Nis}=\ulaNisfin b^* F \in \ulMNSTfin$.
\end{definition}

\begin{remark}\label{rem;Nissheafication}
For $F\in \ulMPST$ and $\fX=(\Xb,\Xinf)\in \ulMCor$ such that $\Xb$ is regular of dimension $1$, we have $\ulaNis F(\fX) = F_{\Nis}(\fX)$.
This follows from \eqref{eq;ulaNisformular} noting that $\Sigmafin\downarrow \fX$ has the unique object $\fX$ by the assumption.
\end{remark}

\medbreak

\begin{lemma}\label{lem;ulMNSTomega}(\cite[Pr. 6.2.1]{kmsy2})
The functors $\ulomega_!$ and $\omega^*$ from \eqref{eq;adjunction} induce 
\[\ulMNST \to \NST \qaq \NST \to \ulMNST\]
which are denoted again by $\ulomega_!$ (resp. $\ulomega^*$) respectively.\footnote{
In \cite{kmsy} these functors are denoted by $\ulomega_{\Nis}$ and $\ulomega^{\Nis}$ respectively.} 
We have
\begin{equation}\label{eq1;omegaNis}
\aVNis\circ \ulomega_! = \ulomega_!\circ \ulaNis,
\end{equation}
where $\aVNis$ is an exact left adjoint to the natural inclusion $\NST\to \PST$ 
constructed by Voevodsky (\cite[Th.3.14]{voetri}).
\end{lemma}

By Lemma \ref{lem;ulMNSTomega}, for $F\in \ulMNST$ we have an isomorphism of sheaves on $X_{\Nis}$:
\begin{equation}\label{eq2;omegaNis}
 (\ulomega_!F)_X \simeq F_{(X,\emptyset)}\qfor X\in \Sm .
\end{equation}

\subsection{$\bcube$-invariance}

\begin{definition}\label{def:ulMCorls}
Let $\ulMCorls\subset \ulMCor$ be the full subcategory of $(\Xb,\Xinf)\in \ulMCor$ with $\Xb\in \Sm$ and $|\Xinf|$ a simple normal crossing divisor on $\Xb$. 
Note that $\ulMCorls$ is stable under $\otimes$ on $\ulMCor$.
\end{definition}

\begin{definition}\label{def:CI}
$F\in \ulMPST$ (resp. $F\in \MPST$) is $\bcube$-invariant at $\fX\in \ulMCor$ (resp. $\fX\in \MCor$) if the map $pr^*: F(\fX)\to F(\fX\otimes\bcube)$ induced by the projection $pr:\fX\otimes\bcube\to\fX$ is an isomorphism. 
$F$ is $\bcube$-invariant (resp. $\lscube$-invariant) if it is so at 
all $\fX\in \ulMCor$ (resp. $\fX\in \ulMCorls$). 
\end{definition}

\begin{remark}\label{def-rem:CI}
Since the zero section $i_0:\Spec(k) \to \bcube$ is the right inverse of $pr$,
$pr^*: F(\fX)\to F(\fX\otimes\bcube)$ is an isomorphism if and only if
$i_0^*: F(\fX\otimes\bcube)\to F(\fX)$ is injective.
\end{remark}

\begin{rk}\label{rk;def:semipure}
Assume the following condition which holds if $\ch(k)=0$:
For any pair $(X,D)$ of a proper scheme $X$ over $k$ and an effective Cartier divisor $D$ on $X$ such that $X - |D|\in \Sm$ and is dense in $X$, there exists a proper birational map $\pi:X'\to X$ such that $X'\in \Sm$ and the support of $D'=\pi^{-1}(D)$ is a simple normal crossing divisor and that $\pi$ is an isomorphism over $X-|D|$. Then the induced map $(X',D') \to (X,D)$ is an isomorphism in $\ulMCor$ by \cite[Pr.1.9.2 c)]{kmsy}. This implies that $F\in \ulMPST$ or $F\in \MPST$ is $\bcube$-invariant if and only if $F$ is $\lscube$-invariant. 
\end{rk}

\begin{lemma}\label{lem;cubeinv} 
The $\bcube$-invariance and $\lscube$-invariance are preserved by
subobjects and quotients.
\end{lemma}
\begin{proof}
The preservance for subobjects follows from Remark \ref{def-rem:CI}.
That for quotients follows then by the five lemma.
\end{proof}

\begin{lemma}\label{lem;CItau}
$F\in \MPST$ is $\bcube$-invariant if and only if $\tau_!F\in \ulMPST$ is $\bcube$-invariant.
\end{lemma}
\begin{proof}
The if-part follows from the fact $\tau_!F(\fX)=F(\fX)$ for $\fX\in \MCor$ and
the only-if-part follows from Lemma \ref{lem;MPST}(5).
\end{proof}

\begin{lemma}\label{lem;omegaHICI}
For $F\in \HI$, $\omega^*F\in \MPST$ is $\bcube$-invariant.
\end{lemma}
\begin{proof}
This follows from the fact that 
for $\fX=(\Xb,\Xinf)\in \MCor$ with $X=\Xb-\Xinf$, we have
$\omega^*F(\fX)=F(X)$ and $\omega^*F(\fX\otimes\bcube)=F(X\times\A^1)$. 
\end{proof}

\begin{definition}\label{def;hcubeM}
For $F\in \ulMPST$, define $\hM F\in \ulMPST$ by:
\begin{equation}\label{eq;hM}
 \hM F(\fX) =\Coker\big(F(\fX\otimes\bcube)\rmapo{i_0^*-i_1^*} F(\fX)\big)
\;\;(\fX\in \ulMCor)\\
\end{equation} 
where $i_\ep^*$ for $\ep=0,1$ is the pullback by the section $i_\ep: \Spec(k) \to \bcube$. 
For $F\in \MPST$, define $\hM F\in \MPST$ in the same way except the cokernel taken in $\MPST$. By Lemma \ref{lem;MPST}(5), we have
\begin{equation}\label{eq;def;hcubeM} 
\tau_! \hM F =\hM {\tau_!F}\qfor F\in \MPST.
\end{equation}
For $\fX\in \ulMCor$ or $\fX\in \MCor$, we write $\hM \fX=\hM{\Ztr(\fX)}$.
\end{definition}

\begin{defn}\label{def;CI}
Let $\CI\subset \MPST$ be the full subcategory of $\bcube$-invariant objects.
\end{defn}

\begin{lemma}\label{lem0;hcube}
For $F\in \MPST$, the following conditions are equivalent.
\begin{itemize}
\item[(i)]
$F\in \CI$.
\item[(ii)]
The natural map $F(\fX)\to \hM F(\fX)$ is an isomorphism. 
\end{itemize}
\end{lemma}
\begin{proof}
Assume (i) and take $\fX\in \MCor$.
By Remark \ref{def-rem:CI} the assumption implies that $i_\epsilon^*: F(\fX\otimes\bcube) \to F(\fX)$ for $\epsilon=0,1$ are both the inverse of $pr^*: F(\fX)\isom F(\fX\otimes\bcube)$ so that $i_0^*-i_1^*=0$, which implies (ii).

Assume (ii). By \eqref{eq;hM} this implies that for any $\fX\in \MCor$ we have
\begin{equation}\label{eq1;lem0;hcube}
 i_0^*=i_1^* : F(\fX\otimes\bcube)\to F(\fX).
\end{equation}
We consider the multiplication map
\[ \mu : \A^1 \times \A^1 \to \A^1; \quad (x, y) \mapsto (xy), \]
Let $\Gamma\subset \A^1\times\A^1\times\A^1$ be the graph of $\mu$.

\begin{claim}\label{claim;lem0;hcube}
We have $\Gamma \in \MCor(\bcube\otimes \bcube, \bcube)$.
\end{claim}

Indeed consider  
\[\pi : P:=\Bl_S(\P^1 \times \P^1) \to \P^1 \times \P^1 \qwith
S:=\{ 0 \times \infty , \infty \times 0 \}.\]
Then one easily checks that 
$\mu$ extends to a morphism $\tilde{\mu} : P \to \P^1$ and that 
\[\pi \times \tilde{\mu} : P \to (\P^1 \times \P^1) \times \P^1\]
is a closed immersion whose image is precisely the closure of $\Gamma$
 in $\P^1 \times \P^1 \times \P^1$.
Now the pull-back of
$(\P^1 \times \P^1) \times \infty$
to $P$ is the strict transform $T$ of 
$\P \times \infty + \infty \times \P^1$,
while that of $(\P^1 \times \infty + \infty \times \P^1) \times \P^1$
is the sum of $T$ and the exceptional divisors. This proves the claim.

By the claim we have a commutative diagram
\[\xymatrix{
F(\fX \otimes \bcube) \ar[r]^{(id_\fX\otimes i_0)^*} \ar[d]_{(id_\fX \otimes \Gamma)^*} & F(\fX) \ar[d]^{pr^*} \\
F(\fX \otimes \bcube\otimes\bcube) \ar[r]^{(id_{\fX\otimes \bcube} \otimes i_0)^*}  &F(\fX \otimes \bcube).
}\]
By the diagram and \eqref{eq1;lem0;hcube}, we get
\[
pr^* (id_\fX\otimes i_0)^* 
= (id_{\fX\otimes \bcube} \otimes i_0)^* \circ (id_\fX \otimes \Gamma)^* 
= (id_{\fX\otimes \bcube} \otimes i_1)^* \circ (id_\fX \otimes \Gamma)^* 
= \id_{\fX \otimes \bcube}^*.
\]
This proves the surjectivity of $pr^*$, which completes the proof of the lemma.
\end{proof}

\begin{lemma}\label{lem0.1;hcube}
For $F\in \MPST$, we have $\hM F\in \CI$.
\end{lemma}
\begin{proof}
By the definition of $\hM F$, for any $\fX\in \MCor$, the map
\[ \hM F(\fX\otimes\bcube) \rmapo{i_0^*-i_1^*} \hM F(\fX) \]
is the zero map so that 
$\hM F(\fX)\simeq \hM{\hM F}(\fX)$.
Hence the lemma follows from Lemma \ref{lem0;hcube}.
\end{proof}

\begin{lemma}\label{lem1;h0cube}
For $F\in \CI$ and $a\in F(\fX)$ with $\fX\in \MCor$, the corresponding map 
$\ta:\Ztr(\fX) \to F$ in $\MPST$ factors through $\hM \fX$.
\end{lemma}
\begin{proof}
Noting that $F \to \hM F$ gives an end-functor on $\MPST$,
$\ta$ induces $\hM \fX \to \hM F$. Hence the assertion follows from Lemma \ref{lem0;hcube}.
\end{proof}

\begin{definition}\label{def;hMMcube}
For $F\in \MPST$, define $\hMM F\in \MPST$ by
\begin{equation}\label{eq;hMM}
\hMM F(\fX) = \Hom_{\MPST}(\hM {\fX},F)\;\; (\fX\in \MCor).
\end{equation}
\end{definition}

\begin{lemma}\label{lem;hMM}
For $F\in \MPST$, $\hMM F$ is the maximal $\bcube$-invariant subobject of $F$. 
The induced functor 
\[h^0_{\bcube} :\MPST \to \CI;\;  F\to \hMM F\]
gives a right adjoint of the inclusion $\CI \inj \MPST$. 
\end{lemma}
\begin{proof}
The fact that $\hMM F$ is a subobject of $F$ follows from \eqref{eq;hMM} and the fact that $\hM \fX$ is a quotient of $\Ztr(\fX)$.
The fact $\hMM F\in \CI$ follows from Lemmas \ref{lem;cubeinv} and \ref{lem0.1;hcube}. Now let $G\subset F$ be a subobject which is in $\CI$.
For $a\in F(\fX)$ with $\fX\in \MCor$, let $\ta:\Ztr(\fX)\to F$ be the corresponding map in $\MPST$ . If $a\in G(\fX)$, $\ta$ factors through $G$ and hence
factors through $\hM \fX$ by Lemma \ref{lem1;h0cube}. Hence $a\in \hMM F(\fX)$ by \eqref{eq;hMM}. This proves $G\subset \hMM F$, which completes the proof of the first assertion. The second assertion follows easily from the first and Lemma \ref{lem;cubeinv}.
\end{proof}

\subsection{$M$-reciprocity and semipurity}

\begin{definition}\label{def:m-rec}
Let $F \in \ulMPST$.
We say $F$ has \emph{$M$-reciprocity} 
if the following equivalent conditions holds:
\begin{enumerate}
\item
the counit map
$\tau_! \tau^*F \to F$ is an isomorphism.
\item
$\tau_! G \cong F$ for some $G \in \MPST$.
\item
For $\fX \in \ulMCor$ we have
$\underset{\fY\in \Comp(\fX)}{\colim} F(\fY) \simeq F(\fX)$.
\end{enumerate}
The equivalence follows easily from Lemma Lemma \ref{lem;MPST}(1) and (5).
\end{definition}

\medbreak

\begin{lemma}\label{lem1;Mrec}
Let $F \in \ulMPST$ have $M$-reciprocity. 
\begin{itemize}
\item[(1)]
Let $\fX=(\Xb,\Xinf)\in \ulMCor$ with $X=\Xb-|\Xinf|$.
Assume $\Xinf=\Xinf^+ + \Xinf^-$, the sum of two effective Cartier divisors.
Put $\Xb^+=\Xb-|\Xinf^-|$ and $\fX^+=(\Xb^+,\Xb^+\cap \Xinf^+)$.
Then the natural map 
\[\underset{n \in \Z_{>0}}{\colim} F(\Xb,\Xinf^+ + n \Xinf^-) \to F(\fX^+)\]
is an isomorphism.
\item[(2)]
For $\fX\in \ulMCor$, $\uHom_{\ulMPST}(\Ztr(\fX),F)$ has $M$-reciprocity, where
$\uHom_{\ulMPST}$ denotes the internal hom in $\ulMPST$.
\item[(3)]
$\ulaNis F$ has $M$-reciprocity.
\end{itemize}
\end{lemma}
\begin{proof}
We first prove (1). By Remark \ref{rem;m-rec}, there is 
$(\Yb,\tXinf^+ +\tXinf^- +\Sigma)\in \Comp(\fX)$ such that 
$\Xb=\Yb-|\Sigma|$ and $\Xinf^\pm = \tXinf^{\pm}\cap \Xb$.
Then
\[(\Yb,\tXinf^+ + n \tXinf^- + m\Sigma)\in \Comp(\fX^+)\;\;\text{for all }n,m\in \Z_{>0},\]
\[(\Yb, \tXinf^+ + n \tXinf^- + m\Sigma)\in \Comp((\Xb,\Xinf^+ +n\Xinf^-))\;\;\text{for all } m\in \Z_{>0}.\]
Hence \eqref{eq;rem;m-rec} implies
\[F(\fX^+)\simeq \underset{m,n\in \Z_{>0}}{\colim} F(\Yb,\tXinf^+ + n\tXinf^- + m\Sigma))
\simeq \underset{n\in \Z_{>0}}{\colim} F(\Xb,\Xinf^+ + n \Xinf^-),\]
which proves (1). 

For $\fY\in \ulMCor$, we have isomorphisms
\begin{multline*}
\tau_!\tau^*\uHom_{\ulMPST}(\Ztr(\fX),F)(\fY) \simeq
\underset{\fZ\in \Comp(\fY)}{\colim}\uHom_{\ulMPST}(\Ztr(\fX),F)(\fZ) 
\\=\underset{\fZ\in \Comp(\fY)}{\colim} F(\fX\otimes\fZ)
\simeq F(\fX\otimes\fY) = \uHom_{\ulMPST}(\Ztr(\fX),F)(\fY), 
\end{multline*}
where the first (resp. second) isomorphism follows from Lemma \ref{lem;MPST}(5)
(resp. (2)). This proves (2).

Finally (3) follows from Theorem \ref{thm;ulMNST}(3) and Definition \ref{def:m-rec}(2).
\end{proof}


\begin{defn}\label{def:semipure}
We say $F\in \ulMPST$ or $\ulMPSTfin$ is \emph{semipure at $\fX\in \ulMCor$} if 
the unit map $F(\fX) \to \ulomega^*\ulomega_! F(\fX)$ is injective.
We say $F\in \ulMPST$ or $\ulMPSTfin$ is \emph{semipure}
if $F\in \ulMPST$ is semipure at any $\fX\in \ulMCor$.
\end{defn}


\begin{lemma}\label{lem1;semipure}
\begin{itemize}
\item[(1)]
For $F\in \ulMPST$ or $\ulMPSTfin$, $F$ is semipure at $(\Xb,\Xinf)\in \ulMCor$ if and only if $F(\Xb,\Xinf) \to F(\Xb-|\Xinf|,\emptyset)$ is injective. 
\item[(2)]
If $F$ is semipure, so is $\uHom_{\ulMPST}(\Ztr(\fX),F)$ for $\fX\in \ulMCor$.
\item[(3)]
If $F\in \ulMPSTfin$ is semipure, so is $\ulaNisfin F$.
\item[(4)]
If $F\in \ulMPST$ is semipure, so is $\ulaNis F$.
\end{itemize}
\end{lemma}
\begin{proof}
(1) is a direct consequence of the formula
\[\ulomega^*\ulomega_! F(\Xb,\Xinf)=F(\Xb-|\Xinf|,\emptyset)\]
from Lemma \ref{lem;MPST}(1). (2) follows directly from (1).
As for (3), we want to show that 
$\ulaNisfin F(\fX) \to \ulaNisfin F(X,\emptyset)$ is injective for $\fX=(\Xb,\Xinf)\in \ulMCor$ with $X=\Xb-|\Xinf|$. By Theorem \ref{thm;ulMNST}(1), this follows from the fact that for any \'etale $\ol{U}\to \Xb$, 
$F(\ol{U},\ol{U}\times_{\Xb} \Xinf) \to F_\Nis(\ol{U}\times_{\Xb}X)$ is injective by the assumption. Finally (4) follows from (3) and \eqref{eq;ulaNisformular}.
\end{proof}

\begin{lemma}\label{lem;Mrecsemipure}
For $F \in \ulMPST$ define $\tilde{F}\in \ulMPST$ as the image of 
$F\to \ulomega^*\ulomega_! F$. Then $\tilde{F}\in \ulMPST$ is semipure.
If $F$ has $M$-reciprocity (resp. is $(ls-)\cube$-invariant), $\tilde{F}$ has $M$-reciprocity (resp. is $(ls-)\cube$-invariant).
\end{lemma}
\begin{proof}
The semipurity of $\tilde{F}$ follows from the exactness of 
$\ulomega^*$ and $\ulomega_!$ (cf. Lemma \ref{lem;MPST}(1)).
The assertion for $\bcube$-invariance follows from Lemma \ref{lem;cubeinv}.
Note $\ulomega^*\ulomega_! \tau_!G\simeq \tau_! \omega^*\omega_! G$ for $G\in \MPST$
by Lemma \ref{lem;MPST}(1). Hence the asssertion for $M$-reciprocity follows from the exactness of $\tau_!$.
\end{proof}

\begin{defn}\label{def;Xi}
\begin{itemize}
\item[(1)]
Let $\CItsp\subset \ulMPST$ be the full subcategory of $\bcube$-invariant objects which are semipure and have $M$-reciprocity 
\item[(2)]
Let $\lsCIt\subset \ulMPST$ be the full subcategory of $\lscube$-invariant objects which have $M$-reciprocity, and $\lsCItsp\subset \lsCIt$ be the full subcategory of \lssemipure objects.
\end{itemize}
\end{defn}

\begin{lemma}\label{lem;def;Xi}
$\lsCIt$ is closed under kernels and cokernels in $\ulMPST$.
\end{lemma}
\begin{proof}
This follows from Lemma \ref{lem;cubeinv} and the fact that 
$\tau_!$ is exact and fully faithful (see Lemma \ref{lem;MPST}(1)). 
\end{proof}


\subsection{Reciprocity presheaves}

\begin{definition}\label{def;hM}
For $\fX\in \MCor$, we define (cf. Definition \ref{def;hcubeM})
\begin{equation*}\label{eq;def;hMp}
\hMw \fX = \omega_! \hM \fX\in \PST.
\end{equation*}
By Lemma \ref{lem;MPST}(1), (3) and (4), $\hMw\fX$ is a quotient of $\Ztr(X)$ 
where $\fX=(\Xb,\Xinf)$ with $X=\Xb-|\Xinf|$, and we have for $Y\in \Sm$ 
\begin{equation}\label{eq;h0sX}
\hMw \fX(Y)=\Coker\big(\MCor((Y,\emptyset)\otimes\bcube,\fX)\rmapo{i_0^*-i_1^*} \Cor(Y,X)\big).
\end{equation} 
\end{definition}

\begin{definition}\label{def;reciprocity}
Let $F\in \PST$ and $X\in \Sm$. We say \emph{$F$ has reciprocity} if 
for any $X\in \Sm$ and $a\in F(X)=\Hom_{\PST}(\Ztr(X),F)$, there exists 
$\fX=(\Xb,\Xinf)\in \MSm(X)$ such that the map $\ta: \Ztr(X) \to F$ associated to 
$a$ factors through $\hMw\fX$. 
We write $\RSC\subset \PST$ for the full subcategory of reciprocity presheaves.
It is easy to see that $\RSC$ is an abelian category closed under subobjects
and quotients in $\PST$.
\end{definition}


\begin{lemma}\label{lem;HIRec}
We have $\HI\subset \RSC$.
\end{lemma}
\begin{proof}
This follows from the fact that $h_0(X)$ is a quotient of $\hMw\fX$ for any $\fX\in \MSm(X)$ in view of \eqref{eq;h0sX} with \eqref{eq;h0X}.
\end{proof}



\begin{lemma}\label{lem;CIRec}
Let $F\in \MPST$. Assume that for any $X\in \Sm$, there exists $\fX\in \MSm(X)$ such 
that $F$ is $\bcube$-invariant at $\fX^{(n)}$ for all $n\in \Z_{>0}$. Then $\omega_!F$ has reciprocity. In particular we have 
$ \omega_!(\CI)\subset \RSC$.
\end{lemma}
\begin{proof}
Take $Y\in \Sm$ and $a\in \omega_!F(Y)$. 
By Lemma \ref{lem;MPST}(4) there is 
$\fY\in \MSm(Y)$ and $\ta\in F(\fY)$ which represents $a$.
It suffices to show the associated map $a:\Ztr(Y)\to \omega_!F$ factor through $h_0(\fY)$. We need show that for any $X\in \Sm$, the induced map
\[ a(X) : \Ztr(Y)(X)=\Cor(X,Y) \to \omega_!F(X)\]
factors through $h_0(\fY)(X)$. Take $\fX\in \MSm(X)$ as in the assumption of 
Lemma \ref{lem;CIRec}. By Lemma \ref{lem;MPST}(3) and (4) we have
\[\Ztr(Y)(X)\simeq \colim_{n>0}\; \Ztr(\fY)(\fX^{(n)}),\quad 
\omega_!F(X) \simeq \colim_{n>0}\; F(\fX^{(n)}),\]
\[\begin{aligned}
h_0(\fY)(X) &\simeq \colim_{n>0}\; \hM{\fY}(\fX^{(n)})\\
&=\colim_{n>0}\; \Coker\big(\Ztr(\fY)(\fX^{(n)}\otimes\bcube) \rmapo{i_0^*-i_1^*} \Ztr(\fY)(\fX^{(n)})\big).\\
\end{aligned}\]
We have a commutative diagram
\[\xymatrix{
\Ztr(\fY)(\fX^{(n)}\otimes\bcube)  \ar[r]^{\hskip 5pti_0^*-i_1^*}\ar[d]_{\ta} &
 \Ztr(\fY)(\fX^{(n)})\ar[d]^{\ta}\\
F(\fX^{(n)}\otimes\bcube) \ar[r]^{i_0^*-i_1^*} & F(\fX^{(n)}) \\
}\]
where the vertical maps are induced by $\ta:\Ztr(\fY) \to F$ associated to
$\ta\in F(\fY)=\Hom_{\MPST}(\Ztr(\fY),F)$. By the assumption the lower $i_0^*-i_1^*$ 
is the zero map so that $\ta: \Ztr(\fY)(\fX^{(n)}) \to F(\fX^{(n)})$ factors through
$\hM \fY ( \fX^{(n)})$.
Taking the colimit over $n\in \Z_{>0}$, this implies the desired assertion.
\end{proof}
\medbreak

Consider the composite functor
\begin{equation}\label{eq;omegaCI}
 \omegaCI: \RSC \rmapo{\omega^*} \MPST \rmapo{h^0_{\bcube}} \CI.
\end{equation}
For $F\in \PST$ we consider the composite map
\begin{equation}\label{eq3;omegaCI}
\omega_!\omegaCI F \to \omega_!\omega^* F \to F,
\end{equation}
where the first (resp. second) map is induced by the inclusion $\hMM F \to F$
(resp. the counit map for the pair $(\omega_!,\omega^*)$
of adjoint functors (cf. \eqref{eq;adjunction})).

\begin{lemma}\label{lem2;CIRec}
For $F\in \MPST$, the following conditions are equivalent.
\begin{itemize}
\item[(1)]
$F\in \RSC$.
\item[(2)]
The map \eqref{eq3;omegaCI} is an isomorphism.
\end{itemize}
In particular we have $\omega_!(\CI)=\RSC$.
\end{lemma}
\begin{proof}
The implication (2)$\Rightarrow$(1) follows from Lemma \ref{lem;CIRec}.
We prove the converse. 
For $X\in \Sm$ we have
\begin{equation*}
\begin{aligned}
\omega_!\omegaCI F(X) &\simeq \colim_{\fX\in \MSm(X)} \Hom_{\MPST}(\hM \fX,\omega^*F)\\
& \simeq \colim_{\fX\in \MSm(X)} \Hom_{\PST}(h_0(\fX),F)\simeq F(X),\\
\end{aligned}
\end{equation*}
where the first (resp. second, resp. last) isomorphism comes from Lemma \ref{lem;MPST}(4) and \eqref{eq;hMM} (resp. \eqref{eq;adjunction} and Definition \ref{def;hM}, resp. Definition \ref{def;reciprocity}). 
This proves the equivalence (1)$\Leftrightarrow$(2).
The last assertion then follows from this and Lemma \ref{lem;CIRec}.
\end{proof}


\section{$V$-pairs}\label{Vpair}

In this section we fix an integral affine $S\in \tSm$.
For $X\in \Sm$ we frequently write $X$ for $(X,\emptyset)\in \ulMCor$.
The following definition is taken from \cite[\S5.1]{bv}.

\begin{defn}\label{def;Vtriple}
A \emph{pre-$V$-pair} over $S$ is a pair $(p:X\to S,Z)$ (or simply denoted by $(X,Z)$) 
where $p$ is a smooth affine morphism of relative dimension one and $Z\subset X$ is
an effective divisor finite over $S$, which satisfy the following conditions:  
\begin{itemize}
\item[$(i)$]
there exists an open immersion $X\hookrightarrow \Xb$ such that $\Xb$ is normal
and $\Xb-X$ is the support of an effective Cartier divisor $\Xinf\subset \Xb$,
and that $p$ extends to a proper morphism $\overline{p}: \Xb\to S$.
\item[$(ii)$]
there exists an affine open $W\subset \Xb$ such that $\Xinf\cup Z\subset W$.
\end{itemize}
The conditions imply that $\Xinf$ is finite over $S$ since it is proper over $S$ and contained in the affine $S$-scheme $W$, and that $\Xinf\cap Z=\emptyset$ since $Z$ is finite over $S$, and that $\overline{p}$ is equidimensional of relative dimension $1$ since so is $p$ and $\Xinf$ is finite over $S$. 
Such $(\Xb,\Xinf)$ is called a \emph{good compactification of $(X,Z)$}. 
A pre-$V$-pair $(p:X\to S,Z)$ over $S$ is called \emph{$V$-pair} if it satisfies
the following condition:
\begin{itemize}
\item[$(iii)$]
The diagonal $Z\hookrightarrow Z\times_S X$ is defined by some
$h\in \Gamma(Z\times_S X,\sO)$.
\end{itemize}
\end{defn}

\begin{rk}\label{rem;def;Vtriple}
\begin{itemize}
\item[(1)]
If $S$ is the spectrum of a field, any given $(X,Z)$ is a pre-$V$-pair, i.e. 
its good compactification always exists.
\item[(2)]
If $(X,Z)$ is a $V$-pair over $S$, then for any morphism $S'\to S$ in $\tSm$
with $S'$ integral affine, the base change $(X,Z)\times_S S'=(X\times_S S',Z\times_S S')$ is a $V$-pair over $S'$.
\item[(3)]
We will encounter a situation (see Lemma \ref{lem3;purity}) where $S$ may not be essentially smooth over $k$ but the conditions $(i)$ and $(ii)$ of Definition \ref{def;Vtriple} are satisfied except that we only require $\Xb\times_S U$ is normal for some dense regular open subset $U\subset S$. Then we call $(X,Z)$ a \emph{quasi-$V$-pair} over $S$.  
\end{itemize}
\end{rk}
\medbreak

In what follows we fix $V$-pairs over $S$ 
\begin{equation}\label{eq;nu}
\nu=(X,Z) \qaq \nu'=(X',Z')
\end{equation}
with a fixed identification $Z=Z'$. 
We also fix a good compactification $(\Xb,\Xinf)$ of $\nu$.

\begin{defn}\label{def;specialfunction}
Let $f$ be a rational function on $X'\times_S X$ and 
$\theta_f=\div_{X' \times_S X}(f)$ be the divisor of $f$ on $X' \times_S X$. 
We call $f$ \emph{admissible for} $(\nu,\nu')$ (or simply \emph{admissible})
if the following conditions are satisfied:
\begin{itemize}
\item[$(1)$]
$f$ is regular in a neighbourhood of $X'\times_S Z$,
\item[$(2)$]
$\theta_f\times_{X} Z =\Delta_Z$, where 
$\Delta_Z: Z\hookrightarrow X'\times_S Z$ is the diagonal,
\item[$(3)$]
$f$ extends to an invertible function on a neighbourhood of $X'\times_S \Xinf$
in $X'\times_S\Xb$.
\end{itemize}
In case $\nu=\nu'$ we call $f$ \emph{strongly admissible} for $\nu$
if $f$ is admissible and $\theta_f$ contains the diagonal 
$\Delta_X: X\hookrightarrow X\times_S X$.
\end{defn}

\begin{rk}\label{rk;specialfunction} For a morphism $\tS\to S$ in $\tSm$ with $\tS$ integral affine, let $\tnu=\nu\times_S\tS$ and $\tnu'=\nu'\times_S\tS$
be as in Remark \ref{rem;def;Vtriple}(2). If $f$ is admissible for $(\nu,\nu')$, then
its pullback $\tilde{f}$ to $(X'\times_S\tS)\times_{\tS}(X\times_S\tS)$ is admissible for $(\tnu,\tnu')$.
\end{rk}

\begin{lemma}\label{lem;specialfunction}
Let $\Xb_\eta$ be the generic fibre of $\Xb/S$. 
The condition $(3)$ of Definition \ref{def;specialfunction}
is equivalent to the following conditions:
\begin{itemize}
\item[$(3')$]
$|\theta_f|$ is finite over $X'$, and 
the support of the divisor of $f$ on $X'\times_S\Xb_\eta$ is contained in $X'\times_SX_\eta$.
\end{itemize}
\end{lemma}
\begin{proof}
First we prove $(3)\Rightarrow (3')$. Clearly $(3)$ implies the second condition of 
$(3')$ and we prove that it implies $|\theta_f|$ is finite over $X'$. 
$(3)$ implies $T=|\div_{X' \times_S \Xb}(f)|$ is contained in $X'\times_S X$ 
so that $T=|\theta_f|$ and it is proper over $X'$.
It is also affine over $X'$ since so is $X'\times_S X\to X'$. 
This implies the desired assertion.

Next assume $(3')$. Since $\Xinf$ is of pure codimension $1$ in $\Xb$,
Definition \ref{def;Vtriple}$(ii)$ implies that any component of $\Xinf$ is finite and surjective over $S$. Hence, by the second condition of $(3')$, $T=|\div_{X' \times_S \Xb}(f)|$ does not contain any component of $X'\times_S \Xinf$. 
Since $|\theta_f|=T\cap (X'\times_S X)$ is finite over $X'$, it is closed in 
$X'\times_S \Xb$. Hence $T\cap (X'\times_S \Xinf)=\emptyset$.
Since $X'\times_S \Xb$ is normal, this implies that $f$ is invertible on 
a neighbourhood of $X'\times_S \Xinf$ in $X'\times_S\Xb$.
This completes the proof of the lemma.
\end{proof}

\begin{defn}\label{def2;specialfunction}
Let $f$ be admissible for $(\nu,\nu')$.
We call $f$ is special
if $f$ satisfies the condition: 
\begin{itemize}
\item[$(\star)$]
$f$ extends to a regular function on a neighbourhood of on $X'\times_S \Xinf$
in $X'\times_S\Xb$ and $f_{|X'\times_S \Xinf}= 1$.
\end{itemize}
\end{defn}

\begin{lemma}\label{lem;specialfunction}
There always exists $f$ admissible and special for $(\nu,\nu')$.
Moreover, if $\nu=\nu'$, one can take $f$ to be strongly admissible. 
\end{lemma}
\begin{proof}
This is shown in \cite[\S5.1]{bv}, which we recall.
Put $T=\Xinf\cup Z$ and let $T\subset W\subset \Xb$ be as in Definition \ref{def;Vtriple}$(ii)$. Since $T$ is finite over $S$, $T \hookrightarrow W$ is a closed immersion of affine schemes and so is $X' \times_S T \to X'\times_S W$. 
Hence one can take $g\in \Gamma(X'\times_S W,\sO)$ such that $g_{|X'\times_S Z}=h'$ with $h'\in \Gamma(X'\times_S Z,\sO)$ as in
Definition \ref{def;Vtriple}$(iii)$ for $(X',Z')=(X',Z)$.
Since $\Xinf\cap Z=\emptyset$, there are $\alpha,\beta\in \Gamma(W,\sO)$ such that
$\alpha_{|\Xinf}=1$, $\alpha_{|Z}=0$, $\beta_{|\Xinf}=0$, $\beta_{|Z}=1$.
Then $f=(\beta g + \alpha)_{|X'\times_S W}\in \Gamma(X'\times_S W,\sO)$ satisfies the desired conditons. 

In case $\nu=\nu'$, let $\overline{\Delta}_{W\cap X}\subset X\times_S W$ be the closure of the diagonal $\Delta_{W\cap X}\subset X\times_S W$.
Since $\Delta_{W\cap X}$ is closed in $X\times_S W$, we have 
\[\overline{\Delta}_{W\cap X}\cap (X\times_S W)=\Delta_{W\cap X}\qaq
\overline{\Delta}_{W\cap X}\cap (X\times_S T)=\Delta_Z.\]
Hence we can take $g$ in the above argument in 
the ideal of $\overline{\Delta}_{W\cap X}\subset X\times_S W$. 
Since $(\overline{\Delta}_{W\cap X}\cup (X\times_S Z)) \cap (X\times_S\Xinf)=\emptyset$, there are 
$\gamma,\delta\in \Gamma(X\times_S W,\sO)$ such that
\[\gamma_{|X\times_S\Xinf}=1,\;\gamma_{|(\overline{\Delta}_{W\cap X}\cup X\times_S Z)}=0,\;\delta_{|X\times_S\Xinf}=0,\;\delta_{|(\overline{\Delta}_{W\cap X}\cup X\times_S Z)}=1.\]
Then $f=(\delta g + \gamma)_{|X\times_S W}\in \Gamma(X\times_S W,\sO)$ satisfies the desired conditons. 
\end{proof}


In what follows we fix an effective Cartier divisor $D\subset S$.
For a scheme $Y$ over $S$, write $D_Y=D\times_S Y$.
 
\begin{defn}\label{def;ulMCorS}
Let $\ulMCor_S\subset \widetilde{\ulMCor}$ (cf. Remark \ref{rk;moduluspair}) be the subcategory of objects of the form
$\fY=(Y,W )$, where $Y$ is an integral scheme equipped with a dominant morphism
$\pi:Y \to S$ and $W\subset Y$ is an effective Cartier divisor whose irreducible
components are all dominant over $S$. For $\fY'=(Y',W')$, define
\begin{equation}\label{eq0;ulMCor}
\ulMCor_S(\fY',\fY)\subset \ulMCor((Y',W'),(Y,W))
\end{equation}
as the subgroup  
generated by those integral cycles on $Y'\times_k Y$ contained in $Y'\times_S Y$. 
\end{defn}

For $\fY=(Y,W )\in \ulMCor_S$, write $\fY_D=(Y, W+D_Y)$.
Note $p_Y^*(D_Y)=p_{Y'}^*(D_{Y'})$,
where $p_Y: Y'\times_S Y \to Y$ and $p_{Y'}: Y'\times_S Y \to Y'$ are the projections.
Hence we have
\begin{equation}\label{eq;ulMCor}
\ulMCor_S(\fY',\fY)\subset \ulMCor(\fY_D,\fY'_D).
\end{equation}

\begin{lemma}\label{lem1;Vtriple}
Let $f$ be admissible for $(\nu,\nu')$. Then
\[\theta_f\in \ulMCor_S((X',\emptyset),(X,\emptyset)),\]
\[\theta_f\in \ulMCor_S((X',Z'),(X, Z)).\]
\end{lemma}
\begin{proof}
By Lemma \ref{lem;specialfunction}, $|\theta_f|$ is finite and surjective over $X'$. The modulus condition follows from \eqref{eq;ulMCor} and Definition \ref{def;specialfunction}(2).
\end{proof}

Let $F\in \ulMPST$ and $f$ be admissible for $(\nu,\nu')$.
By Lemma \ref{lem1;Vtriple} and \eqref{eq;ulMCor}, it induces a map 
\begin{equation}\label{eq;V-triple-thetaf}
 \thetarelf^* : F(X,Z+D_X)/F(X,D_X) \to F(X',Z'+D_{X'})/F(X',D_{X'}).
\end{equation}


\begin{thm}\label{thm2;Vtriple}
Assume $F\in \ulMPST$ has $M$-reciprocity and is $\bcube$-invariant at $X,X',U=X-Z,U'=X'-Z'$
(see Definition \ref{def:CI}).
\begin{itemize}
\item[(1)]
The map $F(X,D_X)\to F(X,Z+D_X)$ is injective.
\item[(2)]
Assume $F$ is semipure at $(X,Z+D_X)$ and $(X',Z'+D_{X'})$ (Definition \ref{def:semipure}). Then $\thetarelf^*$ is independent of $f$ and induces
\[ \phi_{\nu,\nu'} : F(X,Z+D_X)/F(X,D_X) \to F(X',Z'+D_{X'})/F(X',D_{X'})\]
depending only on $(\nu,\nu')$. 
The map $\phi_{\nu,\nu'}$ is an isomorphism and the identity if $\nu=\nu'$. 
For another $V$-pair $\nu''=(X'',Z'')$ over $S$ with an identification $Z''=Z$, 
$\phi_{\nu,\nu''}=\phi_{\nu',\nu''}\circ \phi_{\nu,\nu'}$.
\end{itemize}
\end{thm}

\begin{rk}\label{rk1;thm2;Vtriple}
For $n\in \Z_{>0}$ let $nZ\hookrightarrow X$ be the $n$-th thickening of $Z\hookrightarrow X$. 
Assume that $(X,nZ)$ and $(X',nZ')$ are $V$-pairs over $S$ for all
$n\in \Z_{>0}$. Then the same assertion as Theorem \ref{thm2;Vtriple}(2) holds for 
\[ \phi_{\nu,\nu'} : F(X-Z,D_X)/F(X,D_X) \to F(X'-Z',D_{X'})/F(X',D_{X'})\]
without assuming the semi-purity of $F$.
Indeed we may replace $F$ by $\tilde{F}$ from Lemma \ref{lem;Mrecsemipure} to assume $F$ is semipure. Then the assertion follows from Lemma \ref{lem1;Mrec}(1).
\end{rk}


\begin{defn}\label{def;ulMCorScubeinv}
For $\fX,\fY\in \ulMCor_S$, put
\begin{equation*}\label{eq;}
 \ulMCorScube(\fY,\fX) =\Coker\big(\ulMCor_S(\fY\otimes\bcube,\fX)
\rmapo{i_0^*-i_1^*} \ulMCor_S(\fY,\fX)\big).
\end{equation*} 
\end{defn}

We need preliminary lemmas for the proof of Theorem \ref{thm2;Vtriple}.

\begin{lemma}\label{lem3;Vtriple}
Let $\nu''=(X'',Z'')$ be another $V$-triple over $S$ endowed with 
an identification $Z''=Z$.
Let $f$ (resp. $g$) be  admissible for $(\nu,\nu')$ (resp. $(\nu',\nu'')$).
Then there is $h$  admissible for $(\nu,\nu'')$ such that
$\theta_h= \theta_f\circ\theta_g\in \ulMCorS(X'',X)$.
\end{lemma}
\begin{proof} This is proved in \cite[\S5.1.1]{bv}.
For convenience of the readers, we give an explicit construction of $h$
\footnote{This is suggested by a referee to whom the author is grateful.}.
By definition we can write $\theta_g=\sum_i n_i V_i$ with $n_i\in \Z$, where
$V_i\subset X''\times_S X'$ are closed integral subschemes finite and surjective over $X''$ such that there exists an open $U\subset X''\times_S X'$ containing $X''\times_S \Xinf'$ with $U\cap V_i=\emptyset$ for all $i$, 
and that there exists a unique index $i_0$ such that $n_{i_0}=1$ and 
\[V_{i_0}\cap (X''\times_S Z)=\Delta_Z\;\text{ and }\;
V_{i}\cap (X''\times_S Z)=\emptyset\text{ for } i\not=i_0.\]
Consider the composite map
\[\rho_i:V_i^N\times_S X \to V_i\times_S X\hookrightarrow X''\times_S X'\times_S X \to X'\times_S X,\]
where the first (resp. last) map is induced by the normalization $V_i^N\to V_i$ 
(resp. the projection). Then one can check that
\[ h=\underset{i}{\prod} \Nm_{V_i^N\times_S X/X''\times_S X}(\rho_i^*f)^{n_i}.\]
satisfies the desired condition of the lemma.
\end{proof}

\begin{lemma}\label{lem;MCorScube}
Let $F\in \ulMPST$ be $\bcube$-invariant at $\fY_D$ for $\fY\in \ulMCor_S$ (cf. \eqref{eq;ulMCor}). 
For $\gamma\in \ulMCor_S(\fY,\fX)$, $\gamma^*: F(\fX_D)\to F(\fY_D)$
depends only on the class of $\gamma$ in $\ulMCorScube(\fY,\fX)$
(cf. \eqref{eq0;ulMCor}).
\end{lemma}
\begin{proof} 
It suffices to show $\gamma^*=0: F(\fX_D)\to F(\fY_D)$ if $\gamma=(i_0^*-i_1^*)\Gamma$ for $\Gamma\in  \ulMCor_S(\fY\otimes\bcube,\fX)$. We have a commutative diagram
\[\xymatrix{
F(\fX_D) \ar[r]^{\gamma^*} \ar[rd]_{\Gamma^*} & F(\fY_D) \\
& F(\fY_D\otimes\bcube) \ar[u]_{i_0^*-i_1^*}\\
}\]
By the assumption $F(\fY_D) \overset{pr^*}{\simeq} F(\fY_D\otimes\bcube)$, we have
$i_0^*=i_1^*$, which proves the desired claim.
\end{proof}
\medbreak


Let $f$ be admissible for $(\nu,\nu')$ and let $\theta_f$ denote
the elements in
\[\ulMCor_S(X',(\Xb,\Xinf))\qaq \ulMCor_S((X',Z'),(\Xb,\Xinf+ Z))
\]
which are the composite of $\theta_f$ from Lemma \ref{lem1;Vtriple} and the natural maps
\[X\to (\Xb,\Xinf)\qaq (X,Z)\to (\Xb,\Xinf+ Z)\]
induced by the open immersion $X\hookrightarrow \Xb$.
\medbreak

\begin{lemma}\label{lem2;Vtriple}
Let $f,g$ be  admissible for $(\nu,\nu')$.
There is 
\[\lambda\in \ulMCorS(X',(\Xb,\Xinf+ Z))\]
which makes the triangles in the following diagram commutative in 
$\ulMCorScube(X'-Z',(\Xb,\Xinf+ Z))$ and $\ulMCorScube(X',(\Xb,\Xinf))$:
\[\xymatrix{
X'-Z' \ar[r]^{j'}\ar[d] & X' \ar[dd]^{\theta_f-\theta_g}\ar[ldd]^{\lambda} \\
(X',Z') \ar[d]_{\theta_f-\theta_g} &  \\
(\Xb,\Xinf+ Z) \ar[r]^{\ol{j}} & (\Xb,\Xinf) \\
}\]
\end{lemma}
\begin{proof}
Letting $U'=X'-Z'$, we have a commutative diagram
\[\xymatrix{
\ulMCorScube(U',(\Xb,\Xinf+ Z)) \ar[r]^{\simeq} & 
\Pic(\Xb\times_S U',(\Xinf+ Z)\times_S U')\\
\ulMCorScube(X',(\Xb,\Xinf+ Z)) \ar[r]^{\simeq} \ar[d]^{{\ol{j}}_*}\ar[u]_{{j'}^*}
&  \Pic(\Xb\times_S X',(\Xinf+ Z)\times_S X')\ar[d]\ar[u] \\
\ulMCorScube(X',(\Xb,\Xinf)) \ar[r]^{\simeq} \ar[r] & 
\Pic(\Xb\times_S X',\Xinf\times_S X')\\
}\]
which follows from \cite{ry}.
Consider the line bundle $\sL= (g/f)\sO_{X' \times_S \Xb}\subset \sK$ on 
$X' \times_S \Xb$, where
$\sK$ is the sheaf of meromorphic functions on $X'\times_S\Xb$.
By the conditions $(3)$ and $(2)$ of Definition \ref{def;specialfunction},
$f/g$ is invertible on $X'\times_S (\Xinf+ Z)$.
Thus we have isomorphisms
\[\alpha: \sL_{|(\Xinf+ Z) \times_S X'}\simeq \sO_{(\Xinf+ Z)\times_S X'}\]
given by the multiplication by $f/g$. 
Then the pair $(\sL,\alpha)$ gives an element of 
$\Pic(\Xb\times_S X',(\Xinf+ Z)\times_S X')$
which lifts the classes of $\theta_f-\theta_g$ in 
$\Pic(\Xb\times_S X',\Xinf\times_S X')$ and 
$\Pic(\Xb\times_S U',(\Xinf+ Z)\times_S U')$. This proves Lemma \ref{lem2;Vtriple}.
\end{proof}

\begin{lemma}\label{lemvs;Vtriple}
Assume $\nu=\nu'$ and $f$ is strongly admissible for $\nu$. Then 
\[\theta_f-\Delta_X\in \ulMCorS(X,X- Z)=\ulMCorS(X,(X,Z))\subset \ulMCorS((X,Z),(X,Z)).\]
\end{lemma}
\begin{proof}
By assumption 
$\Gamma:=(\theta_f-\Delta_X)_{|X\times_S W}$ is effective. We have
\[\begin{aligned}
 \Gamma\times_{X\times_S X} (X \times_S Z) 
&= \theta_f\times_{X\times_S X}(X \times_S Z) - 
\Delta_X\times_{X\times_S X}(X \times_S Z) \\
& = \Delta_Z- \Delta_Z =0,
\end{aligned}\]
where the last equality follows from Definition \ref{def;specialfunction}(2).
Hence $|\theta_f-\Delta_X|$ does not meet $X\times_S Z$.
In view of Lemma \ref{lem1;Vtriple} this implies the lemma.
\end{proof}

\begin{lemma}\label{lemH;Vtriple}
Let $f$ be admissible and special for $(\nu,\nu')$. Then $\theta_f=0\in \ulMCorScube(X',(\Xb,\Xinf))$.
\end{lemma}
\begin{proof}
We want to construct 
$\ol{\Theta}_f\in\ulMCorS((X',\emptyset)\otimes\bcube,(\Xb,\Xinf))$ such that
$\partial(\ol{\Theta}_f) =\theta_f$, where
\[\partial=i_0^*-i_1^*: \ulMCorS((X',\emptyset)\otimes\bcube,(\Xb,\Xinf)) \to \ulMCorS((X',\emptyset),(\Xb,\Xinf)).\]
Let $t$ be the coordinate of $\A^1$ and consider the meromorphic function
\[ H_f = tf +1-t \;\text{ on } (X'\times\A^1) \times_S \Xb,\]
and let $\Theta_f=\div_{(X'\times\A^1) \times_S \Xb}(H_f)$.
By Definition \ref{def2;specialfunction},
\[ (H_f)_{|(X'\times \A^1) \times_S \Xinf} = t +1-t=1.\]
Thus $|\Theta_f|\cap \big((X'\times\A^1) \times_S \Xinf\big)=\emptyset$ and hence
\[\Theta_f\in \ulMCor_S(X'\times\A^1,(\Xb,\Xinf)).\]
Let $\ol{\Theta}_f$ be the closure of $\Theta_f$ in $(X'\times\P^1) \times_S \Xb$.
Let $V\subset X'\times_S X$ be an open subset containing $X'\times_S \Xinf$ on which $f$ is regular (cf. Definition \ref{def;specialfunction}(3)).
Putting $u=t^{-1}$, we see 
\[\big(\ol{\Theta}_f\times_{\P^1}(\P^1-0)\big) \cap \big(V\times (\P^1-0)\big) = \div_{V\times(\P^1-0)}(f+u-1).\]
By Definition \ref{def2;specialfunction} this implies  $\ol{\Theta}_f\times_{\P^1}(\P^1-0)\times_{\Xb} \Xinf\subset \{u=0\}$ so that  
$\ol{\Theta}_f\in \ulMCorS((X',\emptyset)\otimes\bcube,(\Xb,\Xinf))$.
It is easy to check $\partial(\ol{\Theta}_f) =\theta_f$ and this completes the proof.
\end{proof}

\begin{lemma}\label{lem4;Vtriple}
Letting $U=X-Z$, there is  
\[\lambda\in \ulMCorS(X,U)=\ulMCorS(X,(X,Z))\]
which makes the following diagram commutative in $\ulMCorScube(X,(\Xb,\Xinf))$:
\[\xymatrix{
U \ar[rd]^{\iota}  & \ar[l]_{\lambda} X \ar[d]^{\tau} \\
 & (\Xb,\Xinf) \\
}\]
where $\iota$ and $\tau$ are induced by the open immersion $X\hookrightarrow \Xb$.
\end{lemma}
\begin{proof}
Consider $\theta_f\in \ulMCorS(X,X)$ for $f$ which is strongly admissible and special (cf. Lemma \ref{lem;specialfunction}). 
By Lemma \ref{lemvs;Vtriple}, $\lambda:= \Delta_X- \theta_f\in \ulMCorS(X,U)$. 
Then 
\[\iota\circ \lambda = \Delta_X - \theta_f=\tau \;\in 
\ulMCorScube(X,(\Xb,\Xinf)),\]
where the last equality follows from Lemma \ref{lemH;Vtriple}.
This completes the proof.
\end{proof}

\begin{cor}\label{cor-lem4;Vtriple}
Assume $F\in \ulMPST$ is $\bcube$-invariant at $X$ with $M$-reciprocity.
Then $j^*: F(X,D_X) \to F(U,D_U)$ is injective, where $j:U=X-Z\to X$ is the open immersion.
\end{cor}
\begin{proof}
For $m\in \Z_{>0}$ let $\tau_m\in \ulMCorS(X,(\Xb,m\Xinf))$ be induced by the open 
immersion $X\hookrightarrow \Xb$.
By Lemma \ref{lem4;Vtriple}, there is 
$\lambda_m\in \ulMCorS(X,U)$ such that the upper triangles in the following diagram
\[\xymatrix{
F(U,D_U) \ar[r]^-{\lambda_m^*}&  F(X,D_X)  \\
F(X,D_X) \ar[u]^{j^*} & \ar[l]^-{\tau_m^*} F(\Xb,m\Xinf+D_{\Xb})\ar[lu]^{\iota^*} 
\ar[u]^{\tau_m^*}\\
}\]
commutes (the commutativity of the lower triangle is obvious). This implies $j^*$ is injective on the image of $\tau_m^*$.
Now the corollary follows from Lemma \ref{lem1;Mrec}(1).
\end{proof}
\medbreak

Let $f$ be admissible for $(\nu,\nu')$.
For $m\in \Z_{>0}$ let
\[ \tthetarelf^* : F(\Xb,m\Xinf + Z+D_{\Xb})/F(\Xb,m\Xinf+D_{\Xb}) 
\to F(X',Z'+D_{X'})/F(X',D_{X'})\]
be the composite of $\thetarelf^*$ from \eqref{eq;V-triple-thetaf} and the map
\begin{equation}\label{eq3;V-triple}
\psi_m^*: F(\Xb,m\Xinf + Z+D_{\Xb})/F(\Xb,m\Xinf+D_{\Xb}) \to F(X,Z+D_X)/F(X,D_X)
\end{equation}
induced by $\psi_m:(X,Z) \to (\Xb,m\Xinf+ Z)$ in $\ulMCor_S$ (cf. \eqref{eq;ulMCor}). 

\begin{lemma}\label{lem5;Vtriple}
Assume $F$ is $\bcube$-invariant at $X'$ and $U'$.
\begin{itemize}
\item[(1)]
If $\nu=\nu'$ and $f$ is strongly admissible for $\nu$, $\thetarelf^*$ is the identity.
\item[(2)]
For $f,g$ admissible for $(\nu,\nu')$, 
$j^*\circ \tthetarelf^* = j^*\circ\tthetarelg^*$, where
\[{j'}^*:F(X',Z'+D_{X'})/F(X',D_{X'}) \to F(U',D_{U'})/F(X',D_{X'})\]
is induced $j': U'=X'-Z'\to (X',Z')$ .
\end{itemize}
\end{lemma}
\begin{proof}
Assume that $f$ is strongly admissible for $\nu$ and consider the commutative diagram
\[\xymatrix{
F(X,D_X)\ar[r]\ar[d]^{\theta_f^*-id} & F(X,Z+D_X) \ar[r]\ar[d]^{\theta_f^*-id}  
& F(X,Z+D_X)/F(X,D_X) \ar[d]^{\thetarelf^*-id} \\
F(X,D_X)\ar[r] & F(X,Z+D_X) \ar[r] & F(X,Z+D_X)/F(X,D_X)\\
}\]
By Lemma \ref{lemvs;Vtriple} the middle vertical map factors through $F(X)$, which implies $\thetarelf^* -id=0$. This proves (1).

To show (2), consider a commutative diagram
\[\xymatrix{
F(\Xb,m\Xinf+D_{\Xb}) \ar[r] \ar[d] & F(\Xb,m\Xinf+ Z+D_{\Xb})  \ar[d] \ar[rd]^{\psi_m^*} &\\
F(X,D_X)\ar[r]\ar[d]^{\theta_f^*-\theta_g^*} 
& F(X,Z+D_X) \ar[r]\ar[d]^{{j'}^*\circ\theta_f^*-{j'}^*\circ\theta_g^*}  
& F(X,Z+D_X)/F(X,D_X) \ar[d]^{{j'}^*\circ\thetarelf^*-{j'}^*\circ\thetarelg^*} \\
F(X',D_{X'})\ar[r] & F(U',D_{U'}) \ar[r] & F(U',D_{U'})/F(X',D_{X'})\\
}\]
By Lemma \ref{lem2;Vtriple} the composite of the middle vertical maps factors 
through $F(X',D_{X'})$, which implies ${j'}^*\circ(\thetarelf^* -\thetarelg^*)\circ  \psi_m^*=0$.
This proves (2).
\end{proof}

\medbreak\noindent
{\it Proof of Theorem \ref{thm2;Vtriple}}:
(1) follows from Corollary \ref{cor-lem4;Vtriple}. We show (2).
Let $f,g$ be  admissible for $(\nu,\nu')$.
The map ${j'}^*$ in Lemma \ref{lem5;Vtriple}(2) is injective by the assumption
on the semi-purity (cf. Lemma \ref{lem1;semipure}(1)).
Hence $\thetarelf^*$ and $\thetarelg^*$ coincide on the image of $\psi_m^*$.
By Lemma \ref{lem1;Mrec}(1) this implies the first assertion.
The last assertion follows from Lemma \ref{lem3;Vtriple}.
It now suffices to show that $\phi_{\nu,\nu}$ is the identity if $\nu=\nu'$.
Choose $g$ strongly admissible for $\nu$.
By Lemma \ref{lem5;Vtriple}(1), $\thetarelg^*$ is the identity on $F(X,Z)/F(X)$.
By the first assertion we have $\phi_{\nu,\nu}=\thetarelg^*$.
This completes the proof.

\medbreak
\begin{cor}\label{cor;thm2;Vtriple}
Let the assumption be as in Theorem \ref{thm2;Vtriple}(2).
Let $\pi:X''\to X'$ be an \'etale morphism which induces an isomorphism $Z'':=X'' \times_{X'} Z' \simeq Z'$.
Assume $\nu''=(X'',Z'')$ is a $V$-pair 
(By Lemma \ref{lem6;Vtriple} below, this is the case if $\nu''$ is a pre-$V$-pair). Then the diagram
\[\xymatrix{
 F(X,Z+ D_X)/F(X,D_X)\ar[r]^{\hskip -10pt\phi_{\nu,\nu'}}\ar[rd]^{\hskip -10pt\phi_{\nu,\nu''}} & 
 F(X',Z'+ D_{X'})/F(X',D_{X'}) \ar[d]^{\pi^*} \\
 & F(X'',Z''+ D_{X''})/F(X'',D_{X''}) \\
 }\]
is commutative and $\pi^*$ is an isomorphism.
\end{cor}
\begin{proof}
Let $f$ be as in Definition \ref{def;specialfunction}. Then the assumption implies that the pullback of $f$ to $X''\times_S X$ is admissible for $(\nu,\nu'')$. This implies the commutativity.
The second assertion follows from this and Theorem \ref{thm2;Vtriple}(2).
\end{proof}

\section{Local injectivity}\label{localinj}

In this section we prove the following.

\begin{thm}\label{thm;localinjectivity}
Let $F\in \lsCIt$ (cf. Definition \ref{def;Xi}).
\begin{itemize}
\item[(1)]
Let $X$ be the semi-localization of an object of $\Sm$, and $V\subset X$ be an open dense subset. Then $F(X)\to F(V)$ is injective.
\item[(2)]
For a dense open immersion $U\subset X$ in $\Sm$,
$F_{\Nis}(X) \to F_{\Nis}(U)$ is injective. 
\end{itemize}
\end{thm}

\begin{remark}\label{rem;thm;localinjectivity}
If $F\in \CIt$, Theorem \ref{thm;localinjectivity} follows from 
\cite[Th. 6]{rec} and \cite[Th. 2]{recII}. 
\end{remark}

\begin{defn}\label{CIrspNisdef}
Let $\ulMNSTls$ be the full subcategory of $\ulMPST$ consisting of such objects
$F$ that for any $\fX=(\Xb,\Xinf) \in \ulMCorls$, $F_\fX$ is a sheaf on $\Xb_\Nis$
(cf. Definition \ref{def:Nissheaves}). 
Put 
\[\lsCItspNis = \lsCItsp\cap \ulMNSTls .\]
\end{defn}

\begin{cor}\label{coro;localinjectivity}
Let $F\in \lsCItspNis$.
For $\fX=(\Xb,\Xinf)\in \ulMCorls$ and a dense open immersion $U\hookrightarrow \Xb$,
the restriction $F(\fX)\to F(\fX_U)$ is injective.
\end{cor}
\begin{proof}
By Lemma \ref{lem1;semipure}(1) we are reduced to the case $\Xinf=\emptyset$, which follows from Theorem \ref{thm;localinjectivity}.
\end{proof}

We need some preliminaries.

\begin{prop}\label{thm;Walker}
Assume $k$ is infinite.
Let $W$ be a smooth affine variety, $D\subset W$ be an effective Cartier divisor 
and $Q\subset D$ be a finite set of points.
Then there exists an affine $S\in \Sm$ and an open neighbourhood $X\subset W$ of $Q$ 
with a smooth morphism $p:X\to S$ such that $(p:X\to S,Z)$ with $Z=X\cap D$ is 
a $V$-pair over $S$.
If $D$ is smooth at $Q$, one can take $Z$ to be \'etale over $S$.
\end{prop}
\begin{proof}
See \cite[Pr.5.3]{bv}.
\end{proof}

\begin{remark}\label{remark0;Walker}
Let $\nu=(X,Z)$ be as above.
\begin{itemize}
\item[(1)]
Let $Q\subset |Z|$ be a finite set of points.
Then $q: X\times_S |Z|\to X$ induced by $Z\to S$ is finite so that $q^{-1}(Q)$ is finite.
Then $\sO_{X\times_S X}(\Delta)$ is trivialized in some open neighbourhood $W$ 
of $q^{-1}(Q)$ in $X\times_S X$.
Since $W\cap (X\times_S |Z|)$ contains $q^{-1}(Q)$ and is open in $X\times_S |Z|$ 
which is finite over $X$, there exists an open neighbourhood $U$ of $Q$ in $X$ such that
$W\cap (X\times_S |Z|)\supset U\times_S |Z|$.
\item[(2)]
Let $nZ\hookrightarrow X$ be the $n$-th thickening of $Z\hookrightarrow X$ for $n\in \Z_{>0}$.
By (1) and the proof of \cite[Pr.5.3]{bv},
one can take $\nu$ in such a way that there exists an open $W\subset X\times_S X$
contaning $X\times_S |Z|$ such that the restriction of 
the diagonal $X\hookrightarrow X \times_S X$ to $W$ is a principal divisor on $W$.
Then $\nu_n=(X,nZ)$ is a $V$-pair over $S$ for any $n\in \Z_{>0}$.
\end{itemize}
\end{remark}

\medbreak\noindent
{\it Proof of Theorem \ref{thm;localinjectivity}}\;\;
By a standard norm argument (see the argument in the last part of \cite[\S7]{rec}),
we may assume $k$ is infinite.
Take $\alpha\in F(X)$ which vanishes in $F(V)$. We need to show $\alpha=0$.
We may assume that $X$ is the localization $W_Q$ of $W\in \Sm$ 
at a finite set $Q$ of points and $V$ is a complement of $D_Q\subset W_Q$ 
for a divisor $D$ on $W$, and that $\alpha$ comes from $\beta\in F(W)$
which vanishes in $F(W-D)$. By Proposition \ref{thm;Walker} we may assume further that
$(W,D)$ is a $V$-pair over some $S\in \Sm$. Then the desired assertion follows from Theorem \ref{thm2;Vtriple}.

To show (2), it suffices to show the injectivity of $F_{\Nis}(X) \to F_{\Nis}(\xi)$, where $\xi$ is the generic point of $X$.
Assume there is a non-zero $f\in F_{\Nis}(X)$ lying in the kernel.
There is a point $x\in X$ such that the image $f_x$ of $f$ in 
$F_{\Nis}(\sO^h_{X,x})=F(\sO^h_{X,x})$ is non-zero, where $\sO^h_{X,x}$ is the 
henselization of $\sO_{X,x}$.
There is a Nisnevich neighbourhood $X'\to X$ of $x$ such that
the image of $f$ in $F_{\Nis}(X')$ comes from $g\in F(X')$.
We have a commutative diagram
\[\xymatrix{
 & \ar[ld] F_{\Nis}(X) \ar[d] \ar[r] & F_{\Nis}(\xi) \ar[d] \\
F_{\Nis}(\sO^h_{X,x})& \ar[l] F_{\Nis}(X') \ar[r] & F_{\Nis}(\xi') \\
F(\sO^h_{X',x})\ar[u]_{\simeq} & \ar[l] F(X') \ar[u] \ar[d] \ar[r]^{\alpha} & F(\xi')\ar[u]_{\simeq} \\
& \ar[lu] F(\sO_{X',x}) \ar[ru]_{\beta}\\
}\]
where $\xi'$ is the generic point of $X'$.
By the assumption on $f$, $\alpha(g)=0$ so that $\beta(g_x)=0$, where
$g_x$ is the image of $g$ in $F(\sO_{X',x})$.
But the assumption $f_x\not=0$ implies $g_x\not=0$ by the above diagram.
This contradicts the injectivity of $\beta$ which follows from (1).
$\square$


\section{Cohomology of $\P^1$ with modulus}\label{CohP1}


\begin{thm}\label{thm-P1}
Let $\eta$ be the generic point of an irreducible object of $\Sm$.
Let $X\subset \P^1_\eta$ be a non-empty affine open subset and $Z\subset X$ be an effective Cartier divisor.
Let $F\in \lsCIt$ be \lssemipure.
\begin{itemize}
\item[(1)]
$F(X,Z)=F_{\Nis}(X,Z)$ (cf. Definition \ref{def:Nissheafication} and Remark \ref{rem;Nissheafication}).
\item[(2)]
$H^i(X_{\Nis},(F_{\Nis})_{(X,Z)})=0$ for $i>0$.
\end{itemize}
\end{thm}
\medbreak

We need some preliminary lemmas. 

\begin{lemma}\label{lemP1-0}
Let $\eta$ be as in Theorem \ref{thm-P1} and 
$X\subset \P^1_\eta$ be non-empty affine open and $Z\subset X$ be an effective Cartier divisor. Put $\Xb=\P^1_\eta$ and $\Xinf=\P^1-X$.
Then $(X,Z)$ is a $V$-pair over $\eta$ and $(\Xb,\Xinf)$ is a good compactification of $(X,Z)$.
\end{lemma}
\begin{proof}
We only check the condition $(iii)$ of Definition \ref{def;Vtriple} for $(X,Z)$.
Other conditions are easily checked. 
Put $K=k(\eta)$. We may assume $X\subset \A_\eta^1$. Then one can write
\[ X=\Spec K[t][1/f(t)],\;\; Z=\Spec K[t]/(g(t))\]
for some $f(t),g(t)\in K[t]$ such that $(f,g)=1$. Then 
\[ X\times_\eta Z = \Spec K[t,s]/(g(t))[1/f(s)]\]
and $Z\hookrightarrow X\times_\eta Z$ is the divisor of $t-s$.
This proves the desired assertion.
\end{proof}

\begin{lemma}\label{lem6;Vtriple}
Let $(X,Z)$ be a $V$-pair over $S\in \Sm$ and $f:X'\to X$ be 
an \'etale morphism such that $Z':=X'\times_X Z \isom Z$. 
If $(X',Z')$ is a pre-$V$-pair over $S$,
it is a $V$-pair over $S$ (cf. Definition \ref{def;Vtriple}).
\end{lemma}
\begin{proof}
We need check the condition $(iii)$ of Definition \ref{def;Vtriple} for $(X',Z')$.
It suffices to show the following diagram
\[\xymatrix{
Z' \ar[r]^{\hskip -20pt\Delta_{Z'}} \ar[d]^{f_{|Z'}} & X'\times_S Z' \ar[d]^{f} \\
Z \ar[r]^{\hskip -15pt\Delta_{Z}} & X\times_S Z \\
}\]
 is Cartesian, which can be seen from 
$Z'=X'\times_X Z \isom Z$.
\end{proof}

\begin{lemma}\label{lemP1-1}
Let $\eta$ be as in Theorem \ref{thm-P1} and 
$(X,Z)$ be a $V$-pair over $\eta$.
Let $f:X'\to X$ be an \'etale map such that $X'\times_X |Z| \isom |Z|$. 
Let $F\in \lsCIt$ (cf. Definition \ref{def;Xi}). Put $Z'=X'\times_X Z$.
\begin{itemize}
\item[(1)]
The natural maps $F(X) \to F(X,Z)$ and $F(X') \to F(X',Z')$ are injective.
\item[(2)]
Assume further that $F$ is \lssemipure. Then the natural map
\[F(X,Z)/F(X) \to F(X',Z')/F(X')\]
is an isomorphism.
\item[(3)]
Assume $(X,nZ)$ is a $V$-pair over $\eta$ for any $n\in\Z_{>0}$.
Then the natural map
\[ F(X-Z)/F(X) \to F(X'-Z')/F(X') \]
is an isomorphism.
\end{itemize}
\end{lemma}
\begin{proof}
The assumption implies $Z' \isom Z$, and
$(X',Z')$ is a $V$-pair over $\eta$ by Lemma \ref{lem6;Vtriple} and Remark \ref{rem;def;Vtriple}(1). Hence (1) and (2) follow from Theorem \ref{thm2;Vtriple} and
Corollary \ref{cor;thm2;Vtriple}.
To show (3) we may replace $F$ by $\tilde{F}$ from Lemma \ref{lem;Mrecsemipure}
to assume $F$ is semipure. Then it follows from (2) in view of Lemma \ref{lem1;Mrec}.
\end{proof}

{\it Proof of Theorem \ref{thm-P1}}\; 
Let $Z\subset X\subset \P^1_\eta$ be as in Theorem \ref{thm-P1} and $\xi$ be the generic point of $\P^1_\eta$.
For a closed point $x\in X$ let $\hen X x$ be the henselization of $X$ at $x$ and $Z_x=Z\times_X \hen X x$. By \cite[Ch. XI, Th. 1]{ray}, we have
$\hen X Z= \prod_{x\in Z} \hen X x$. 
In view of Lemma \ref{lemP1-0}, by passing to the limit over all Nisnevich neighbourhoods of $Z\subset X$ in Lemma \ref{lemP1-1}, one obtains isomorphisms
\begin{equation}\label{eq;lemP1-0}
 F(X,Z)/F(X) \simeq  \underset{x\in Z}{\bigoplus}\;
F(\hen X x,Z_x)/F(\hen X x),
\end{equation}
\begin{equation}\label{eq2;lemP1-0}
 F(X-W)/F(X) \simeq  \underset{x\in W}{\bigoplus}\;
F(\hen X x-x)/F(\hen X x)
\end{equation}
for any finite set $W\subset X$ of closed points.
Putting $U=X-|Z|$ we have a commutative diagram
\[\xymatrix{
0\ar[r] & F(X) \ar[r]\ar[d] & F(\xi) \ar[r]\ar[d] & 
\underset{x\in X}{\bigoplus}\;\frac{F(\hen X x -x)}{F(\hen X x )}\ar[r]\ar[d] & 0\\
& F(X,Z) \ar[r] & F(\xi) \ar[r]& 
\underset{x\in U}{\bigoplus}\;\frac{F(\hen X x -x)}{F(\hen X x )} \oplus 
\underset{x\in Z}{\bigoplus}\;\frac{F(\hen X x -x)}{F(\hen X x,Z_x)}. \\
}\]
The upper horizontal sequence is exact and obtained by taking the colimit of \eqref{eq2;lemP1-0} over all finite sets $W\subset X$.
Write $G=F_{(X,Z)}$. Noting \eqref{eq;lemP1-0} an easy diagram chase implies the exactness of  
\begin{equation}\label{eq;thm-P1}
0\to G(X) \to G(\xi) \to \underset{x\in X}{\bigoplus}\; G(\hen X x -x)/G(\hen X x) \to 0.
\end{equation}
On the other hand, we have a localization exact sequence
\begin{multline*}
0\to \underset{x\in X}{\bigoplus}\; H^0_x(X_{\Nis},G_{\Nis}) \to G_{\Nis}(X) \to G_{\Nis}(\xi) \to \\
\underset{x\in X}{\bigoplus}\; H^1_x(X_{\Nis},G_{\Nis})  \to H^1(X_{\Nis},G_{\Nis}) \to 0.
\end{multline*} 
We have
\[
\begin{aligned}
 H^1_x(X_{\Nis},G_{\Nis})
&=\Coker\big(G_{\Nis}(\hen X x) \to G_{\Nis}(\hen X x -x)\big)\\
&=\Coker\big(G(\hen X x) \to G(\hen X x -x)\big),\\
\end{aligned}
\]
\[
\begin{aligned}
 H^0_x(X_{\Nis},G_{\Nis})
&=\Ker\big(G_{\Nis}(\hen X x) \to G_{\Nis}(\hen X x -x)\big)\\
&=\Ker\big(G(\hen X x) \to G(\hen X x -x)\big)=0\\
\end{aligned}
\]
where the last equality follows from the semipurity of $F$ and Theorem \ref{thm;localinjectivity}. Thus we get an exact sequence
\begin{equation}\label{eq2;thm-P1}
0\to G_{\Nis}(X) \to G(\xi) \to \underset{x\in X}{\bigoplus}\; G(K^h_x)/G(\hen X x) \to 
H^1(X_{\Nis},G_{\Nis}) \to 0.
\end{equation}
Comparing \eqref{eq;thm-P1} and \eqref{eq2;thm-P1}, this proves Theorem \ref{thm-P1}.

\section{Contractions}\label{contraction}

\def\ul0{\underline{0}}

In this section we fix an integral affine $S\in \tSm$ and an effective Cartier divisor $D\subset S$.
For a $S$-scheme $X$ write $D_X=D\times_S X$. 

\begin{defn}\label{def;niceVpair}
A $V$-pair $\nu=(X,Z)$ over $S$ is \emph{nice} if 
\begin{itemize}
\item[$(i)$]
$Z$ is reduced and \'etale over $S$.
\item[$(ii)$]
$\nu_n=(X,nZ)$ is a $V$-pair over $S$ for all $n\in \Z_{>0}$.,
\end{itemize}
Here $nZ\hookrightarrow X$ is the $n$-th thickening of $Z\hookrightarrow X$. 
\end{defn}

\begin{lemma}\label{lem;;niceVpair}
Let $X\subset \A^1=\Spec k[t]$ be a dense open subset containing the origin. For any $n\in \Z_{>0}$ 
\[ \alpha_n := (S\times X,S\times \Lambda_n).\]
is a nice $V$-pair over $S$, where $\Lambda_n=\Spec k[t]/(t^n)$. 
\end{lemma}
\begin{proof}
This follows from Lemma \ref{lemP1-0} and Remark \ref{rem;def;Vtriple}(2).
\end{proof}

\begin{defn}\label{def;formalframe}
A formal frame of a nice $V$-pair $\nu=(X,Z)$ over $S$ is an isomorphism over $S$  of formal schemes
\[ \hep: \hat{X}_{|Z} \simeq Z\times \Spf k[[t]],\]
where $\hat{X}_{|Z}$ is the formal completion of $X$ along $Z$.
It gives rise to a compatible system of isomorphisms for $n\in \Z_{>0}$:
\begin{equation}\label{eq;def;formalframe} 
{\hep}_n: n Z \simeq Z\times \Lambda_n \qwith \Lambda_n=\Spec k[t]/(t^n).
\end{equation}
\end{defn}

\def\pFcont{\sigma(F)}
\def\pFcontt#1{\sigma^{(#1)}(F)}
\def\ptFcont{\sigma_{\A^1}(F)}
\def\ptFcontt#1{\sigma_{\A^1}^{(#1)}(F)}

\def\pGcont{\sigma(G)}
\def\pGcontt#1{\sigma^{(#1)}(G)}
\def\ptGcont{\sigma_{\A^1}(G)}
\def\ptGcontt#1{\sigma_{\A^1}^{(#1)}(G)}

\begin{defn}\label{def2p;contraction}
For $F\in \ulMPST$ and $n\in \Z_{>0}$, we define $\pFcont, \pFcontt n \in \ulMPST$ by:
\[ \pFcontt n(\fX)=\Coker\big(F(\fX\otimes (\P^1,\infty))\to 
F(\fX\otimes (\P^1, n0 +  \infty))\big),\]
\[ \pFcont(\fX)=\Coker\big(F(\fX\otimes (\P^1,\infty))\to 
F(\fX\otimes (\P^1-0,\infty))\big),\]
where $\fX\in \ulMCor$. 
\end{defn}
\medbreak

Assume $F$ is $\lscube$-invariant. Then, for $\fX\in\ulMCorls$, we have natural isomorphisms
\[ \pFcontt n (\fX)=\Coker\big(F(\fX) \rmapo {pr^*} 
F(\fX\otimes (\P^1,n0 + \infty))\big),\]
\[ \pFcont(\fX)=\Coker\big(F(\fX) \rmapo {pr^*} F(\fX\otimes (\P^1-0,\infty))\big)\]
with $pr^*$ induced by the projection $pr: \fX\otimes (\P^1, n0 +  \infty)\to\fX$.
Note that $pr^*$ is split injective and its inverse is given by the pullback along
a section of $\P^1-\{0,\infty\} \to \Spec(k)$. Hence there are natural isomorphisms
\begin{equation}\label{eq1;def2;contraction} 
\begin{aligned}
 &F(\fX\otimes (\P^1, n0 +  \infty)) \simeq \pFcontt n(\fX) \oplus F(\fX),\\
 &F(\fX\otimes (\P^1-0,\infty)) \simeq \pFcont(\fX) \oplus F(\fX)
\end{aligned}
\end{equation}
which are functorial in $\fX\in \ulMCorls$.
\medbreak

\begin{lemma}\label{lem1;contraction}
Let $F\in \lsCItsp$ (cf. Definition \ref{def;Xi}). We have $\pFcontt n,\pFcont\in \lsCItsp$.
If moreover $F\in \ulMNSTls$, then $\pFcontt n,\pFcont\in \ulMNSTls$ (cf. Definition \ref{CIrspNisdef}).
\end{lemma}
\begin{proof}
This follows from \eqref{eq1;def2;contraction} noting Definitions \ref
{def:m-rec} and \ref{def:semipure} together with Lemmas \ref{lem1;Mrec}(1) and \ref{lem1;semipure}.
\end{proof}

\begin{defn}\label{def2;contraction}
For $F\in \lsCItspNis$ (cf. Definition \ref{CIrspNisdef}) and $n\in \Z_{>0}$, write 
\[ \Fcont = {\pFcont}\qaq \Fcontt n ={\pFcontt n}.\]
By Lemma \ref{lem1;contraction} the association $F\to \Fcont$ gives an endofunctor on $\lsCItspNis$. We define $\Fccc i$ for $i>0$ inductively by $\Fccc i=(\Fccc {(i-1)})_{-1}$.

\end{defn}
\medbreak

\begin{lemma}\label{lem0;thm2-P1}
Let $\phi:(X,0_X) \to (\A^1,0)$ be an affine Nisnevich neighbourhood such that
$\phi^{-1}(0)=\{0_X\}$.
Let $\phib:\Xb\to \P^1$ be the normalization of $\P^1$ in $X$.
For any effective Cartier divisor $\Xinf\subset \Xb$ such that $|\Xinf|=\Xb-X$, 
there exists a rational function $f$ on 
$\A^1 \times X$ admissible for $\big((X,n0_X),(\A^1,n 0)\big)$ 
for all $n\in \Z_{>0}$ (cf. Definition \ref{def;specialfunction}) such that 
\begin{equation}\label{eq;lem0;thm2-P1}
\theta_f=\div_{\A^1\times X}(f) \in \ulMCor\big((\P^1, \infty+ n0),(\Xb,\Xinf + n0_X))\big).
\end{equation}
\end{lemma}
\begin{proof}
First we note that $(\A^1,0)$ and $(X,0_X)$ are $V$-pairs over $k$ 
by Lemmas \ref{lemP1-0} and \ref{lem6;Vtriple}. 
Write
\[\A^1=\Spec k[t],\; \P^1-\{1\}=\Spec k[s]\]
for variables $t,s$ with $t=s/(s-1)$.
Take an affine open neighbourhood $W=\Spec B$ of $\{0_X\}\cup \Xinf$ in $\Xb$
and $\pi_0,\piinf\in B$ such that 
\[ \div_W(\pi_0) = 0_X \qaq \div_W(\piinf) = \Xinf\]
and set $\tau=\phi^*(t)\in B[1/\piinf]$. By the condition $\phi^{-1}(0)=\{0_X\}$, 
\begin{equation}\label{eq4;lem0;thm2-P1}  
\tau=\pi_0 b /\pi_\infty^r\;\text{ for some }r>0,\; b\in B, \text{ and }
B/(\pi_0^n) = B[1/\piinf]/(\tau^n).
\end{equation}
Take $u,v\in B$ such that $u\pi_0^n+v\piinf^r=1$ (this is possible since $0_X\not\in \Xinf$). Put
\begin{multline*}
 f= v\piinf^r(t-v\pi_0b)  +  u\pi_0^n = 1+ v\piinf^r( t -1- v\pi_0 b)
\in B[t] =\Gamma(\A^1\times W,\sO).
\end{multline*}

\begin{claim}\label{claim;lem0;thm2-P1}
$f$ is admissible for $\big((X,n0_X),(\A^1,n0)\big)$.
\end{claim}

Indeed, $f$ clearly satisfies Definition \ref{def;specialfunction}(1).
We have $f\equiv 1\mod (\piinf)\subset B[t]$ so that $f$ satisfies Definition \ref{def;specialfunction}(3). 
Thanks to \eqref{eq4;lem0;thm2-P1}, we have 
\begin{equation}\label{eq3;lem0;thm2-P1}  B[t]/(\pi_0^n, f)= 
B[1/\piinf][t]/(\tau^n, t-v\pi_0b)= B[1/\piinf][t]/(\tau^n, t-\tau),
\end{equation}
where the last equality follows from
\[ t-v\pi_0b = t-v\piinf^r\tau = t-(1-u\pi_0^n)\tau\equiv t-\tau\mod (\tau^n)\subset B[1/\piinf][t],\]
where the last congruence follows from \eqref{eq4;lem0;thm2-P1}.
Since the diagonal $\Delta_{n0_X}: n0_X \to \A^1\times n0_X$ is induced by the map
\[B[1/\piinf][t]/(\tau^n) \to B[1/\piinf]/(\tau^n)\;;\; t\to \tau,\] 
\eqref{eq3;lem0;thm2-P1} implies that $f$ satisfies Definition \ref{def;specialfunction}(2) and the claim is shown. 

Note 
\[f=\frac{1}{s-1}\big(s-1+ v\piinf^r (1-(s-1)v\pi_0b)\big) \]
so that the closure $\Theta$ of $(\theta_f)_{|(\A^1-\{1\})\times W}$ in
$(\P^1-\{1\})\times W$ is
\[\div_{(\P^1-\{1\})\times W}\big(s-1+ v\piinf^r (1-(s-1)v\pi_0b)\big)\]
and hence
$(s-1)/\piinf=-v\piinf^{r-1}(1-(s-1)v\pi_0b)$ is regular on $\Theta$.
This proves that $f$ satisfies \eqref{eq;lem0;thm2-P1} and the proof of Lemma \ref{lem0;thm2-P1} is complete.
\end{proof}
\bigskip

For $F\in \ulMPST$ and $n\in \Z_{>0}$, we define $\ptFcont, \ptFcontt n \in \ulMPST$ by:
\[ \ptFcontt n(\fX)=\Coker\big(F(\fX\otimes (\A^1,\emptyset))\to F(\fX\otimes (\A^1, n0))\big),\]
\[ \ptFcont(\fX)=\Coker\big(F(\fX\otimes (\A^1,\emptyset))\to F(\fX\otimes (\A^1-0,\emptyset))\big),\]
where $\fX\in \ulMCor$. 
We have natural maps in $\ulMPST$:
\[ \pFcont \to \ptFcont,\quad \pFcontt n \to \ptFcontt n\]
induced by the maps $(\A^1,n0) \to (\P^1, n0 +  \infty)$ and 
$\A^1-0 \to (\P^1-0,\infty)$ in $\ulMCor$.

\begin{lemma}\label{lem2;contraction}
Take $F\in \lsCItsp$ and $\fZ=(Z,D)\in \ulMCorls$.
Then we have isomorphisms
\begin{equation}\label{eq2;thm;contraction}
\pFcontt n (\fZ) \isom  \ptFcontt n (\fZ),
\end{equation}
\begin{equation}\label{eq2-1;thm;contraction}
\pFcont (\fZ)\isom  \ptFcont(\fZ).
\end{equation}
\end{lemma}
\begin{proof}
\eqref{eq2-1;thm;contraction} is deduced from \eqref{eq2;thm;contraction} using Lemma \ref{lem1;Mrec}(1). We show the injectivity of \eqref{eq2;thm;contraction}.
We claim that this is reduced to the case $\fZ=(\Spec k,\emptyset)$.
Indeed consider $G=\uHom_{\MPST}(\Ztr(\fZ),F)$.
By Lemma \ref{lem1;Mrec}(2) and Lemma \ref{lem1;semipure}(2), we have
$G\in \lsCItsp$. Hence the claim follows from the natural isomorphisms
\[\pFcontt n (\fZ)\simeq \pGcontt n (\Spec k,\emptyset),\quad
\ptFcontt n (\fZ) \simeq \ptGcontt n (\Spec k,\emptyset).\]
We have a commutative diagram
\[\xymatrix{
F(\P^1, n0+\infty)/F(\P^1,\infty) \ar[r]^{\hskip 20pt\eqref{eq2;thm;contraction}}\ar[d]^{\hookrightarrow} &
F(\A^1,n0)/F(\A^1,\emptyset)\ar[d]^{\hookrightarrow}  \\
F(\P^1-0,\infty)/F(\P^1,\infty) \ar[r]^{\hskip 20pt\eqref{eq2-1;thm;contraction}} &
F(\A^1-0,\emptyset)/F(\A^1,\emptyset) \\}\]
By the semipurity of $F$ the vertical maps are injective by Lemma \ref{lem1;semipure}(1). Hence it suffices to show the injectivity of \eqref{eq2-1;thm;contraction}. The square
\[\xymatrix{
(\A^1-0,\emptyset) \ar[r]\ar[d] & (\A^1,\emptyset)\ar[d]\\
(\P^1-0,\infty) \ar[r] & (\P^1,\infty)\\
}\]
comes by pullbacks from the elementary Nisnevich square in $\Sm$:
\[\xymatrix{
\A^1-0 \ar[r]\ar[d] & \A^1\ar[d]\\
\P^1-0 \ar[r] & \P^1\;.\\
}\]
Thus, if $F\in \ulMNST$, the sheaf condition for $F_{(\P^1,\infty)}$ implies that the sequence
\begin{multline*}
 0\to F(\P^1,\infty) \to F(\P^1-0,\infty)\oplus F(\A^1,\emptyset) \to F(\A^1-0,\emptyset), \\
\end{multline*}
is exact, which implies the desired injectivity.
Even if $F$ is not necessarily in $\ulMNST$, the above sequence is exact: Indeed the second and third terms do not change after replacing $F$ by $F_\Nis$ thanks to Theorem \ref{thm-P1}(1). The same is true for the first term by the cube-invariance of $F$. This completes the proof of the injectivity.

Next we prove the surjectivity of \eqref{eq2;thm;contraction} (in case $\fZ=(\Spec k,\emptyset)$). For this  we use the $V$-pairs $(C,0_C)$ and 
$(\A^1,0)$ over $\Spec k$, where $C=\A^1$ and $0_C=0\in \A^1$ is the origin (cf. Lemma \ref{lemP1-0}). 
Let $\Cbb=\P^1$ be the smooth compactification of $C$.
Thanks to the semipurity of $F$ noted above, 
Theorem \ref{thm2;Vtriple}(2) implies that for any rational function $f$ on $\A^1\times C$ admissible for $\big((C, n0_C),(\A^1,n0)\big)$ (cf. Definition \ref{def;specialfunction}),
\[ \theta_f=\div_{\A^1\times C}(f) \in \ulMCor\big((\A^1,n0)),(C, n0_C)\big)\]
induces an isomorphism
\begin{equation*}\label{eq1.1;thm;contraction}
\theta_f^* : F(C,n0_C)/F(C,\emptyset)\isom F(\A^1,n0)/F(\A^1,\emptyset),
\end{equation*}
which is independent of $f$.
By Lemma \ref{lem0;thm2-P1} applied to the identity $C\to \A^1$,
for any effective divisor $\Cinf$ supported in $\infty=\Cbb-C$, 
there exists such $f$ as above that
\begin{equation*}
 \theta_f\in \ulMCor((\P^1, n0 +  \infty), (\Cbb, n0_C + \Cinf)).
\end{equation*}
Note that one cannot take $f$ to be the diagonal since $\Cinf$ may not be reduced
even though it is supported on $\infty$.
This implies that the image of the composite map
\begin{multline*}
F(\Cbb, n0_C  +  \Cinf) \to F(C,n0_C) \rmapo{\theta_f^*} F(\A^1,n0)/F(\A^1,\emptyset)
\end{multline*}
is contained in the image of \eqref{eq2;thm;contraction}, where the first map is induced by the map $(C,n0_C)\to (\Cbb,n0_C  +  \Cinf)$ in $\ulMCor$.
This implies the desired assertion since $\theta_f^*$ is independent of $f$ and   
we have by Lemma \ref{lem1;Mrec}(1)
\[ F(C,n0_C) \simeq  \indlim {\Cinf} F(\Cbb, n0_C  +  \Cinf),\]
where the limit is over all such $\Cinf$ that $|\Cinf|=\Cbb-C$.
\end{proof}

\begin{lemma}\label{lem3;contraction}
Let $\fZ=(Z,D)$ be as Lemma \ref{lem2;contraction} and $F\in \lsCItsp$.
Then $\Fccc i(\fZ)$ for $i>0$ is naturally isomorphic to 
\[\Coker\big(\underset{1\leq j\leq i}{\bigoplus}\;
F((\A^1-0)^{j-1}\times \A^1 \times(\A^1-0)^{i-j}\otimes \fZ)\to  F((\A^1-0)^i\otimes \fZ)\big).\]
\end{lemma}
\begin{proof}
This is easily shown by using Lemma \ref{lem2;contraction} repeatedly.
\end{proof}


\begin{thm}\label{thm;contraction}
Let $\nu=(X,Z)$ be a nice $V$-pair over $S$ equipped with a formal frame ${\hep}$. 
Assume that $Z\to S$ is an isomorphism and that $|D|$ is a simple normal crossing divisor on $S$. 
Put $\fZ=(Z,D_Z)\in \ulMCorls$. Take $F\in \lsCItspNis$.
\begin{itemize}
\item[(1)]
There exist a compatible system of isomorphisms for $n\in \Z_{>0}$:
\[ \theta_{{\hep}_n}: F(X,nZ+D_X)/F(X,D_X)\simeq \Fcontt n(\fZ)\]
depending on ${\hep}_n$ from \eqref{eq;def;formalframe}. Moreover, for $m\geq n$,
the diagram
\[\xymatrix{
  F(X,nZ+D_X)/F(X,D_X) \ar[r]\ar[d] & \Fcontt n(\fZ)\ar[d] \\
  F(X,mZ+D_X)/F(X,D_X) \ar[r] & \Fcontt m(\fZ) \\
}\]
is commutative, where the vertical maps are the natural ones.
\item[(2)]
There exist an isomorphism
\[ \theta_{\hep}: F(X-Z)/F(X)\simeq \Fcont (Z)\]
depending on ${\hep}$.
\end{itemize}
\end{thm}
\begin{proof}
This is a direct consequence of Theorem \ref{thm2;Vtriple}(2),
Lemmas \ref{lem;;niceVpair} and \ref{lem2;contraction}.

\end{proof}


\begin{rk}\label{rm;thm;contraction}
Let the assumption be as in Theorem \ref{thm;contraction}.
\begin{itemize}
\item[(1)]
Let $\gamma:S'\to S$ be a morphism in $\tSm$ with $S'$ integral affine and 
$D'\subset S'$ be an effective Cartier divisor such that $|D'|$ is a simple normal crossing divisor on $S'$ and $D'\geq \gamma^* D$. 
Then the base change $(X',Z')=(X,Z)\times_S S'$ is a nice $V$-pair (cf. Remark \ref{rem;def;Vtriple}) and a formal frame ${\hep}$ of $(X,Z)$ induces a formal frame ${\hep}'$ of $(X',Z')$. Put $D'_{X'}=X'\times_{S'} D'$.
In view of Remark \ref{rk;specialfunction} the following diagram is commutative:
\[\xymatrix{
 F(X,nZ+ D_Z)/F(X,D_X)\ar[r]^{\hskip 20pt\theta_{\hep}}\ar[d]^{\gamma^*} & 
\Fcontt n (Z,D_Z) \ar[d]^{\gamma^*} \\
 F(X',nZ'+ D'_{X'})/F(X',D'_{X'})\ar[r]^{\hskip 30pt\theta_{{\hep}'}} & \Fcontt n (Z',D_{Z'})  \\
}\]
\item[(2)]
Let $\pi:X'\to X$ be an \'etale morphism which induces an isomorphism
$Z':=X'\times_X Z \simeq Z$.
Assume $(X',Z')$ is a nice $V$-pair 
(By Lemma \ref{lem6;Vtriple} and Remark \ref{rem;def;Vtriple}(1), this is the case if $S$ is the spectrum of a field).
The formal frame ${\hep}$ of $(X,Z)$ induces a formal frame ${\hep}'$ of $(X',Z')$.
By Corollary \ref{cor;thm2;Vtriple} the diagram
\[\xymatrix{
 F(X,nZ+ D_Z)/F(X,D_X)\ar[r]^{\hskip 20pt\theta_{\hep}}\ar[d]^{\pi^*} & 
\Fcontt n(Z,D_Z) \ar[d]^{\pi^*} \\
 F(X',nZ'+ D_{X'})/F(X',D_{X'})\ar[r]^{\hskip 30pt\theta_{{\hep}'}} & 
\Fcontt n (Z',D_{Z'})  \\
}\]
 is commutative and the vertical maps are isomorphisms.
\end{itemize}
\end{rk}

\section{Fibrations with coordinate}\label{fibration}

\def\ul0{\underline{0}}
In this section we assume $k$ is perfect and infinite.

\begin{lemma}\label{lem;etalepoint}
Take $W\in \Sm$ and a point $e\in W$ with $K=k(e)$.
Let $Z\subset W$ be a reduced closed subscheme containing $e$. 
Let $\sX=\hen W e$ be the henselization of $W$ at $e$.
Then there is an isomorphism $\sX\simeq \Spec K\{x_1,\dots,x_d\}$ such that
there exists an open dense subset $\sU$ of $\sZ:=Z\times_W \sX$ for which 
$\sU \to \sX \to \Spec(K)$ is essentially smooth.
\end{lemma}
\begin{proof}
Let $E\subset W$ be the closure of $e$.
Let $r=\dim(E)$ and $d:=\codim_W(E)$. 
We may assume $W$ is affine and choose an closed embedding $W\hookrightarrow \A^N$,
where $\A^N$ is the $N$-dimensional affine space over $k$.
By \cite[Pr.5.5]{bv} there exist linear projection $\A^N \to \A^r$ satisfying
the following conditions:
\begin{itemize}
\item
the induced map $W \to \A^r$ is smooth at $e$,
\item
the induced map $E \to \A^r$ is finite and generically \'etale,
\item
there is an open dense $U\subset Z$ smooth over $\A^r$.
\end{itemize}
Put $L=k(\A^r)$ and consider the composite
\[\Spec(K)\rmapo{e} W\times_{\A^r} L\to \Spec(L),\]
where the second map is induced by $W\to \A^N\to \A^r$.
By the construction the map is \'etale so that there is an affine Nisnevich neighbourhood $(V,e_V)\to (W\times_{\A^r} L,e)$ and a smooth morphism 
\[\pi: V\to \Spec(K)\]
such that $e_V: \Spec(K)\to V$ is a section of $\pi$.
Then $U_V=U\times_W V$ is smooth over $K$. 
Choose a closed embedding $V \hookrightarrow \A^M_K=\A^M\otimes_k K$ which maps $e_V$ to the origin. 
By \cite[Pr.5.5]{bv} there exist linear projection
$\A^M_K \to \A^{d}_K$ such that the induced map $V\to \A^{d}_K$ is finite 
and \'etale at $e_V$. Noting that $\sX=\hen W e\simeq \hen V {e_V}$,
this proves the lemma.
\end{proof}
\medbreak

Let $\sX$ be as in Lemma \ref{lem;etalepoint} and $0_\sX\in \sX$ be the closed point. 

\begin{defn}\label{goodfibration}
\begin{itemize}
\item[(1)]
A fibration of $\sX$ with coordinate $t$ is an essentially \'etale morphism 
$\psi_t: \sX \to \A^1_\sS=\sS[t]$ where
\[\sS=\Spec K\{y_1,\dots,y_{d-1}\}\]
such that $\psi_t$ induces $\sX\simeq \hen{(\A^1_\sS)} {(0_\sS,t)}$,
where $0_\sS$ is the closed point.
We write $\psi: \sX\to \sS$ for the composite of $\psi_t$ and the projection $\A^1_\sS\to\sS$.
\item[(2)]
Let $\sH\subset \sX$ be a regular divisor.
A fibration $\psi_t:\sX\to \sS[t]$ with coordinate $t$ is a $\sH$-fibration if
there exists a regular divisor $H_\sS\subset \sS$ such that $\sH=\psi^{-1}(H_\sS)$.
\item[(3)]
A reduced closed subscheme $\sZ\subset \sX$ is admissible for $\psi$ if $\sZ$ does not contain $\psi^{-1}(0_\sS)$.
If moreover every irreducible component of $\sZ$ is generically \'etale over 
its image in $\sS$, we say $\sZ$ is strongly admissible for $\psi$.
\end{itemize}
\end{defn}  

\medbreak

\begin{lemma}\label{lem;existencegoodfibration}
Let $\sX$ be as in Lemma \ref{lem;etalepoint} and $\sH\subset \sX$ be a regular divisor.
Let $\sZ\subset \sX$ be a reduced closed subscheme such that $\sH\not\subset \sZ$.
Then there exists a $\sH$-fibration $\psi_t:\sX\to \sS[t]$ with coordinate $t$
such that $\sZ$ is strongly admissible for $\psi$.
\end{lemma}
\begin{proof}
Let $\sX\simeq  \Spec K\{x_1,\dots,x_d\}$ be as in Lemma \ref{lem;etalepoint}.
Take a regular system of parameters $(y_1,\dots,y_d)$ of $\sX$ such that $\sH$ is the divisor of
$y_1$. Then $\sX$ is the henselization of $L:=\A^d_K=\Spec K[y_1,\dots,y_d]$ at 
the origin and $\sH$ is the pullback of the hyperplane $H=\{y_1=0\}\subset L$.
Let $Z\subset \A^d_K$ be the closure of the image of $\sZ$.
By Lemma \ref{lem;etalepoint} we may suppose that there is an open dense subscheme of $Z$ which is smooth over $K$. It then suffices to find a linear projection
$\phi: L \to L'\simeq \A^{d-1}_K$ satisfying the conditions:
\begin{itemize}
\item
$Z$ does not contain $\phi^{-1}(0_{L'})$, where $0_{L'}$ is the origin of $L'$,
\item
$H$ is the pullback of a hyperplane in $L'$,
\item
every irreducible component of $Z$ is generically \'etale over its image in $L'$.
\end{itemize}
Consider the open immersion $L\hookrightarrow \Lb=\P^d_K$ and put $\Linf=\Lb-L$. 
Let $\Zb\subset \Lb$ be the closure and $\Zinf=\Zb-Z$.
We define $\Hinf=\Hb-H$ similarly. Let $\Linf(K)$ be the set of $K$-rational points of
$\Linf$. For $v\in \Linf(K)$ we have a linear projection 
\[\phi_v: L \to L_v:= L/\ell_v\cap L \simeq \A^{d-1}_K,\]
where $\ell_v\subset \Lb$ is the line passing through $v$ and the origin $0\in L$. Note $\phi_v^{-1}(0_{L_v})=\ell_v\cap L$ and that 
if $v\in \Hinf$, $H$ is the pullback of a hyperplane in $L_v$. 
By the assumption, for any irreducible component $Z_\lambda$ of $Z$, we have
$Z_\lambda\cap H\subsetneq H$. Hence there is a dense open subscheme $\Hinf^\circ\subset \Hinf$ such that $Z$ does not contain $\phi_v^{-1}(0_{L_v})=\ell_v\cap L$ for $v\in \Hinf^\circ(K)$.
Since $K$ is infinite and $\Hinf\simeq \P^{d-1}_K$, $\Hinf^\circ(K)$ is non-empty.

To ensure that $\phi_{v}$ induces a generically \'etale map from $Z$ to its image in $L_v$, choose an irreducible component $Z_\lambda$ of $Z$.
%
Choose a smooth closed point $x$ on $Z_\lambda$ such that $K(x)$ is separable over $K$ and put $y=\phi_v(x)$. It suffices to show the following. 

\begin{claim}\label{calim;existencegoodfibration}
After possibly changing the choice of $x\in Z_\lambda$ as above and that of a regular system of parameters $(y_1,y_2,\dots,y_d)$ of $\sX$ in the beginning, there exists a dense open subscheme $U_\lambda\subset \Hinf$ such that
for every $v\in U_\lambda(K)$, the map $\phi_{v,Z_\lambda}: Z_\lambda\to L_v$ induced by $\phi_v$ is \'etale (resp. an immersion) at $x$ if $\dim(Z_\lambda)=d-1$ (resp. $\dim(Z_\lambda)<d-1$).
\end{claim}
\def\tZl{\tilde{Z}_\lambda}

We follow the argument from \cite[5.5]{bv}.
For our purpose, we may replace $L$ and $Z$ by its base changes via $K\to K(x)$ so that we may assume $x$ is a $K$-rational point of $Z$.
Let $T_x$ be the tangent space of $Z_\lambda$ at $x$ considered as an affine subspace of $L$. We define $\Txinf=\ol{T}_{x}-T_x\subset\Linf$ as before. 
Let $Z_\lambda-\{x\} \to \Linf$ be the map which sends $z$ to the intersection of $\Linf$ and the line in $\Lb$ passing through $x$ and $z$ and let $\tZl\subset\Linf$ be the closure of its image. For $v\in \Linf(K)$, if $v\not\in \Txinf$, 
$d\phi_{v,\lambda}: T_x \to  T_{L_v,y}$ is injective so that $\phi_v:Z_\lambda\to L_v$ is unramified at $x$. In case $\dim(Z_\lambda)=d-1$, $\phi_{v,Z_\lambda}$ is dominant so that it is \'etale at $x$ by \cite[Exp. 1, Cor. 9.11]{SGA1}.
In case $\dim(Z_\lambda)<d-1$, if $v\not\in \tZl\cup \Txinf$, $\phi_v:Z_\lambda\to L_v$ is an immersion at $x$. Note
\[\dim(\Txinf)<\dim(Z_\lambda),\;
\dim(\tZl)\leq \dim(Z_\lambda),\; \dim(\Hinf)=d-2.\] 
Thus, if $\dim(Z_\lambda)<d-2$, we may take 
$U_\lambda=\Hinf^\circ\backslash (\tZl\cup \Txinf)$.
It remains to treat the following cases:
\begin{itemize}
\item[(i)]
$\dim(Z_\lambda)=d-1$ and $\Hinf\subset \Txinf$ for all smooth closed points of $Z_\lambda$ such that $K(x)$ is separable over $K$.
\item[(ii)]
$\dim(Z_\lambda)=d-2$ and $\Hinf=\tZl$ (note that $\tZl$ is irreducible).
\end{itemize}
In case (i), $Z_\lambda\subset L=\A^d_K$ is defined by a polynomial 
$F\in K[y_1,\dots,y_d]$, and $\Hinf\subset \Txinf$ implies
$\partial F/\partial y_i \in \fm_x$ for $i=2,\dots, d$,
where $\fm_x\subset K[y_1,\dots,y_d]$ is the maximal ideal of $x$.
If this holds for all smooth closed points $x$ of $Z_\lambda$ such that $K(x)$ is separable over $K$, we get
\[\partial F/\partial y_i\equiv 0\in K[y_1,\dots,y_d] \qfor i=2,\dots, d\]
since such points are dense in $Z_\lambda$.
If $\ch(K)=0$, this implies $F=c y_1$ with $c\in K$, which contradicts 
the assumption that $H\not\subset Z_\lambda$. If $\ch(K)=p>0$, we can write
\[ F = a_0^p + a_1^p y_1 + \cdots a_n^p y_1^n\qfor a_i\in K[y_2,\dots,y_d].\]
The last case is excluded by changing the choice of a regular system of parameters $(y_1,y_2,\dots,y_d)$ of $\sX$ in the beginning of the proof,
e.g. changing it to $(y_1(1+y_2),y_2,\dots,y_d)$.

In case (ii), $Z_\lambda$ is contained in the hyperplane in $\Lb$ containing $\Hinf$ and $x$.
This is also excluded by changing the choice of a regular system of parameters as above.
This completes the proof of Lemma \ref{lem;existencegoodfibration}.
\end{proof}
\medbreak

Let $\psi:\sX \to \sS[t]$ be as in Definition \ref{goodfibration}. We have
\begin{equation}\label{eq0;lem1;purity}
\sX \simeq \lim_{(X,0_X)}\; X,
\end{equation}
where $(X,0_X)$ range over all Nisnevich neighbourhoods of $(\A^1_\sS,0)$.
We have induced maps 
\[\psi_X: (\sX,0_\sX) \to (X,0_X)\qaq p_X: (X,0_X)\to (\sS,0_S).\] 

\begin{lemma}\label{lem1;purity}
For a closed subscheme $T\subset X$, $\psi_X^{-1}(T)\not=\emptyset$ 
if and only if $T$ contains $0_X\in X$. 
\end{lemma}
\begin{proof}
The if-part is obvious since $0_X\in X$ is in the image of $\psi_X$. 
Assume $T\times_X \sX\not=\emptyset$.
Then it is a closed subscheme of $\sX$ which is local. Hence it must contain 
$0_\sX\in \sX$. Hence $T$ contains $0_X=\psi_X(0_\sX)$.
\end{proof}

Consider the following condition for $(X,0_X)$ in \eqref{eq0;lem1;purity}:
\begin{enumerate}
\item[$(*)$]
$p_X^{-1}(0_\sS)$ is irreducible and there exist an open immersion $X\hookrightarrow \Xb$ such that $p_X:X\to \sS$ extends to a proper morphism $\pb_X: \Xb \to \sS$ such that
$p_X^{-1}(0_\sS)$ is dense in $\pb_X^{-1}(0_\sS)$.
\end{enumerate}

By Lemma \ref{lem;Levine} those $(X,0_X)$ satisfying $(*)$ form a cofinal system
of Nisnevich neighbourhoods of $(\A^1_\sS,0)$.
\medbreak

Put $\sX\ttimes_\sS \sX=\Spec K\{y_1\dots,y_{d-1},t_1,t_2\}$ and let 
$\Delta_\sX \hookrightarrow \sX\ttimes_\sS \sX$ be the closed immersion given by
the ideal $(t_1-t_2)$. 
For a Nisnevich neighbourhood $(X,0_X)$ of $(\A^1_\sS,0)$, put
\[D=\Delta_{\A^1_\sS}\times_{\A^1_\sS\times_\sS\A^1_\sS} (X\times_\sS X),\quad
E=\Delta_X\times_{X\times_{\sS}X}(\sX\ttimes_\sS \sX).\]
Note that $D$ is a Cartier divisor on $X\times_\sS X$.
We have a commutative diagram
\begin{equation}\label{eq;proof-lem3;purity}
\xymatrix{
E \ar[r] \ar[d]& \Delta_{\sX} \ar[r]\ar[d] \ar[ld]_{\hskip 20pt \psi_X} & \sX\ttimes_\sS\sX \ar[d]\\
\Delta_X \ar[r]\ar[rd]_{\beta} & D \ar[r]\ar[d]^{\alpha} & X\times_\sS X \ar[d]^{\Phi}\\
& \Delta_{\A^1_\sS} \ar[r] & \A^1_\sS\times_\sS\A^1_\sS \\
}\end{equation}
where all squares are Cartesian. 
All horizontal maps are closed immersions, and $\alpha$ and $\beta$ are \'etale.
Hence $\Delta_X\hookrightarrow D$ is a closed and open immersion so that
$D=\Delta_X \coprod \Xi_X$ for a closed and open immersion $\Xi_X\hookrightarrow D$.
Then $E\hookrightarrow \Delta_\sX$ is also a closed and open immersion so that $E=\Delta_\sX$ since $\Delta_\sX$ is irreducible. Thus we must have
\begin{equation}\label{eq0;claim1;purity}
\Xi_X \times_{X\times_\sS X} (\sX\ttimes_\sS\sX) =\emptyset.
\end{equation}

\begin{lemma}\label{claim1;purity}
Assume that $(X,0_X)$ satisfies $(*)$.
Let $T\subset \sS$ be an integral closed subscheme and
$W,W'\subset X$ be integral closed subschemes which are generically finite over $T$ and contain $0_X$.
Assume that $W$ does not contain $X_{0}=p_X^{-1}(0_\sS)$.
Then $\Xi_X$ does not contain an irreducible component of 
$W\times_{\sS} W'$ which is dominant over $T$.
\end{lemma}
\begin{proof}
Assume $\Xi_X$ contains an irreducible component $\sC$ of $W\times_{T} W'$ 
dominant over $T$. We claim $(0_X,0_X)\in \sC$. Then 
$\sC\times_{X\times_\sS X} (\sX\ttimes_\sS\sX) \not =\emptyset$, which contradicts 
\eqref{eq0;claim1;purity}. 
Let $W_{0}\subset W$ and $W'_{0}\subset W'$ be the fibres over $0_\sS\in \sS$.
By the assumption $(*)$ and Lemma \ref{lemA3}, $|W_{0}|=0_X$ and $W$ is 
finite over $\sS$, and $|W'_{0}|$ is either $0_X$ or $X_{0}$.
This implies that the projection $W\times_\sS W'\to W'$ is finite so that
$\sC\to W'$ is finite and surjective since it is dominant by the assumption. 
This implies $\sC_{0} \to W'_{0}$ is surjective. The last map factor as 
\[\sC_{0}\to (W\times_{\sS} W')_{0}= 0_X\times_\sS W'_{0} \to W'_{0}\]
where the second map is the projection. Hence we must have $(0_X,0_X)\in \sC$.
\end{proof}

\medbreak

\begin{lemma}\label{lem3;purity}
Assume $K$ is infinite. 
Let $T\subset \sS$ be an integral closed subscheme and $X_T=X\times_\sS T$.
Let $W\subset X_T$ be a reduced closed subscheme such that every irreducible component
of $W$ is generically finite and surjective over $T$ and contains $0_X$ but does not contain $p_X^{-1}(0_\sS)$.
Then there is a Nisnevich neighbourhood $(\tX,0_{\tX})$ of $(X,0_X)$ such that 
$(\tX, n \tW)$ is a quasi-$V$-pair over $T$ for all $n\in \Z_{>0}$, where
$\tW=W\times_X\tX$ (cf. Remark \ref{rem;def;Vtriple}(3)).
\end{lemma}
\begin{proof}
Note that if $(X',0_{X'})$ is a Nisnevich neighbourhood $(\tX,0_{\tX})$ of $(X,0_X)$, then
$W'=W\times_X X'$ satisfies the same condition as $W$.  
By Lemma \ref{lem;Levine} we may assume that there exists 
a closed immersion $X\hookrightarrow \A^N_\sS$ over $\sS$ such that
letting $\Xb$ be the normalization of the closure of $X$ in $\P^N_\sS$, 
$\Xinf=\Xb-X$ is finite over $\sS$. Let $\pb_X:\Xb\to \sS$ be the projection.
Noting $0_X\notin \Xinf$, Lemmas \ref{lemA3} (applied to the irreducible components of $W$) imply that $W$ is closed in $\Xb$ and every irreducible component of $W$ is finite and surjective over $\sS$ and $W\cap \pb_X^{-1}(0_\sS)=\{0_X\}$. Let $\Xi_{X}\subset X\times_\sS X$ be as \eqref{eq0;claim1;purity}.

\begin{claim}\label{claim1;lem3;purity}
Every irreducible component of $\Xi_X \cap (X\times_\sS W)$ is dominant and generically finite over $T$. 
\end{claim}
\begin{proof}
Let $W_\sS\subset \A^1_\sS$ be the image of $W$ under $X\to \A^1_\sS$.
Since $W$ is finite over $\sS$, $W_\sS$ is closed in $\A^1_\sS$ and every irreducible component of $W_\sS$ is finite and surjective over $\sS$.
By Lemma \ref{lemA3}, every irreducible component of $W_{\sS}\times_{\A^1_\sS} X$ which intersects $W$ must contain $0_X$ and be finite over $\sS$. Noting that any irreducible component of $\Xi_X \cap (X\times_\sS W)$ is an irreducible component of $\Xi_X \cap (X\times_\sS W_i)$ for some irreducible component $W_i$ of $W$, we may replace $W$ by the union of those irreducible components of $W_{\sS}\times_{\A^1_\sS} X$ that intersect with $W$.
Then the natural inclusion $W\hookrightarrow W_\sS\times_{\A^1_\sS} X$ is a closed and open immersion. Hence we get a closed and open immersion
\[X\times_\sS W \hookrightarrow X\times_\sS W_\sS\times_{\A^1_\sS} X\subset X\times_\sS X.\]
Note
\[ D=\Phi^{-1}(\Delta_{\A^1_\sS})\qaq 
X\times_\sS W_\sS\times_{\A^1_\sS} X=\Phi^{-1}(\A^1_\sS\times_\sS W_\sS),\]
where $\Phi:X\times_\sS X \to \A^1_\sS\times_\sS \A^1_\sS$ is from \eqref{eq;proof-lem3;purity}. Since
\[ \Delta_{\A^1_\sS}\cap (\A^1_\sS\times_\sS W_\sS)=\Delta_{W_\sS}
\subset W_\sS\times_\sS W_\sS\subset \A^1_\sS\times_\sS\A^1_\sS,\]
we get
\[ D\cap \big(X\times_\sS W_\sS\times_{\A^1_\sS} X)\big)=\Phi^{-1}(\Delta_{W_\sS}).\]
In view of the closed and open immersion $\Xi_X\hookrightarrow D$,
we conclude that $\Xi_X \cap (X\times_\sS W)$ is closed and open in $\Phi^{-1}(\Delta_{W_\sS})$. 
Since $\Phi$ is \'etale and every irreducible component of $W_\sS$ is finite and surjective over $T$, every irreducible component of $\Phi^{-1}(\Delta_{W_\sS})$ is dominant and generically finite over $T$. This proves the desired claim.
\end{proof}

Put 
\[  \Sigma=pr\big(\Xi_{X} \cap (X\times_\sS W)\big)\subset X_T,\]
where $pr: X\times_\sS W=X_T\times_T W \to X_T$ is the projection.
Since $W$ is finite and surjective over $T$, $pr$ is finite so that $\Sigma$ is closed in $X$.

\begin{claim}\label{claim2;lem3;purity}
We have $0_X\not\in \Sigma$.
\end{claim}
\begin{proof}
Assume there is an irreducible component $C$ of $\Sigma$ containing $0_X$.
Since $pr$ is finite and surjective, there is an irreducible component $\tC$ of $\Xi_X \cap (X_T\times_T W)$ such that $C=pr(\tC)$ and $\tC$ is an irreducible component of $C\times_T W$.
Since $\tC$ is dominant and generically finite over $T$ by Claim \ref{claim1;lem3;purity}, so must be $C$. This contradicts Lemma \ref{claim1;purity}.
\end{proof}

Let $\ol{\Sigma}\subset \Xb$ be the closure of $\Sigma$. 
Noting $\ol{\Sigma}\cap X =\Sigma$,
Claim \ref{claim2;lem3;purity} implies $0_X\not\in \ol{\Sigma}$. 
By \cite[Th. 5.1]{gll}, there exists hypersurfaces $H\subset\P^N_\sS$ and $H_{\inf}\subset\P^N_\sS$ over $\sS$ satisfying the following conditions:
\begin{itemize}
\item[$(i)$]
$H_{\inf}\supset \Xinf\cup \ol{\Sigma}$ and $0_X\not\in H_{\inf}$.
\item[$(ii)$]
$H\cap \big((H_{\inf}\cup W) \cap \pb_X^{-1}(0_\sS)\big) =\emptyset$.
\end{itemize}
Note $0_X\not\in H_{\inf}$ implies $H_{\inf}\cap \pb_X^{-1}(0_\sS)$ is finite.
Recall $W\cap \pb_X^{-1}(0_\sS)=\{0_X\}$. By Lemma \ref{lemA1}, $H_{\inf}$ is finite over $\sS$ and we have
\[(\Xb\cap H_{\inf})\cap W=\emptyset\qaq
(\Xb\cap H) \cap (H_{\inf}\cup W) =\emptyset.\]
Hence $(\Xb\cap H_{\inf})\cup W$ is contained in the affine scheme 
$\Xb\backslash H$. Putting 
\[\tX=X\backslash H_{\inf}=\Xb\backslash (\Xinf\cup H_{\inf}),\]
this implies $W\subset \tX$ and that $(\Xb,W+ (\Xinf + H_{\inf}))$ is a good compactification of $(\tX,W)$ over $\sS$.
By the definition of $\Xi_{X}$, $\Delta_{X}$ is principal on 
$(X\times_\sS X)\backslash \Xi_{X}$. 
Since $\big(\Xi_X\cap (\tX\times_\sS W)\big)=\emptyset$ by the construction,
this implies that the diagonal $n W\hookrightarrow \tX\times_\sS n W $ is principal for all $n\in \Z_{>0}$.
This completes the proof of Lemma \ref{lem3;purity}.
\end{proof}

\begin{rk}\label{rem;lem3;purity}
Let the assumption be as in Lemma \ref{lem3;purity}.
Assume further that $W$ is generically \'etale over a non-empty open $U\subset T$
such that $U\in \tSm$. Then $(\tX_U, \tW_U)$ is a nice $V$-pair over $U$,
where $\tX_U=\tX \times_\sS U$ and $\tW_U=W\times_X \tX \times_\sS U$. 
\end{rk}

\section{Gysin maps}\label{Gysin}

\begin{lemma}\label{lem4;purity}
Let $\psi_t: \sX \to \A^1_\sS=\sS[t]$ be as in Definition \ref{goodfibration}.
Let $i: \sZ\hookrightarrow \sX$ be the closed immersion defined by $t=0$ and
$j:\sU=\sX-\sZ\hookrightarrow \sX$ be the open complement.
Let $\pi: \sS' \to \sS$ be an essentially smooth morphism and 
$i':\sZ'\to \sX'$ (resp. $j':\sU'\to \sX'$) be the base change of $i$ (resp. $j$)
by $\sS' \to \sS$.
Let $D\subset \sS$ (resp. $D'\subset \sS'$) be an effective Cartier divisor such that
such that $|D|$ (resp. $|D'|$) is a simple normal crossing divisor on $S$ (resp. $S'$)
and $D'\geq \pi^* D$.
Write $D'_\sT=D'\times_{\sS'} \sT$ for a $\sS'$-scheme $\sT$.
Take $F\in \lsCItspNis$.
\begin{itemize}
\item[(1)]
There is an exact sequence of sheaves on $\sX'_{\Nis}$
\begin{equation*}\label{eq3;lem4;purity}
 0\to F_{(\sX',D'_{\sX'})} \to F_{(\sX',n\sZ'+D'_{\sX'})} \to 
i'_*(\Fcontt n)_{(\sZ',D'_{\sZ'})} \to 0,
\end{equation*}
which is natural in $(\sS',D') \to (\sS,D)$. For for $m\geq n$, the diagram
\[\xymatrix{
 0\ar[r]& F_{(\sX',D'_{\sX'})} \ar[r]\ar[d]^{=} & F_{(\sX',n\sZ'+D'_{\sX'})} 
\ar[r]\ar[d] &  i'_*(\Fcontt n)_{(\sZ',D'_{\sZ'})} \ar[r] \ar[d] &0\\
 0\ar[r]& F_{(\sX',D'_{\sX'})} \ar[r] & F_{(\sX',m\sZ'+D'_{\sX'})} \ar[r] 
& i'_*(\Fcontt m)_{(\sZ',D'_{\sZ'})} \ar[r] &0,
}\]
is commutative, where the vertical maps are the natural ones.
\item[(2)]
There is an exact sequence of sheaves on $\sX'_{\Nis}$
\[ 0\to F_{(\sX',D'_{\sX'})} \to j'_* F_{(\sU',D'_{\sU'})} \to 
i'_*(\Fcont)_{(\sZ',D'_{\sZ'})} \to 0,\]
and an isomorphism of sheaves on $\sZ'_{\Nis}$:
\[ \theta_{\psi}: (\Fc)_{(\sZ',D'_{\sZ'})} \isom R^1{i'}^!F_{(\sX',D'_{\sX'})},\]
which are natural in $(\sS',D') \to (\sS,D)$.
\end{itemize}
\end{lemma}
\begin{proof}
The exact sequence of (2) follows from that of (1) by taking the colimit over $n$,
where one uses isomorphisms
\[ j'_* F_{(\sU',D'_{\sU'})} \simeq \indlim n F_{(\sX',n\sZ'+D'_{\sX'})}\qaq
(\Fcont)_{(\sZ',D'_{\sZ'})} \simeq \indlim n(\Fcontt n)_{(\sZ',D'_{\sZ'})}\]
coming from Lemma \ref{lem1;Mrec}(1) (Note that $\Fc$ and $\Fcontt n$ have $M$-reciprocity by Lemma \ref{lem1;contraction}). 
The isomorphism of (2) follows from the exact sequence of (2) in view of an exact sequence
\begin{equation}\label{eq2;lem4;purity}
 0 \to F_{(\sX',D'_{\sX'})} \to {j'}_*F_{(\sU',D'_{\sU'})} \to 
{i'}_* R^1{i'}^!F_{(\sX',D'_{\sX'})}\to 0,
\end{equation}
where the injectivity of the first map comes from Theorem \ref{thm;localinjectivity}. 

It remains to show the exact sequence in (1).
Since $F_{(\sX',n\sZ'+D'_{\sX'})}/F_{(\sX',D'_{\sX'})}$ is supported on $\sZ'$,
it suffices to construct an isomorphism of sheaves
\begin{equation}\label{eq4;lem4;purity}
(\Fcontt n)_{(\sZ',D'_{\sZ'})} \isom 
{i'}^*\big(F_{(\sX',n\sZ'+D'_{\sX'})}/F_{(\sX',D'_{\sX'})}\big).
\end{equation}
which is natural in $(\sS',D')$. By Lemma \ref{lem3;purity},
there exists a cofinal system of Nisnevich neighbourhoods $(X,\ul0)$ of 
$(\A^1_\sS=\sS[t],\ul0)$ 
such that $(X,\sZ_X)$ are all nice $V$-pairs over $\sS$,
where $\sZ_X\subset X$ is the closed subscheme defined by $t=0$.
Put $X'=X\times_\sS\sS'$ and $\sZ'_X=\sZ_X\times_\sS\sS'$.
For $U\to \sS'$ \'etale, these provide nice $V$-pairs over $U$:
\[(X'\times_{\sS'} U,\sZ'_X\times_{\sS'} U)\]
by the base change. 
Note that $\sX'\to X'\to \sS'$ induce isomorphisms $\sZ'\simeq \sZ'_{X}\simeq \sS'$ and
that the formal completions of $X'$ along $\sZ'_{X}$ are naturally isomorphic to 
$\sZ'\times_K \Spf K[[t]]$. Thus Theorem \ref{thm;contraction}(1) implies an exact sequence
\[0\to F((X',D'_{X'}) \times_{\sS'} U) \to F\big((X',n\sZ'_{X}+D'_{X'})\times_{\sS'} U\big)\to \Fcontt n ((\sZ'_X,D'_{\sZ'_X})\times_{\sS'} U)\to 0.\]
By taking the colimit over $X$, we get an exact sequence
\begin{equation}\label{eq4-1;lem4;purity}
0\to F((\sX',D'_{\sX'}) \times_{\sS'} U) \to F\big((\sX',n\sZ'+D'_{\sX'})\times_{\sS'} U\big) \to \Fcontt n((\sZ',D'_{\sZ'}) \times_{\sS'} U)\to 0,
\end{equation}
which is natural in $U\to \sS'$ by Remark \ref{rm;thm;contraction}(1).
It gives an isomorphism $\gamma: (\Fcontt n)_{(\sZ',D'_{\sZ'})} \simeq \Theta$ of presheaves on $\sZ'_{\et}$, where $\Theta$ is defined by
\begin{equation*}\label{eq5;lem4;purity}
\Theta(U)=F\big((\sX',n\sZ'+D'_{\sX'})\times_{\sS'} U\big)/F((\sX',D'_{\sX'}) \times_{\sS'} U).
\end{equation*}
Since $U\to \sZ'\rmapo{i'} \sX'$ factors naturally as
$U\to \sX'\times_{\sS'} U\to \sX'$,
we have a natural map 
$\Theta \to {i'}^*\big(F_{(\sX',n\sZ'+D'_{\sX'})}/F_{(\sX',D'_{\sX'})}\big)$
of presheaves on $\sZ'_{\et}$.
Composing this with $\gamma$, we get the desired map \eqref{eq4;lem4;purity}. 
Its naturality in $(\sS',D') \to (\sS,D)$ follows from that of \eqref{eq4-1;lem4;purity}.

We now show that \eqref{eq4;lem4;purity} is an isomorphism at stalks over 
a point $x\in \sZ'$. 
Put $\sT=\hen {\sS'} y$ where $y\in \sS'$ is the image of $x$.
We have a natural isomorphism $\hen{\sX'} x\simeq \sT\{t\}$ such that $\hen {\sZ'} x\subset \hen{\sX'} x$ is given by $\{t=0\}$, where $\sT\{t\}$ is the henselization of $\A^1_\sT$ at the origin $\ul0$ of the closed fiber. 
By the same argument as above, we get an exact sequence
\begin{equation*}\label{eq1;lem4;purity}
0\to F(\hen{\sX'} x,D'_{\hen{\sX'} x})\to 
F(\hen{\sX'} x,n\hen{\sZ'} x+ D'_{\hen{\sX'} x}) \to 
\Fcontt n (\hen {\sZ'} x,D'_{\hen {\sZ'} x})\to 0.
\end{equation*}
Since the first (resp. second) term is the stalk of ${i'}^*F_{(\sX',D'_{\sX'})}$
(resp. ${i'}^*F_{(\sX',n\sZ'+D'_{\sX'})}$) at $x$, this implies the desired assertion and hence the exact sequence in (1).
Its naturality in $\sS'\to \sS$ follows from that of \eqref{eq4;lem4;purity}
and the commutativity of the diagram in (1) is obvious.
This completes the proof of Lemma \ref{lem4;purity}.
\end{proof}

\subsection{Construction of Gysin maps}\label{subsecgtionGysin}

In this subsection we assume $k$ is perfect and infinite.
Take $W\in \Sm$ and $e\in \cd W n$ with $n\in \Z_{>0}$.
Let $\sX=\hen{W}{e}$ and $K=k(e)$ be the residue field of the closed point 
$\ul0_\sX\in \sX$. 
Choose a regular system of parameters $\tul=(t_1,\dots,t_n)$ and an isomorphism
from Lemma \ref{lem;etalepoint}:
\[\ep :\sX\simeq \Spec K\{\tul\} = \Spec K\{t_1,\dots,t_n\}.\]
For chosen $\tul$, $\ep$ is determined by a choice of map $\sX\to \Spec(K)$
inducing the identity on $\ul0_\sX$ (cf. the proof of Lemma \ref{lem;etalepoint}).
For $0<i<n$ put $\sS_i=\Spec K\{t_{i+1},\dots,t_n\}$.
By convention put $\sS_0=\sX$ and $\sS_n=\Spec(K)$. 
Let $\sZ_i=\{t_i=0\}\subset \sX$ and $\sZ_{[1,i]}=\sZ_1\cap \sZ_2\cap \cdots\sZ_i$
with the closed immersion $\iota_{[1,i]}:\sZ_{[1,i]}\to \sX$ .
We then have a diagram:
\begin{equation}\label{eq1;Gysin}
\xymatrix{
\sX \ar[d]^{\psi_1} & \ar[l]_{{\iota}_1}\ar[ld] \sZ_1\ar[d]^{\psi_2}  & \ar[l]_{{\iota}_2}\ar[ld] \sZ_{[1,2]}\ar[d]^{\psi_3} &\ar[l]\cdots& \ar[l]\sZ_{[1,n-1]}\ar[d]^{\psi_n}  & \ar[l]_{{\iota}_n}\sZ_{[1,n]}\ar[ld] \\
\sS_1 & \sS_2 & \sS_3 & \cdots & \sS_n \\
}\end{equation}
where $\psi_i$ is a fibration with coordinate $t_i$ and ${\iota}_i$ is 
the closed immersion defined by $t_i=0$ and all slanting arrow are isomorphisms.
\medbreak

Take $F\in \lsCItspNis$.
By Lemma \ref{lem4;purity}(2) we obtain a series of isomorphisms of Nisnevich sheaves:
\begin{equation}\label{eq2;Gysin}
\begin{aligned}
& \theta_{\psi_1}: (\Fc)_{\sZ_1} \isom R^1{\iota}^!_1 F_{\sX},\\
& \theta_{\psi_2}: (\Fcc)_{\sZ_{[1,2]}} \isom R^1{\iota}^!_2 (\Fc)_{\sZ_1},\\
& \hskip 30pt \cdots\cdots\cdots\cdots\cdots\cdots\\ 
& \theta_{\psi_n}: (\Fccc n)_{\sZ_{[1,n]}} \isom 
R^1{\iota}^!_n (F_{-(n-1)})_{\sZ_{[1,n-1]}}.\\
\end{aligned}
\end{equation}
\medbreak

For an integer $n>0$ consider the following condition:
\begin{enumerate}
\item[$(\clubsuit)_n$]
For any $X\in \Sm$ and $x\in \cd X q$ with $q\leq n$, and 
any regular divisor $\sH\subset\sX:=\hen{X} {x}$ with $\sU=\sX-\sH$, and 
for $F\in \lsCItspNis$, 
\begin{equation}\label{eq-1;Gysin}
H^i(\sU_{\Nis}, F_\sU)=0 \qfor i>0.
\end{equation}
\end{enumerate}

Noting $H^0_\sH(\sX_{\Nis}, F_\sX)=0$ by Theorem \ref{thm;localinjectivity}, \eqref{eq-1;Gysin} is equivalent to 
\begin{equation}\label{eq0;Gysin}
H^i_\sH(\sX_{\Nis}, F_\sX)=0 \qfor i\not=1.
\end{equation}
$(\clubsuit)_n$ implies that for any closed immersion $\iota:Z\hookrightarrow X$ of codimension $q\leq n$ in $\Sm$, we have
\begin{equation}\label{eq0;puritycontraction}
R^i\iota^! F_X=0 \qfor i\not=q.
\end{equation}
Indeed, looking at the stalks, it suffices to show \eqref{eq0;puritycontraction}
in case $X$ is equal to $\sX$ in $(\clubsuit)_n$.
Then one can take
\[\sX \lmapo {\iota_1} \sZ_1 \lmapo{\iota_2} \sZ_{[1,2]}\gets \cdots
\lmapo{\iota_q}\sZ_{[1,q]} \]
as in \eqref{eq1;Gysin} with $Z=\sZ_{[1,q]}$.
By \eqref{eq0;Gysin} we get
\begin{equation*}\label{eq0a;puritycontraction}
R^\nu\iota^!_i F_{\sZ_{[1,i-1]}} =0\qfor\nu\not= 1\text{ and } 1\leq i\leq q.
\end{equation*}
\eqref{eq0;puritycontraction} then follows from this and Lemma \ref{lem4;purity}(2) since $(\clubsuit)_n$ applies also to $\Fc$ thanks to 
Lemma \ref{lem1;contraction}. It also implies the isomorphism
\begin{equation}\label{eq0-1;Gysin}
 R^q{\iota}_\sZ^!F_{\sX} \simeq R^1{\iota}_q^! \cdots R^1{\iota}_2^! R^1{\iota}_1^! F_{\sX},
\end{equation}
where $\iota_\sZ:\sZ\hookrightarrow \sX$ is the closed immersion.
\medbreak

We now return to the set-up in the beginning of \S\ref{subsecgtionGysin} and
assume $(\clubsuit)_n$.
By \eqref{eq0-1;Gysin}, the isomorphisms $\theta_{\psi_\nu}$ for $\nu=1,\dots,i$ induce isomorphisms 
\begin{equation}\label{eq;puritycontraction}
\theta^i_{\epsilon}: (\Fccc i)_{\sZ_{[1,i]}} \isom R^i{\iota}_{[1,i]}^!F_{\sX},
\end{equation}
depending on $\epsilon:\sX\simeq \Spec K\{t_1,\dots,t_n\}$. In particular we get
\begin{equation}\label{eq1;puritycontraction}
\theta_{\epsilon}=\theta^n_{\epsilon}: (\Fccc n)_{\ul0_\sX} \isom R^n \iota^!F_{\sX},
\end{equation}
where $\iota:\ul0_\sX \simeq \sZ_{[1,n]}\hookrightarrow \sX$.
We have the following commutative diagram
\begin{equation}\label{eq2;puritycontraction}
\xymatrix{
\Fccc i ( {\sZ}_{[1,i]}-{\sZ}_{[1,i+1]}) \ar[r]^{\partial} 
\ar[d]^{\theta^i_{\epsilon}}& \Fccc {(i+1)} ( {\sZ}_{[1,i+1]}) 
\ar[d]^{\theta^{i+1}_{\epsilon}} \\
H^i_{{\sZ}_{[1,i]}-{\sZ}_{[1,i+1]}}({\sX}_{\Nis},F_{\sX}) \ar[r]^{\delta} & 
H^{i+1}_{{\sZ}_{[1,i+1]}}({\sX}_{\Nis},F_{\sX}),\\}
\end{equation}
where the vertical maps are induced by the maps \eqref{eq;puritycontraction}
in view of \eqref{eq0;puritycontraction}, $\delta$ is a boundary map in the localization sequence for 
${\sZ}_{[1,i+1]}\hookrightarrow {\sZ}_{[1,i]}$, and $\partial$ is the composite
\begin{multline}\label{eq2-1;puritycontraction}
 \Fccc i ( {\sZ}_{[1,i]}-{\sZ}_{[1,i+1]}) \rmapo{\delta} 
H^1_{{\sZ}_{[1,i+1]}}({\sZ}_{[1,i]},(\Fccc i)_{{\sZ}_{[1,i]}}) 
\rmapo{\theta_{\psi_{i+1}}^{-1}} \Fccc {(i+1)}({\sZ}_{[1,i+1]}),\\
\end{multline}
where the second map is the inverse of the isomorphism from \eqref{eq2;Gysin}
in view of \eqref{eq0;puritycontraction}.
The commutativity of \eqref{eq2;puritycontraction} follows from the 
(obvious) commutativity of the diagram
\begin{equation}\label{eq2-3;puritycontraction}
\xymatrix{
H^0({\sZ}_{[1,i]}-{\sZ}_{[1,i+1]},(\Fccc i)_{{\sZ}_{[1,i]}}) \ar[r]^{\delta} 
\ar[d]^{\theta^i_\ep}& H^1_{{\sZ}_{[1,i+1]}}({\sZ}_{[1,i]},(\Fccc i)_{{\sZ}_{[1,i]}})) 
\ar[d]^{\theta^i_\ep} \\
H^0({\sZ}_{[1,i]}-{\sZ}_{[1,i+1]},R^i{\iota}_{[1,i]}^!F_{\sX}) \ar[r]^{\delta} 
& H^1_{{\sZ}_{[1,i+1]}}({\sZ}_{[1,i]},R^i{\iota}_{[1,i]}^!F_{\sX}) ,
\\}
\end{equation}
where the vertical maps are induced by the map \eqref{eq;puritycontraction}.
In fact, $\theta^{i+1}_\ep$ in \eqref{eq2;puritycontraction}
is the composite of $\theta_{\psi_{i+1}}$ in \eqref{eq2-1;puritycontraction}
and the map
\[ H^1_{{\sZ}_{[1,i+1]}}({\sZ}_{[1,i]},(\Fccc i)_{{\sZ}_{[1,i]}}) 
\rmapo{\theta^i_\ep} 
H^1_{{\sZ}_{[1,i+1]}}({\sZ}_{[1,i]},R^i{\iota}_{[1,i]}^!F_{\sX})\overset{\alpha}{\simeq}
H^{i+1}_{{\sZ}_{[1,i+1]}}({\sX}_{\Nis},F_{\sX}),\]
where the first map is the one from \eqref{eq2-3;puritycontraction} and the isomorphism $\alpha$ comes from \eqref{eq0;puritycontraction}.
Hence $\theta^{i+1}_\ep\circ\partial$ (resp. $\delta\circ\theta^i_\ep$) in \eqref{eq2;puritycontraction} is identified via $\alpha$ with $\theta^i_\ep\circ\delta$ (resp. $\delta\circ\theta^i_\ep$) in \eqref{eq2-3;puritycontraction}.

\medbreak

Let $\xi_{[1,i]}$ be the generic point of ${\sZ}_{[1,i]}$ and
$\hen {\sX} {\xi_{[1,i]}}$ be the henselization of $\sX$ at $\xi_{[1,i]}$
with the closed immersion $\iota_{\xi_{[1,i]}}: \xi_{[1,i]}\to 
\hen {\sX} {\xi_{[1,i]}}$.
We have the following commutative diagram
\begin{equation}\label{eq2-2;puritycontraction}
\xymatrix{
(\Fccc i)_{\xi_{[1,i]} } \ar[r]^{\hskip -20pt\theta^i_{\epsilon}}
 \ar[rd]_{\theta^i_{\ep_{\xi_{[1,i]}}}} 
&(R^i{\iota}_{[1,i]}^!F_{\sX})_{\xi_{[1,i]}} \ar[d]^{\simeq}\\
& R^i{\iota}_{\xi_{[1,i]}}^!F_{\hen {\sX} {\xi_{[1,i]}} }
 \\
}\end{equation}
where $\ep_{\xi_{[1,i]}}: \hen {\sX} {\xi_{[1,i]}}\simeq \Spec k(\xi_{[1,i]})\{t_1,\dots,t_i\}$ is the isomorphism induced by $\ep$, and
$\theta^i_{\ep_{\xi_{[1,i]}}}$ is \eqref{eq1;puritycontraction} for 
$(\sX,\ep)=(\hen {\sX} {\xi_{[1,i]}},\ep_{\xi_{[1,i]}})$.
The commutativity can be checked by using Remark \ref{rm;thm;contraction}.
\medbreak

We will use the following variant of \eqref{eq1;puritycontraction}. 
Let $\sS' \to \Spec(K)$ be an essentially smooth morphism and let
\begin{equation*}
\xymatrix{
\sX' \ar[d]^{\psi'_1} & \ar[l]_{{\iota'}_1}\ar[ld] \sZ'_1\ar[d]^{\psi'_2}  & \ar[l]_{{\iota'}_2}\ar[ld] \sZ'_{[1,2]}\ar[d]^{\psi'_3} &\ar[l]\cdots& \ar[l]\sZ'_{[1,n-1]}\ar[d]^{\psi'_n}  & \ar[l]_{{\iota'}_n}\sZ'_{[1,n]}\ar[ld] \\
\sS'_1 & \sS'_2 & \sS'_3 & \cdots & \sS'_n \\
}\end{equation*}
be the base changes of \eqref{eq1;Gysin} by $\sS' \to \Spec(K)$. 
By Lemma \ref{lem4;purity}(2) we obtain a series of isomorphisms as in \eqref{eq2;Gysin}:
\begin{equation}\label{eq3;Gysin}
\begin{aligned}
& \theta_{\psi_1}: (\Fc)_{\sZ'_1} \isom R^1{\iota'}^!_1 F_{\sX'},\\
& \theta_{\psi_2}: (\Fcc)_{\sZ'_{[1,2]}} \isom R^1{\iota'}^!_2 (\Fc)_{\sZ'_1},\\
& \hskip 30pt \cdots\cdots\cdots\cdots\cdots\cdots\\ 
& \theta_{\psi_n}: (\Fccc n)_{\sZ'_{[1,n]}} \isom 
R^1{\iota'}^!_n (F_{-(n-1)})_{\sZ'_{[1,n-1]}}.\\
\end{aligned}
\end{equation}
Assume now 
\begin{equation}\label{eq3.5;Gysin}
\text{$(\clubsuit)_{\dim(\sX)+\dim(\sS')}$ holds.}
\end{equation}
The same argument as above gives an isomorphism
\begin{equation}\label{eq4;puritycontraction}
\theta_{\epsilon,\sS'}: (\Fccc n)_{\sS'} \isom R^n {\iota'}^!F_{\sX'},
\end{equation}
where $\iota'=\iota\times_K \sS'$ for $\iota:\ul0_\sX \hookrightarrow \sX$.
Let $\eta$ be the generic point of $\sS'$ which gives rise to the point
$(\ul0_\sX,\eta)\in \sX'=\sX\times_K \sS'$ denoted also by $\eta$.
Let $\hen {\sX'} {\eta}$ be the henselization of $\sX'$ at $\eta$
with the closed immersion $\iota_\eta: \eta \to \hen {\sX'} {\eta}$.
We have the following commutative diagram
\begin{equation}\label{eq4-2;puritycontraction}
\xymatrix{
(\Fccc n)_{\eta } \ar[r]^{\hskip -20pt\theta_{\epsilon,\sS'}}
 \ar[rd]_{\theta_{\ep_\eta} }
&(R^n{\iota'}^!F_{\sX'})\ar[d]^{\simeq})_{\eta} \\
& R^n{\iota}_\eta^!F_{\hen {\sX'} {\eta} }
 \\
}\end{equation}
where $\ep_{\eta}: \hen {\sX'} {\eta} \simeq \Spec k(\eta)\{t_1,\dots,t_n\}$ is 
the isomorphism induced by $\ep$, and $\theta_{\ep_\eta}$ is \eqref{eq1;puritycontraction} for 
$(\sX,\ep)=(\hen {\sX'} {\eta},\ep_{\eta})$.
The commutativity can be checked again by using Remark \ref{rm;thm;contraction}.


\section{Vanishing theorem}\label{vanishing}

In this section we assume $k$ is perfect.

\begin{thm}\label{thm;vanishing}
Take $X\in \Sm$ and $x\in X$ and 
a regular divisor $\sH\subset\sX:=\hen{X} {x}$ with $\sU=\sX-\sH$. 
For $F\in \lsCItspNis$ we have
\begin{equation}\label{eq;thm;vanishing}
H^i(\sU_{\Nis}, F_\sU)=0 \qfor i>0.
\end{equation}
\end{thm}

\medbreak

Let $\sU$ be a noetherian scheme and $F$ be a sheaf on $U_{Nis}$. Put 
\[ \sCFU q = \underset{x\in \cd \sU q}{\bigoplus}\; H^q_x(\sU_{Nis},F)
\qfor q\in\Z_{\geq 0}\]
and let $\sCFU \bullet$ denote the Cousin complex:
\begin{equation}\label{eq;cousincomplex}
\sCFU 0 \rmapo{\partial^o} \sCFU 1 \rmapo{\partial^1}  \cdots \to 
\sCFU q \rmapo{\partial^q}\cdots  ,
\end{equation}
where $\partial^q$ is the boundary map arising from localization theory.

\medbreak\noindent
{\it Proof of Theorem \ref{thm;vanishing}}\;
By a standard norm argument we may assume $k$ is infinite.
We prove $(\clubsuit)_d$ from \S\ref{Gysin} by induction on $d$.
Assume $(\clubsuit)_{d-1}$ holds. Take $\sX=\Spec K\{x_1,\dots,x_d\}$ and a regular divisor $\sH\subset\sX$ with $\sU=\sX-\sH$ (cf. Lemma \ref{lem;etalepoint}), and 
$F\in \lsCItspNis$.
We write $\sCFU \bullet$ for $\sC_{\sU}^\bullet(F_\sU)$, where
$F_\sU$ is the sheaf on $\sU_{\Nis}$ induced by $F$.
Noting $\dim(\sU)=d-1$, the induction hypothesis and \eqref{eq0;puritycontraction} imply
\begin{equation}\label{eq0;proofvanishing}
H^i_x(\sX_{\Nis},F_\sX)=0 \;\;\text{for $x\in \cd \sU q$ and $i\not=q$}.
\end{equation}
Then \eqref{eq;thm;vanishing} follows from the acyclicity: 
\begin{equation}\label{eq0.1;proofvanishing}
H^q(\sCFU \bullet)=0 \qfor q>0.
\end{equation}

\medbreak
In what follows we fix $q$ with $1\leq q\leq d-1$.
Take successive fibrations with coordinate:
\begin{equation}\label{eq1;proofvanishing}
\xymatrix{ 
\sX \ar[r]^{\psi_1} & \sS_1 \ar[r]^{\psi_2} & \sS_2 \ar[r] & \cdots 
\ar[r]^{\psi_{d-q}} & \sS_{d-q}\\
\sH_0 \ar[r] \ar[u]_{\hookrightarrow}& \sH_1 \ar[r]\ar[u]_{\hookrightarrow} & 
\sH_2 \ar[r]\ar[u]_{\hookrightarrow} & \cdots \ar[r]
& \sH_{d-q} \ar[u]_{\hookrightarrow}\;. \\
}
\end{equation}
Here $\sH_0=\sH$ and $\sH_{i}\subset \sS_i$ is a regular divisor 
and $\psi_i$ is a $\sH_{i-1}$-fibration with coordinate $t_i$ 
such that $\sH_{i-1}=\psi_i^{-1}(\sH_i)$ for $1\leq i\leq d-q$.
Put $\sV_i=\sS_i-\sH_i$. 
Let $\phi_i: \sX \to \sS_i$ be the induced map.
For $w\in \sS_i$ let $\sX_w=\phi_i^{-1}(w)$ and $\xi_w$ be its generic point.

\begin{claim}\label{claim0;vanishing}
For $w\in \cd {\sV_i} q$ we have an exact sequence 
\begin{equation*}\label{eq1.0;proofvanishing}
0\to H^q_{\sX_w}(\sX_{Nis},F_\sX) \to 
H^q_{\xi_w}(\sX_{Nis},F_\sX) \rmapo{\partial} 
\underset{z\in \cd {\sX_w} 1}{\bigoplus}\;H^{q+1}_{z}(\sX_{Nis},F_\sX),
\end{equation*}
where $\partial$ is induced by $\partial^q$ in \eqref{eq;cousincomplex}.
\end{claim}
\begin{proof}
In view of the spectral sequence
\[ E_1^{a,b}=\underset{x\in \cd {\sX} a\cap \sX_w}{\bigoplus}\; H^{a+b}_x(\sX_{Nis},F_\sX) 
\Rightarrow H^{a+b}_{\sX_w}(\sX_{Nis},F_\sX),\]
the claim follows from \eqref{eq0;proofvanishing}.
\end{proof}

Put
\[\sDFUphii q i = \underset{w\in \cd {\sV_i} {q}}{\bigoplus}\;
H^q_{\sX_w}(\sX_{Nis},F_\sX).\]
By Claim \ref{claim0;vanishing}, $\partial^q$ in \eqref{eq;cousincomplex} induces
\begin{equation}\label{eq1.1;cousincomplex}
 \dphi q: \sDFUphii q i\to \sDFUphii {q+1} i,
\end{equation}
which give rise to a subcomplex $(\sDFUphii \bullet i,\dphi \bullet)$ of 
$(\sCFU \bullet,\partial^\bullet)$.
For $w\in \cd{\sV_{i}} q$ let $(\sS_{i-1})_w=\sS_{i-1}\times_{\sS_{i}} w$ 
and $\eta_w$ be its generic point and $\sX_{\eta_w}=\sX\times_{\sS_{i-1}} \eta_w$.
By definition we have
\[ \cd {(\sX_w)} 1 = \cd {(\sX_{\eta_w})} 1 
\cup \{\xi_y\;|\; y\in \cd {(\sS_{i-1})_w} 1\},\]
where $\xi_y$ is the generic point of $\sX_y=\sX\times_{\sS_{i-1}}y$, and 
$\sX_{\eta_w}$ and $\sX_w$ have the same generic point $\xi_w$.
By Claim \ref{claim0;vanishing} we have an exact sequence
\begin{equation}\label{eq1.2;cousincomplex}
0\to H^q_{\sX_w}(\sX_{Nis},F_\sX) \to
H^q_{\sX_{\eta_w}}(\sX_{Nis},F_\sX) \rmapo {\parphiw q} 
\underset{y\in \cd {(\sS_{i-1})_w} 1}{\bigoplus}\;H^{q+1}_{\sX_y}(\sX_{Nis},F_\sX),
\end{equation}
where $\parphiw q$ is induced by $\partial^q$ in \eqref{eq;cousincomplex}.
Putting
\[ \sDFUphih q i= \underset{w\in \cd {\sV_i} {q-1}}{\bigoplus}
\underset{y\in \cd {(\sS_{i-1})_w} 1}{\bigoplus}\; H^q_{\sX_y}(\sX_{Nis},F_\sX),\]
\[ \sDFUphiv q i = \underset{w\in \cd {\sV_i} {q}}{\bigoplus}\; 
H^q_{\sX_{\eta_w}}(\sX_{Nis},F_\sX),\]
we get a map
\begin{equation}\label{eq1.3;cousincomplex}
   \parphi q=\underset{w\in \cd {\sV_i} q}{\bigoplus}\parphiw q : 
\sDFUphiv q i \to  \sDFUphih  {q+1} i
\end{equation}
such that $\sDFUphii q i=\Ker(\parphi q)$. 
Noting
\[\cd {\sV_{i-1}} q = \underset{w\in \cd {\sV_i} {q-1}}{\cup} \cd {(\sS_{i-1})_w} 1
\;\cup\; \{\eta_w\;|\; w\in \cd {\sV_i} q\},\]
we have an exact sequence
\[ 0\to   \sDFUphiv q i  \rmapo{f_q} \sDFUphii q {i-1} \rmapo{g_q}  \sDFUphih q i\to 0\]
and $\parphi q$ is identified with $g_{q+1}\circ \delta_{\phi_{i-1}}^q \circ f_q$.
We have a commutative diagram
\begin{equation}\label{eq1.3.5;cousincomplex}
 \xymatrix{
\sDFUphii {q-1} i \ar[r]^{\hookrightarrow}  \ar[d]^{\delta^{q-1}_{\phi_i}}& 
\sDFUphiv {q-1} i \ar[r]^{ f_{q-1}}
\ar@/^9ex/[rrd]^{\partial^{q-1}_{\phi_i}}  & \sDFUphii {q-1} {i-1} 
\ar[d]^{\delta^{q-1}_{\phi_{i-1}}} \\ 
 \sDFUphii q i\ar[r]^{\hookrightarrow} &\sDFUphiv q i \ar[r]^{ f_q} 
\ar@/_9ex/[rrd]_{\partial^q_{\phi_i}} & \sDFUphii q {i-1} \ar[r]^{g_q} \ar[d]^{\delta^q_{\phi_{i-1}}}&  \sDFUphih q i \\
  && \sDFUphii {q+1} {i-1} \ar[r]^{g_{q+1}} &  \sDFUphih {q+1} i 
}\end{equation}
Thus there are inclusions of complexes: 
\begin{equation}\label{eq1.3.8;cousincomplex}
(\sDFUphii \bullet i,\dphi \bullet) \hookrightarrow (\sDFUphii \bullet {i-1},\dphii \bullet) \qfor i\geq1,
\end{equation}
which induce maps
\begin{equation}\label{eq1.4;cousincomplex}
  H^q(\sDFUphii \bullet {i}) \to H^q(\sDFUphii \bullet {i-1}),
\end{equation}
where $(\sDFUphii \bullet i,\dphi \bullet) =(\sCFU \bullet,\partial^\bullet)$ for
$i=0$ by convention.
An easy diagram chase on \eqref{eq1.3.5;cousincomplex} shows the following.
 
\begin{claim}\label{claim0-1;vanishing}
Let $\alpha\in \Ker(\sDFUphii q {i-1}\rmapo{\dphii q} \sDFUphii {q+1}{i-1} )$ and 
assume $g_q(\alpha)\in \Image(\parphi {q-1})$.
Then the class of $\alpha$ in $H^q(\sDFUphii \bullet {i-1})$ lies in the image of 
the map \eqref{eq1.4;cousincomplex}.
\end{claim}

Now the key point of the proof of \eqref{eq0.1;proofvanishing} is the following.

\begin{claim}\label{claim;vanishing}
Let $\alpha\in \sDFUphii q {i-1}$ and $\Lambda_\alpha\subset \cd {\sV_{i-1}} q$ be the finite subset such that the $w$-component of $\alpha$ is $0$ for 
$w\not\in \Lambda_\alpha$.
Assume that the closure $\ol{\Lambda}_\alpha$ of $\Lambda_\alpha$ in $\sS_{i-1}$ is strongly admissible for $\psi_{i}$ (cf. Definition \ref{goodfibration}(3)). 
Then $g_q(\alpha)\in \Image(\parphi {q-1})$.
\end{claim}

\medbreak

Admitting Claim \ref{claim;vanishing}, we finish the proof of \eqref{eq0.1;proofvanishing}. Take 
\[\alpha\in \Ker(\sCFU q\rmapo{\partial^q} \sCFU {q+1})\]
and let $\Lambda_\alpha\subset \cd {\sU} q$ be the finite subset such that the $v$-component of $\alpha$ is $0$ for $v\not\in \Lambda_\alpha$.
By Lemma \ref{lem;existencegoodfibration} we can choose $\psi_1(=\phi_1):\sX \to \sS_1$ in \eqref{eq1;proofvanishing} so that the closure $\ol{\Lambda}_\alpha$ of $\Lambda_\alpha$ in $\sX$ is strongly admissible for $\psi_1$. By Claims \ref{claim0-1;vanishing} and \ref{claim;vanishing}, there exists 
\[\alpha_1\in \Ker(\sDFUphii q 1\rmapo{\delta_{\phi_1}^q} \sDFUphii {q+1} 1)\] 
whose image under \eqref{eq1.3.8;cousincomplex} for $i=1$ and $\alpha$ have the same class in $H^q(\sCFU \bullet)$.
Applying the same argument to $\alpha_1$ in place of $\alpha$, we can choose $\psi_2:\sS_1\to \sS_2$ in such a way that there exists 
\[\alpha_2\in \Ker(\sDFUphii q 2\rmapo{\delta_{\phi_2}^q} \sDFUphii {q+1} 2)\] 
whose image under \eqref{eq1.3.8;cousincomplex} for $i=2$ and $\alpha_1$ have the same class in $H^q(\sDFUphii \bullet 1)$.
Repeating the same argument, this finishes the proof of \eqref{eq0.1;proofvanishing}
since $\cd {\sV_i} q=\emptyset$ and $\sDFUphii \bullet {i}=0$ for $i=d-q$.
\bigskip

In what follows we prove Claim \ref{claim;vanishing}.
We fix $i\geq 1$ and write $\sT=\sS_i$ and $\sS=\sS_{i-1}$, and
\[\psi=\psi_i:\sS\to \sT \qaq \phi=\phi_i:\sX\to \sT,\]
where $\sS=\sX$ if $i=1$ by convention. We also fix $w\in \cd \sT q\cap \sV_i$ and choose a regular system of parameters $\tauul_w=(\tau_{w,1},\dots,\tau_{w,q})$ of $\sT$ at $w$ and an isomorphism
\[\lambda_w: \sT^{h}_{|w} \simeq \Spec k(w)\{\tauul_w\}.\]
Let $\sS_w=\sS\times_\sT w$ and $\iota_w:\sX_w=\sX\times_\sS \sS_w \to \sX$ be the immersion.

\begin{lemma}\label{claim;henselianframe}
The map $\lambda_w$, $\phi$ and $\psi$ determine isomorphisms
\begin{equation}\label{eq6-2;cousincomplex}
\ep_{\phi,\lambda_w} : \sX^{h}_{|\sX_w} \simeq
\big(\sT^{h}_{|w} \times_{k(w)} \sX_w \big)^h_{|\sX_w}
\end{equation}
\begin{equation}\label{eq6-1;cousincomplex}
\ep_{\psi,\lambda_w} : \sS^{h}_{|\sS_w} \simeq 
\big(\sT^{h}_{|w} \times_{k(w)} \sS_w \big)^h_{|\sS_w}.
\end{equation}
\end{lemma}
\begin{proof}
We prove only the first isomorphism (the second is a special case of the first).
Let $\phi_w:\hen{\sX}{\sX_w}\to \hen{\sT}{w}$ be induced by $\phi:\sX\to\sT$.
It suffices to construct a map 
$\sigma: \hen{\sX}{\sX_w}\to \sX_w$ such that $\sigma\circ\iota_w=\id_{\sX_w}$ and that the diagram
\begin{equation}\label{eq1;claim;cousincomplex}
\xymatrix{
\hen{\sX}{\sX_w}\ar[r]^{\sigma}\ar[d]^{\phi_w}&  \sX_w\ar[d]\\
\hen \sT w\ar[r] & \Spec k(w)\\}\end{equation}
commutes, where $\hen{\sT}{w}$ is viewed as a scheme over $k(w)$ via $\lambda_w$. Indeed note that letting $\tul=(t_1,\dots,t_{i})$, $\phi_w$ factors as 
\[\hen{\sX}{\sX_w}\rmapo{\upsilon} \hen{\sT}{w}[\tul]\to \hen{\sT}{w},\]
where $\upsilon$ is induced by $\sS_{\nu-1}\to \sS_\nu[t_\nu]$ from \eqref{eq1;proofvanishing} for $1\leq \nu\leq i$ and it is essentially \'etale.
We get a map
\[ \gamma=(\sigma,\phi_w): \hen{\sX}{\sX_w}\to \sX_w\times_{k(w)} \hen{\sT}{w}=
\sX_w\times_{k(w)[\tul]} \hen{\sT}{w}[\tul].\]
The map is essentially \'etale since its composite with the projection
$\sX_w\times_{k(w)[\tul]} \hen{\sT}{w}[\tul] \to \hen{\sT}{w}[\tul]$ coincides with
$\upsilon$. Moreover the fibres over $w\in \hen{\sT}{w}$ of $\hen{\sX}{\sX_w}$ and 
$\sX_w\times_{k(w)} \hen{\sT}{w}$ are both equal to $\sX_w$ and $\gamma$ induces
the identity on it. Hence it induces the desired isomorphism $\ep_{\phi,\lambda_w}$. 

To construct $\sigma$, it suffices to construct a section 
$\rho$ of the projection 
\[pr: \hen{\sX}{\sX_w}\times_{k(w)} \sX_w\to \hen{\sX}{\sX_w}\]
such that the restriction of $\rho$ to the fibers over $w\in \hen{\sT}{w}$ is the diagonal $\Delta_{\sX_w}: \sX_w\to \sX_w\times_{k(w)}\sX_w$.
Here $\hen{\sX}{\sX_w}$ is viewed as a scheme over $k(w)$ via
\[ \hen{\sX}{\sX_w} \rmapo{\phi_w}\hen{\sT}{w}\rmapo{\lambda_w} 
\Spec k(w)\{\tauul_w\} \to \Spec k(w).\]
Then the composite $\sigma=pr\circ \rho$ satisfies
$\sigma\circ\iota_w=\id_{\sX_w}$ and the commutativity of \eqref{eq1;claim;cousincomplex}.

Let $\widehat{\sX}_{|\sX_w}$ (resp. $\widehat{\sT}_{|w}$) be the formal completion of $\sX$ along $\sX_w$ (resp. $\sT$ at $w$). Note that $\lambda_w$ induces an isomorphism
$\widehat{\sT}_{|w}\simeq \Spf k(w)[[\tauul_w]]$.
By \cite[Ch. 0, (19.7.1.5)]{EGA4}, there exists an isomorphism 
$ \widehat{\sX}_{|\sX_w} \simeq \widehat{\sT}_{|w} \times_{k(w)} \sX_w$
of formal schemes over $\widehat{\sT}_{|w}$. It gives rise to a map
$\hat{\sigma}:  \widehat{\sX}_{|\sX_w} \to \sX_w$ such that
$\hat{\sigma}\circ\iota_w=\id_{\sX_w}$.
Put
\[\hat{\rho}=(id_{\widehat{\sX}_{|\sX_w}},\hat{\sigma}) : \widehat{\sX}_{|\sX_w}\to \widehat{\sX}_{|\sX_w}\times_{k(w)} \sX_w.\]
By the construction its restriction to the fibers over $w\in \widehat{\sT}_{|w}$ is the diagonal $\Delta_{\sX_w}$, and $\hat{\rho}$ is a section of the projection
$\hat{pr}:\widehat{\sX}_{|\sX_w}\times_{k(w)} \sX_w \to \widehat{\sX}_{|\sX_w}$,
which is the base change via $\widehat{\sX}_{|\sX_w}\to \hen{\sX}{\sX_w}$ of $pr$.
Take a factorization $\sX_w\to X_w\to w$ of $\sX_w\to w$, where $X_w$ is of smooth of finite type over $k(w)$ and $\sX_w\to X_w$ is essentially \'etale.
Then $\hat{\rho}$ induces a section $\hat{\rho}_{X_w}$ of the projection
$\hat{pr}_{X_w}:\widehat{\sX}_{|\sX_w}\times_{k(w)} X_w \to \widehat{\sX}_{|\sX_w}$.
By \cite[II Th.2 bis]{elk} there exists a section
$\rho_{X_w}$ of the projection
$pr_{X_w}: \hen{\sX}{\sX_w}\times_{k(w)} X_w \to \hen{\sX}{\sX_w}$ such that
the restrictions of $\rho_{X_w}$ and $\hat{\rho}_{X_w}$  to the fibers over $w\in \widehat{\sT}_{|w}$ coincide. Let $\sigma_{X_w}:\hen{\sX}{\sX_w}\to X_w$ be the composite of $\rho_{X_w}$ and the projection $\hen{\sX}{\sX_w}\times_{k(w)} X_w\to X_w$. The map $\sX_w\to X_w$ induces
\[\pi : \hen{\sX}{\sX_w} \times_{X_w} \sX_w \to \hen{\sX}{\sX_w}, \]
where the fiber product is the base change of $\sX_w\to X_w$ via $\sigma_{X_w}$.
The composite of $\pi$ and $(\iota_w,id_{\sX_w}): \sX_w\to \hen{\sX}{\sX_w} \times_{X_w} \sX_w$ coincides with $\iota_w$, where $\iota_w:\sX_w\to \hen{\sX}{\sX_w}$
is the closed immersion. Since $\pi$ is essentially \'etale, there is an open and closed subscheme $Q\subset \hen{\sX}{\sX_w} \times_{X_w} \sX_w$ containing the image of $(\iota_w,id_{\sX_w})$ such that $\pi$ induces an isomorphism $Q\simeq \hen{\sX}{\sX_w}$. Thus we get a map 
\[\rho: \hen{\sX}{\sX_w} \simeq Q\hookrightarrow \hen{\sX}{\sX_w} \times_{X_w} \sX_w
\to \hen{\sX}{\sX_w} \times_{k(w)}\sX_w.\]
Then $\rho$ satisfies the desired property. 
This completes the proof.

\end{proof}

Using \eqref{eq6-2;cousincomplex} we apply \eqref{eq4;puritycontraction} to 
$(\sX,\ep,\sS')=(\hen{\sT} {w},\lambda_w,\sX_w)$ (note \eqref{eq3.5;Gysin} is satisfied since $\dim(\hen{\sT} {w})+\dim(\sX_w)=\dim(\hen {\sX} {\sX_w})\leq \dim(\sU)=d-1$). We then get an isomorphism of sheaves on $(\sX_w)_{\Nis}$:
\begin{equation}\label{eq6-3;cousincomplex}
 \theta_{\phi,\lambda_w}: (\Fccc q)_{\sX_w} \isom R^q{\iota_w}^! F_\sX
\end{equation}
as the composite map
\begin{multline*}
(\Fccc q)_{\sX_w} \underset{\eqref{eq4;puritycontraction}}
{\rmapo{\theta_{\lambda_w,\sX_w}}} R^q\iota_w^! F_{|\hen{\sT} w \times_w \sX_w}
\rmapo{exc} R^q\iota_w^! F_{|(\hen{\sT} w \times_w \sX_w)^h_{|\sX_w}}\\
\underset{\eqref{eq6-2;cousincomplex}}{\rmapo{\epsilon_{\phi,\lambda_w}}} R^q{\iota_w}^! F_{\hen{\sX} {\sX_w}}
 \rmapo{exc^{-1}} R^q{\iota_w}^! F_\sX
\end{multline*}
where $exc$ means excison isomorphisms.
By \eqref{eq0;puritycontraction} $R^\nu{\iota_w}^! F_\sX =0$ for $\nu\not=q$.
Hence \eqref{eq6-3;cousincomplex} induces isomorphisms
\begin{equation}\label{eq6-5;cousincomplex}
\begin{aligned}
&\theta_{\phi,\lambda_w}:\Fccc q(\sX_w)\simeq H^q_{\sX_w}(\sX_{Nis},F_\sX),\\
&\theta_{\phi,\lambda_w}:\Fccc q(\sX_{\eta_w})\simeq H^q_{\sX_{\eta_w}}(\sX_{Nis},F_\sX).\\
\end{aligned}
\end{equation}

For each $y\in \cd {\sS_w} 1$ choose a local parameter $\tau_y$ of $\sS$ at $y$ such 
that $(\tauul_w,\tau_y)$ forms a regular system of parameters of $\sS$ at $y$.
Choose also an isomorphism
\begin{equation}\label{eq6-5-1;cousincomplex}
\lambda_y:\hen{(\sS_w)} y \simeq \Spec k(y)\{\tau_y\}.
\end{equation}
\eqref{eq6-1;cousincomplex} and $\lambda_y$ determine an isomorphism
\begin{equation}\label{eq6-9;cousincomplex}
\ep_{\psi,\lambda_w,\lambda_y}: \hen{\sS} y \simeq \Spec k(y)\{\tauul_w,\tau_y\}.
\end{equation}
By the same argument as in the proof of Lemma \ref{claim;henselianframe},
$ \lambda_y$ and $\phi$ determine isomorphisms
\begin{equation}\label{eq6-6;cousincomplex}
\ep_{\phi,\lambda_y} : \hen {(\sX_w)} {\sX_y} \simeq
\big(\hen{(\sS_w)} y \times_{k(y)} \sX_y\big)^h_{|\sX_y}\overset{\lambda_y}{\simeq} 
\big(\Spec k(y)\{\tau_y\}\times_{k(y)} \sX_y\big)^h_{|\sX_y},
\end{equation}
\begin{equation}\label{eq6-8;cousincomplex}
\ep_{\phi,\lambda_w,\lambda_y} : \hen {\sX} {\sX_y} \simeq 
\big(\hen{\sS} y \times_{k(y)}\sX_y\big)^h_{|\sX_y}\overset{\eqref{eq6-9;cousincomplex}}{\simeq} 
\big(\Spec k(y)\{\tauul_w,\tau_y\}\times_{k(y)} \sX_y\big)^h_{|\sX_y}.
\end{equation}
which are compatible with \eqref{eq6-2;cousincomplex} in an obvious sense.

We have the following commutative diagram
\begin{equation}\label{eq7-1;cousincomplex}
\xymatrix{
\Fccc q(\sX_w) \ar[r]\ar[d]^{\theta_{\phi,\lambda_w}}_{\simeq} & 
\Fccc q(\sX_{\eta_w}) \ar[r]^{\hskip -30pt\partial^q}\ar[d]^{\theta_{\phi,\lambda_w}}_{\simeq}_{\simeq} & {\bigoplus}_{y\in \cd {\sS_w} 1}\;
\Fccc {(q+1)}(\sX_y) \ar[d]^{\theta_{\phi,\lambda_w,\lambda_y}}_{\simeq}_{\simeq}\\
H^q_{\sX_w}(\sX_{Nis},F_\sX) \ar[r]^{\gamma}&
H^q_{\sX_{\eta_w}}(\sX_{Nis},F_\sX) \ar[r]^{\hskip -20pt\partial^q_{\phi,w}} 
& {\bigoplus}_{y\in \cd {\sS_w} 1}\;
H^{q+1}_{\sX_y}(\sX_{Nis},F_\sX)\\
}\end{equation}
where $\partial^q_{\phi,w}$ is from \eqref{eq1.2;cousincomplex}, and 
$\theta_{\phi,\lambda_w,\lambda_y}$ is the isomorphism induced by $\theta_{\ep,\sS'}$ in \eqref{eq4;puritycontraction} for 
($\sX,\ep,\sS')=(\hen{\sS} {y},\ep_{\psi,\lambda_w,\lambda_y},\sX_y)$ using \eqref{eq6-8;cousincomplex}. 
Letting $\eta_{w,y}$ be the generic point of $\hen{(\sS_w)}{y}$,
$\partial^q$ is the sum over $y\in \cd {\sS_w} 1$ of the map 
\begin{multline}\label{eq7-1-1;cousincomplex}
\partial^q_{w,y}: \Fccc q(\sX_{\eta_w}) = \Fccc q(\sX_w\times_{\sS_w} \eta_w) \to 
\Fccc q(\sX_w\times_{\sS_w} \eta_{w,y}) \\
\to H^1_{\sX_y}(\sX_w\times_{\sS_w}\hen{(\sS_w)}{y},\Fccc q) \simeq (\Fccc q)_{-1}(\sX_y)=\Fccc {(q+1)}(\sX_y) ,
\end{multline}
where the isomorphism is induced by $\theta_{\ep,\sS'}$ in \eqref{eq4;puritycontraction} for 
$(\sX,\ep,\sS')=(\hen{(\sS_w)} y, \lambda_y,\sX_y)$ using \eqref{eq6-6;cousincomplex}. Note
\begin{equation}\label{eq7-1-2;cousincomplex}
\Ker(\partial^q_{w,y})\supset \Fccc q(\sX_w\times_{\sS_w} \sS_{w,y}),
\end{equation}
where $\sS_{w,y}$ is the localization of $\sS_w$ at $y$.
The commutativity of the left square is obvious.
To show that of the right square, note the following  commutative diagrams:
\begin{equation}\label{eq7-2;cousincomplex}
\xymatrix{
\Fccc q(\sX_{\eta_w}) \ar[r]^{\hskip -15pt\theta_{\phi,\lambda_w}} \ar[d] &
H^q_{\sX_{\eta_w}}(\sX_{Nis},F_\sX)\ar[d]\\
\Fccc q(\xi_w) \ar[r]^{\hskip -15pt\theta_1} & H^q_{\xi_w}(\sX_{Nis},F_\sX),\\
}\end{equation}
\begin{equation}\label{eq7-3;cousincomplex}
\xymatrix{
\Fccc {(q+1)}(\sX_y) \ar[r]^{\hskip -10pt\theta_{\phi,\lambda_w,\lambda_y}} \ar[d] &
H^{q+1}_{\sX_y}(\sX_{Nis},F_\sX)\ar[d]\\
\Fccc {(q+1)}(\xi_y) \ar[r]^{\hskip -15pt\theta_2} & H^{(q+1)}_{\xi_y}(\sX_{Nis},F_\sX),\\
}\end{equation}
where $\xi_w$ (resp. $\xi_y$) is the generic point of $\sX_w$ (resp. $\sX_y$) and $\theta_1$ (resp. $\theta_2$) is the isomorphism induced by $\theta_\ep$ in
\eqref{eq1;puritycontraction} for $\sX=\hen{(\sX)} {\xi_w}$ and 
$\ep:\hen{(\sX)} {\xi_w}\simeq \Spec k(\xi_w)\{\tauul_w\}$ induced by \eqref{eq6-2;cousincomplex} (resp. for $\sX=\hen{(\sX)} {\xi_y}$ and 
$\ep:\hen{(\sX)} {\xi_y}\simeq \Spec k(\xi_y)\{\tauul_w,\tau_y\}$ induced by \eqref{eq6-8;cousincomplex}). The commutativity follows from \eqref{eq4-2;puritycontraction}.
It also implies the commutativity of the following diagram:
\begin{equation}\label{eq7-4;cousincomplex}
\xymatrix{
\Fccc {q}(\sX_{\eta_w}) \ar[r]^{\hskip -15pt\partial^q_{w,y}} \ar[d] 
&\Fccc {(q+1)}(\sX_y)\ar[d]\\
\Fccc {q}(\xi_w) \ar[r]^{\hskip -15pt\delta^q_{w,y}}  &\Fccc {(q+1)}(\xi_y),\\
}\end{equation}
where $\delta^q_{w,y}$ is the composite map (cf. \eqref{eq7-1-1;cousincomplex})
\begin{equation*}
\Fccc q(\xi_w) \to H^1_{\xi_y}(\hen{(\sX_w)} {\xi_y},\Fccc q) \simeq 
(\Fccc q)_{-1}(\xi_y)=\Fccc {(q+1)}(\xi_y) ,
\end{equation*}
where the isomorphism is induced by $\theta_\ep$ in 
\eqref{eq1;puritycontraction} for $\sX=\hen{(\sX_w)} {\xi_y}$ and 
$\ep: \hen{(\sX_w)} {\xi_y}\simeq \Spec k(\xi_y)\{\tau_y\}$ 
induced by \eqref{eq6-6;cousincomplex}.
The right vertical arrows in \eqref{eq7-2;cousincomplex} and
\eqref{eq7-3;cousincomplex} are injective by \eqref{eq0;proofvanishing}. 
Thus the desired commutativity follows from that of 
\[\xymatrix{
\Fccc {q}(\xi_w) \ar[r]^{\hskip -15pt\theta_1} \ar[d]^{\delta^q_{w,y}} &
H^{q}_{\xi_w}(\sX_{Nis},F_\sX)\ar[d]\\
\Fccc {(q+1)}(\xi_y) \ar[r]^{\hskip -15pt\theta_2} & H^{q+1}_{\xi_y}(\sX_{Nis},F_\sX),\\
}\]
which follows from \eqref{eq2;puritycontraction}
for $\sX=\hen{(\sX)} {\xi_y}\simeq \Spec k(\xi_y)\{\tauul_w,\tau_y\}$
with $\sZ_{[1,i]}=\{\tauul_w=0\}$ and $\sZ_{[1,i+1]}=\{\tauul_w=\tau_y=0\}$.

\medbreak\noindent{\it Proof of Claim \ref{claim;vanishing}:}\;
Write $\Lambda=\Lambda_\alpha\subset \sS$ for simplicity.
Let $\Lambda_\sT\subset \sT$ be the image of $\Lambda\subset \sS$.
Note $\Lambda_\sT\subset \cd \sT {q-1}$ since the induced map $\ol{\Lambda}\to \sT$ is finite.
For $w\in \Lambda_\sT$ write $\Lambda_{w}=\Lambda\times_\sT w$.
By \eqref{eq7-1-2;cousincomplex} the composite of the natural map
\[\Fccc {(q-1)} (\sX_w-\underset{y\in \Lambda_{w}}{\cup}\sX_y)\to
\Fccc {(q-1)} (\sX_{\eta_w})\]
and $\partial_{w,y}^{q-1}$ from \eqref{eq7-1-1;cousincomplex} is the zero map for 
$y\not\in \Lambda_w$. Hence we get the induced map 
\[ \partial^{q-1}_{\Lambda_w} : 
\Fccc {(q-1)} (\sX_w-\underset{y\in \Lambda_{w}}{\cup}\sX_y)\to
\underset{y\in \Lambda_{w}}{\bigoplus}\;\Fccc {q}(\sX_y).\]
By the commutativity of \eqref{eq7-1;cousincomplex} with $q$ replaced by $q-1$,
it now suffices to show that
$\partial^{q-1}_{\Lambda_w} $ is surjective. Let 
\[\sY=\sT\{\tul\}=\sT\{t_1,\dots,t_{i-1}\}=\hen {(\sT[\tul])} {(\ul0_\sT,\tul)}.\] 
Writing $s=t_i$, we have $\sX=\sY\{s\}=\sS\{\tul\}$ and $\sS=\sT\{s\}$.
We have an essentially \'etale map 
$\pi:\sX \to \sS\times_\sT \sY$ fitting in a commutative diagram
\begin{equation}\label{eq7-5a;cousincomplex}
\xymatrix{
\sX \ar[r]^{\hskip -20pt\pi}\ar[d]^{\phi_t} & \sS\times_\sT \sY \ar[r] \ar[d] & \sY[s] \ar[r]\ar[d] &\sY\ar[d]\\
\sS[\tul] \ar[r] &\sS \ar[r]^{\psi_s} &  \sT[s] \ar[r] & \sT \\
}
\end{equation}
and $\sX$ is identified with $\hen {(\sY[s])}{(\ul0_\sY,s)}$,
where $\ul0_\sY\in \sY$ is the closed point.
Since $\ol{\Lambda}\to \sT$ is finite by the assumption of Claim \ref{claim;vanishing}, 
$\pi$ induces an isomorphism 
$\sX\times_\sS \ol{\Lambda} \simeq \ol{\Lambda}\times_\sT \sY$.
Letting $\sY_w=\sY\times_\sT w$, we get isomorphisms for $y\in \Lambda_w$:
\begin{equation}\label{eq7-5;cousincomplex}
 \sX_y =\sX\times_\sS y \simeq y\times_\sT\sY = y\times_w \sY_w .
\end{equation}
The map $\pi$ from \eqref{eq7-5a;cousincomplex} induces 
$\pi_w: \sX_w \to \sS_w\times_w \sY_w$. 
By the definition of $\partial^{q-1}_{\Lambda_w}$ using \eqref{eq6-6;cousincomplex}, we have a commutative diagram:
\[\xymatrix{
\Fccc {(q-1)}(\sS_w\times_w \sY_w-\Lambda_w\times_w \sY_w) \ar[r]^-{\pi_w^*} 
\ar[d]^{\pi_{\sS_w}^*}
& \Fccc {(q-1)} (\sX_w-\underset{y\in \Lambda_{w}}{\cup}\sX_y)
\ar[ddd]^-{\partial^{q-1}_{\Lambda_w}} \\ 
\frac{\Fccc {(q-1)}((\sS_w)^h_{|\Lambda_w} \times_w \sY_w-\Lambda_w\times_w \sY_w)}{\Fccc {(q-1)}((\sS_w)^h_{|\Lambda_w} \times_w \sY_w) } 
\ar[d]^{\simeq}_{(*1)} \\
\underset{y\in \Lambda_{w}}{\bigoplus}\;
\frac{\Fccc {(q-1)}((\sS_w)^h_{|y} \times_w \sY_w-y\times_w \sY_w) }
{\Fccc {(q-1)}((\sS_w)^h_{|y} \times_w \sY_w) }
\ar[d]^{\simeq}_{(*2)}\\
\underset{y\in \Lambda_{w}}{\bigoplus}\;
\frac{\Fccc {(q-1)}((\sS_w)^h_{|y} \times_y \sX_y-y\times_y \sX_y)}
{\Fccc {(q-1)}((\sS_w)^h_{|y} \times_y \sX_y) } 
\ar[r]_-{\oplus_y \theta_y}^-{\simeq} 
& \underset{y\in \Lambda_{w}}{\bigoplus}\;\Fccc {q}(\sX_y).
}\]
where $\pi_{\sS_w}: (\sS_w)^h_{|\Lambda_w} \to \sS_w$ is the natural map and
$\theta_y$ is the isomorphism $\theta_{\ep,\sS'}$ in \eqref{eq4;puritycontraction} for 
$(\sX,\ep,\sS')=(\hen{(\sS_w)} y, \lambda_y,\sX_y)$.
The isomorphism $(*1)$ (resp. $(*2)$) follows from the fact 
$(\sS_w)^h_{|\Lambda_w} = \prod_{y\in\Lambda_w} (\sS_w)^h_{|y}$ by \cite[Ch. XI, Th. 1]{ray} (resp. \eqref{eq7-5;cousincomplex}).
Hence it suffices to show the surjectivity of $\pi_{\sS_w}^*$.

For a Nisnevich neighbourhood $(S,\ul0_S)\to (\sT[s],(\ul0_\sT,s))$, let 
$\Lambda_S\subset S$ be the image of $\Lambda$ and put 
$\Lambda_{S,w}=\Lambda_S \times_\sT w$ for $w\in \Lambda_\sT$.
By Lemma \ref{lem3;purity} there exists a cofinal system of \'etale neighbourhoods $(S,\ul0_S)\to (\sT[s],(\ul0_\sT,s))$ such that
$(S_w, \Lambda_{S,w})$ is a nice $V$-pair over $w$.
Hence $(S_w\times_w \sY_w, \Lambda_{S,w}\times_w \sY_w)$ is a nice $V$-pair over $\sY_w$ by Remark \ref{rem;def;Vtriple}(2).
By Corollary \ref{cor;thm2;Vtriple} and Lemma \ref{lem6;Vtriple}, the natural map
\begin{small}
\[\frac{ \Fccc {(q-1)} (S_w\times_w \sY_w-\Lambda_{S,w}\times_w \sY_w)}
{\Fccc {(q-1)} (S_w\times_w \sY_w) } \to 
\frac{\Fccc {(q-1)} ((S_w)^h_{|\Lambda_{S,w}} \times_w \sY_w-\Lambda_{S,w}\times_w \sY_w)}{\Fccc {(q-1)} ((S_w)^h_{|\Lambda_{S,w}} \times_w \sY_w) } \]
\end{small}
is an isomorphism.
The desired surjectivity of $\pi_{\sS_w}^*$ follows by taking the colimit over the above cofinal system. This completes the proof of Claim \ref{claim;vanishing}.
\medbreak

\begin{cor}\label{cor;vanishing}
Take $F\in \lsCItspNis$.
\begin{itemize}
\item[(1)]
Let $X\in \Sm$ and $x\in \cd X n$ with $n\in \Z_{>0}$ and $K=k(x)$.
Then 
\[  H^i_x(X_{\Nis},F_X)=0\qfor i\not=n,\]
and there exists an isomorphism
\[ \theta_\ep: \Fccc n(x) \simeq H^n_x(X_{\Nis},F_X),\]
which depends on an isomorphism (see Lemma \ref{lem;etalepoint}):
\[\ep :\hen{X}{x}\simeq  \Spec K\{t_1,\dots,t_n\}.\]
\item[(2)]
Let $X\in \Sm$ and $x\in \cd X n$ with $n\in \Z_{>0}$.
Let $D\subset X$ be a regular closed subscheme.
For $e\in \Z_{\geq 0}$ we have 
\[  H^i_x(X_{\Nis},F_{(X,eD)})=0\qfor i\not=n.\]
\item[(3)]
Let $(X,D)\in \ulMCorls$ and $i:Z\hookrightarrow  X$ be a closed immersion of pure codimension $q$ such that $Z$ is regular and transversal with $D_1\cap\cdots\cap D_r$ for any $r>0$ and any distinct irreducible components $D_1,\dots,D_r$ of $D$.
Then $R^\nu i^! F_{(X,D)} =0$ for $\nu\not=q$.
\end{itemize}
\end{cor}
\begin{proof}
(1) follows from Theorem \ref{thm;vanishing} and \eqref{eq0;puritycontraction} and
\eqref{eq1;puritycontraction}. As for (2), we may assume $x\in D$ and replace $X,D$ by its henselization at $x$. By Lemma \ref{lem4;purity}(1), there is an exact sequence of sheaves on $X_{\Nis}$:
\begin{equation}\label{eq3;lem4;purity}
 0\to F_X \to F_{(X,eD)} \to \iota_*(\Fcontt e)_{D} \to 0,
\end{equation}
where $\iota:D\to X$ is the closed immersion. 
Note $\Fcontt e\in \lsCItspNis$ by Lemma \ref{lem1;contraction}. Looking at the long exact sequence of cohomology arising from \eqref{eq3;lem4;purity}, (2) gives 
$H^i_x(X_{\Nis},F_{(X,eD)})=0$ for $i\not=n,n-1$. It also gives isomorphisms
\[ \theta_\epsilon: H^n_x(X_{\Nis},F_X)\simeq \Fccc n(x),\;\;
\theta_{\epsilon_D}: H^n_x(D_{\Nis},(\Fcontt e)_{D} )\simeq (\Fcontt e)_{-(n-1)}(x) ,\]
depending on a chosen isomorphism
\[\ep :\hen{X}{x}\simeq  \Spec K\{t_1,\dots,t_n\}\]
such that $\hen D x\subset \hen X x$ is defined by $t_1=0$, and on the induced isomorphism
\[\ep_D :\hen{D}{x}\simeq  \Spec K\{t_2,\dots,t_n\}.\]
This gives rise to an exact sequence
\begin{multline*}
0\to H^{n-1}_x(X_{\Nis},F_{(X,eD)}) \to (\Fcontt e)_{-(n-1)}(x) \rmapo{\partial}
\Fccc n(x) \\ \to H^{n}_x(X_{\Nis},F_{(X,eD)}) \to 0,
\end{multline*}
where $\partial$ is the boundary map coming from \eqref{eq3;lem4;purity}.
By \eqref{eq2;puritycontraction} (and the same argument as the proof of \eqref{eq7-1;cousincomplex}), $\partial$ is identified with the map induced by $\Fcontt e \to \Fcont$ under the identification $\Fccc n=(\Fcont)_{-(n-1)}$. 
Since the latter map is injective by the semipurity of $F$
and \eqref{eq1;def2;contraction} ,
$\partial$ is injective and hence $H^{n-1}_x(X_{\Nis},F_{(X,eD)})=0$.
This completes the proof of (2).

As for (3), if $D=\emptyset$, it follows from Theorem \ref{thm;vanishing} and \eqref{eq0;puritycontraction}. 
In general we proceed by the induction on the number of the irreducible components of $D$. We may replace $X$ by its henselization at $x$. 
Let $E$ be one of the irreducible components of $D$ and write
$D=D' + e E$ with $e>0$ and $E\not\subset D'$.
By Lemma \ref{lem4;purity}(1), there is an exact sequence of sheaves on $X_{\Nis}$:
\begin{equation}\label{eq3b;lem4;purity}
 0\to F_{(X,D')} \to F_{(X,D)} \to \iota_*(\Fcontt e)_{(E,E\cap D')} \to 0,
\end{equation}
where $\iota:E\to X$ is the closed immersion. It gives rise to a long exact sequence of sheaves on $Z_\Nis$:
\begin{equation}\label{eq3c;lem4;purity} 
\dots \to R^\nu i^! F_{(X,D')} \to R^\nu i^! F_{(X,D)} \to R^\nu i^! \iota_*(\Fcontt e)_{(E,E\cap D')} \to \cdots.
\end{equation}
and we have an isomorphism
\[ R^\nu i^! \iota_*(\Fcontt e)_{(E,E\cap D')} \simeq
(\iota_Z)_*R^\nu i_E^! (\Fcontt e)_{(E,E\cap D')},\]
where $i_E: Z\cap E\to E$ and $\iota_Z:  Z\cap E \to Z$ are the closed immersions.
Noting $\Fcontt e\in \lsCItspNis$ and that the triple $(E,D\cap E,Z\cap E)$ satisfies the same assumption as that of $(X,D,Z)$, we get
\[ R^\nu i^! F_{(X,D')} =R^\nu i_E^! (\Fcontt e)_{(E,E\cap D')}=0\qfor \nu\not=q\]
by the induction hypothesis. This implies (3) by \eqref{eq3c;lem4;purity}.

\end{proof}

\section{$\bcube$-invariance of cohomology presheaves}\label{cubeinvcoh}

\begin{thm}\label{thm2-P1}
Let $\eta$ be the generic point of an integral $S\in \Sm$ and $Z\subset \P^1_\eta$ be an effective Cartier divisor such that $\infty\in |Z|$. For $F\in\lsCItspNis$, we have
\[H^i((\P^1_\eta)_{\Nis},F_{(\P^1,Z)})=0 \qfor i>0.\]
\end{thm}

We need a preliminary lemma for the proof.

\begin{lemma}\label{lem;thm2-P1}
Let $\eta$ be as in Theorem \ref{thm2-P1}.
For $F\in\lsCItspNis$ and $n\in \Z_{>0}$, the natural map
\[ \frac{F(\P^1_\eta, n0+\infty)}{F(\P^1_\eta,\infty)} \to 
\frac{F(\hen{(\A^1_\eta)} 0,n0)}{F(\hen{(\A^1_\eta)} 0)}\]
is an isomorphism.
\end{lemma}
\begin{proof}
By Lemmas \ref{lemP1-0} and Lemma \ref{lemP1-1}(2) the natural map
\[  F(\A^1_\eta,n0)/F(\A^1_\eta)\to F(\hen{(\A^1_\eta)} 0,n0)/F(\hen{(\A^1_\eta)} 0)\]
is an isomorphism. Lemma \ref{lem;thm2-P1} now follows from Lemma \ref{lem2;contraction}.
\end{proof}

We now prove Theorem \ref{thm2-P1}.
For simplicity write $\P^1_\eta=\P^1$ and $\A^1_\eta=\A^1$.
We only need to prove the vanishing of $H^1$. 
By the semipurity of $F$,
Corollay \ref{coro;localinjectivity} implies $F_{(\P^1,\infty)}\to F_{(\P^1,Z)}$ is injective. Its cokernel is supported on $Z$. Hence it suffices to show the vanshing in case $Z=\infty$. 
We have a localization exact sequence
\begin{multline*}
F(\P^1-0,\infty) \rmapo{\delta} 
\frac{F(\hen{(\A^1)} 0 -0)}{F(\hen{(\A^1)} 0)}\to 
H^1((\P^1)_{\Nis},F_{(\P^1,\infty)})  \\ \to H^1((\P^1-0)_{\Nis},F_{(\P^1-0,\infty)}).
\end{multline*}
The last term vanishes by Theorem \ref{thm-P1}(2) and $\delta$ is surjective by
Lemma \ref{lem;thm2-P1} and Lemma \ref{lem1;Mrec}(1). This proves the desired vanishing.
$\square$
\medbreak



\bigskip

\begin{thm}\label{thm2;purityM}
For $F\in\lsCItspNis$ and $(X,D)\in \ulMCorls$, 
\begin{equation}\label{eq;thm2;purityM} 
\pi^*: H^q(X_{\Nis},F_{(X,D)}) \isom H^q((X\times\P^1)_{\Nis},F_{(X,D)\otimes \bcube}) 
\end{equation}
induced by the projection $\pi:(X,D)\otimes \bcube \to (X,D)$.
\end{thm}
\begin{proof}
We proceed by induction on $\dim(X)$. If $\dim(X)=0$, \eqref{eq;thm2;purityM} follows from Theorem \ref{thm2-P1}. Assume $\dim(X)>0$.
By considering the Leray spectral sequence for $\pi$, we may replace $X,D$ by its henselization at $x$. 

\begin{claim}\label{claim;thm2;purityM}
We may assume $D=\emptyset$.
\end{claim}
\begin{proof}
Let $D_1,\dots,D_r$ be the irreducible components of $D$ and $e_i$ be the multiplicity of $D_i$ in $D$. Put $D'=\underset{2\leq i\leq r}{\sum} e_i D_i$ and $E=D_1\cap D'$.
By Lemma \ref{lem4;purity}(1), there are exact sequences of sheaves on $X_{\Nis}$
and $(X\times\P^1)_{\Nis}$ respectively:
\begin{equation*}
\begin{aligned} 
& 0\to F_{(X,D')} \to F_{(X,D)} \to \iota_*(\Fcontt {e_1})_{(D_1,E)} \to 0,\\
& 0\to F_{(X,D')\otimes\bcube } \to F_{(X,D)\otimes\bcube} \to \iota_*(\Fcontt {e_1})_{(D_1,E)\otimes\bcube} \to 0,\\
\end{aligned}
\end{equation*}
which are compatible in an obvious sense, 
where $\iota:D_1\hookrightarrow X$ is the closed immersion. Note $\Fcontt e\in \lsCItspNis$ by Lemma \ref{lem1;contraction}. By the induction hypothesis \eqref{eq;thm2;purityM} holds for
$(D_1,E)$ and hence we are reduced to showing it for $(X,D')$.
Repeating the same argument, this proves the claim.
\end{proof}

Now assume $D=\emptyset$. The idea of the following proof comes from \cite[Lecture 24]{mvw}. Take a regular divisor $i:Z\hookrightarrow X$ and let $j:U=X-Z\hookrightarrow X$ be the open immersion. We have a commutative diagram
\[\xymatrix{
H^q((X\times\P^1)_\Nis,F_{X\otimes\bcube})\ar[r]^{j^*} \ar[d]^{i_0^*} &
H^q((U\times\P^1)_\Nis,F_{U \otimes\bcube})\ar[d]^{i_0^*}  \\
H^q(X_\Nis,F_X) \ar[r] & H^q(U_\Nis,F_U)\\
}\]
where the vertical maps are the pullback along the $0$-section of $\P^1$.
It suffices to show the injectivity of $i_0^*$ on the left hand side since it is a left inverse of $\pi^*$.
Noting $\dim(U)<\dim(X)$ ($X$ is local) and using the Leray spectral sequence, 
the induction hypothesis implies that the right vertical map is an isomorphism. 
Thus we are reduced to showing the injectivity of $j^*$.
By Corollary \ref{cor;vanishing}(3), we get
\[R^\nu j_* F_{U\otimes\bcube}\simeq R^{\nu+1} i^! F_{X\otimes\bcube}=0
\qfor \nu\not=0.\]
By Lemma \ref{lem4;purity}(1), there is an exact sequence of sheaves on $(X\times\P^1)_{\Nis}$:
\begin{equation*}
 0\to F_{X \otimes\bcube } \to j_* F_{U\otimes\bcube} \to i_*(\Fcont)_{Z\otimes\bcube} \to 0.
\end{equation*}
Hence the desired injectivity follows from the surjectivity of
\[H^{q-1}((U\times\P^1)_\Nis,F_{U \otimes\bcube})\to 
H^{q-1}((Z\times\P^1)_\Nis,(\Fcont)_{Z \otimes\bcube}).\]
We have a commutative diagram
\[\xymatrix{
H^{q-1}((U\times\P^1)_\Nis,F_{X\otimes\bcube})\ar[r]\ar[d]^{i_0^*} &
H^{q-1}((Z\times\P^1)_\Nis,(\Fcont)_{Z\otimes\bcube})\ar[d]^{i_0^*}  \\
H^{q-1}(U_\Nis,F_U) \ar[r]^{\beta} & H^{q-1}(Z_\Nis,(\Fcont)_Z)\\
}\]
where the right vertical map is an isomorphism by the induction hypothesis.
If $q>1$, $H^{q-1}(Z_\Nis,(\Fc)_Z)=0$ since $Z$ is henselian local.
If $q=1$, $\beta$ is surjective by Lemma \ref{lem4;purity}(1).
This completes the proof.
\end{proof}

\section{Sheafication preserves $\bcube$-invariance}\label{cubeNis}

In this section we prove the following.

\begin{thm}\label{thm;sheafication}
If $F\in \CItsp$, $\ulaNis F\in \CItsp$ (cf. Definition \ref{def;Xi})
\end{thm}
\medbreak

We need a preliminary for the proof.
For $F\in \MPST$ put
\[\hF=\tau_!\hMMnp\omega^*\omega_! F\;\in \ulMPST
\text{ (cf. \eqref{eq;adjunction} and \eqref{eq;hMM})}.\]
By Lemma \ref{lem;MPST}(5) and (1), for $M\in \ulMCor$ we have
\begin{equation}\label{eq1;sheafication}
\begin{aligned}
\hF(M) &= \colim_{N\in \Comp(M)} \Hom_{\MPST}(\hM N,\omega^*\omega_! F)\\
& \simeq \colim_{N\in \Comp(M)} \Hom_{\PST}(h_0(N),\omega_! F)\\
& \simeq \colim_{N\in \Comp(M)} \Hom_{\PST}(h_0(N),\ulomega_! \tau_!F),
\end{aligned}
\end{equation}
where $h_0(N)=\omega_!h_0^\bcube(N)\in \PST$. By Lemma \ref{lem;MPST}(1)
the unit map $u:F \to \omega^*\omega_! F$ from \eqref{eq;adjunction} induces 
\[\gamma: \tau_!F \to \tau_!\omega^*\omega_! F=\ulomega^*\ulomega_!\tau_!F.\]
Assume $F\in \CI$ (cf. Definition \ref{def:CI}). By Lemma \ref{lem;hMM},
$u$ factors as
\[ F \to \hMMnp\omega^*\omega_! F \to \omega^*\omega_! F,\]
where the second map is injective. It gives a factorization of $\gamma$ as
\begin{equation}\label{eq1.5;sheafication}
\tau_!F \rmapo{\hgamma}  \hF \rmapo{\hgamma'} \ulomega^*\ulomega_!\tau_!F
\end{equation}
where $\hgamma'$ is injective by the exactness of $\tau_!$ (cf. Lemma \ref{lem;MPST}(1)). Assume further that $\tau_!F$ is \lssemipure. Then
\begin{equation}\label{eq2;sheafication}
 \tau_!F(\fX)  \rmapo{\hgamma}  \hF(\fX)\; \text{ is injective for any $\fX\in \ulMCorls$.}
\end{equation}
By Lemma \ref{lem;MPST}(2), we have $ \ulomega^*\ulomega_!\tau_!F(X)\simeq F(X)$ for $X\in \Sm$. By \eqref{eq1.5;sheafication} we get isomorphisms
\begin{equation}\label{eq3;sheafication}
\tau_! F(X) \simeq \hF(X)\simeq \ulomega^*\ulomega_!\tau_!F(X) \quad(X\in \Sm).
\end{equation}


\begin{lemma}\label{lem1;sheafication}
Let $G\in \CItsp$ and $\eta$ be the generic point of $S\in \Sm$ and $(X,0_X) \to (\A_\eta^1,0)$ be a
Nisnevich neighbourhood of $0\in \A_\eta^1$. Consider
\[ \phi: G(\P_\eta^1-0,\infty) \to G(X-0_X)/G(X). \]
Then we have $G(\eta) \simeq \Ker(\phi)$.
\end{lemma}
\begin{proof}
By Lemma \ref{lem;CItau}, $G=\tau_!F$ with $F\in \CI$. 
In view of \eqref{eq2;sheafication} and \eqref{eq3;sheafication},
it suffices to show that the assertion for $G=\hF$. By \eqref{eq1;sheafication}
\[ \hF(\P_\eta^1-0,\infty) \simeq \colim_{n\in \Z_{>0}} 
\Hom_{\PST}(h_0(\P^1_\eta,\infty+n 0),\ulomega_! \tau_!F).\]
Put $E=\ulomega_! \tau_!F\in \PST$. We claim
\begin{equation}\label{eq4;sheafication}
E(U)\simeq E_{\Nis}(U)\quad \text{for any open $U\subset \A^1_\eta$.}
\end{equation}
where $E_{\Nis}=\aVNis E \in \NST$ (cf. \eqref{eq1;omegaNis}).
Indeed we have
\[ E_{\Nis}(U) \simeq \ulomega_! \ulaNis(\tau_!F) (U)
 = \ulaNis(\tau_!F) (U,\emptyset) \simeq \tau_!F(U,\emptyset)= E(U),\]
where the first (resp. second) isomorphism follows from \eqref{eq1;omegaNis}
(resp. Theorem \ref{thm-P1}(1) and Remark \ref{rem;Nissheafication}).
For $U=\A^1_\eta-0$, we have a commutative diagram
\[\xymatrix{
E(U) \ar[r]^{\hskip -40pt\simeq} \ar[d]^{\simeq}_{\eqref{eq4;sheafication}} 
&\Hom_{\PST}(\Ztr(U),E)\ar[d] ^{\simeq}
& \ar[l]_{\iota} \Hom_{\PST}(h_0(\P^1_\eta,\infty+ n0),E)\ar[d]\\
E_{\Nis}(U) \ar[r]^{\hskip -40pt\simeq} & \Hom_{\NST}(\Ztr(U),E_{\Nis})
& \ar[l] \Hom_{\NST}(h_0(\P^1_\eta,\infty+ n0)_{\Nis},E_{\Nis})
}\]
where $\iota$ is injective since $h_0(\P^1_\eta,\infty+n 0)$ is a quotient in $\PST$ of $\Ztr(\A^1_\eta-0)$. Hence the right vertical map of the diagram is injective.
Hence we get a natural injection
\begin{equation}\label{eq5;sheafication}
 \hF(\P^1_\eta-0,\infty) \hookrightarrow \colim_{n\in \Z_{>0}} 
\Hom_{\NST}(h_0(\P^1_\eta,\infty+n 0)_{\Nis},E_{\Nis}).
\end{equation}
We have a commutative diagram
\begin{equation}\label{eq6;sheafication}
\xymatrix{
&&\hF(\P^1_\eta-0,\infty) \ar[r]^{\hskip -20pt\phi}\ar[d]^{\hookrightarrow} & \hF(X-0_X)/\hF(X)\ar[d]^{\eqref{eq3;sheafication}}_{\simeq} \\
0\ar[r] & \tau_!F(\A^1_\eta) \ar[r] &\tau_!F(\A^1_\eta-0) \ar[r] & \tau_!F(X-0_X)/\tau_!F(X)\;,
}
\end{equation}
where the lower sequence is exact by Lemmas \ref{lemP1-1} and \ref{lemP1-0}
and the left vertical injection is induced by $\hgamma'$
from \eqref{eq1.5;sheafication} in view of Lemma \ref{lem;MPST}(2).
We have a commutative diagram
\[\xymatrix{
\colim_{m\in \Z_{>0}} \Hom_{\PST}(h_0(\P^1_\eta,m\infty),E) \ar[r]\ar[d]^{\alpha} &
\colim_{m\in \Z_{>0}} \Hom_{\NST}(h_0(\P^1_\eta,m\infty)_{\Nis},E_{\Nis}) 
\ar[d]^{\beta} \\
\Hom_{\PST}(\Ztr(\A^1_\eta),E) \ar[d]^{\simeq} & 
\Hom_{\NST}(\Ztr(\A^1_\eta),E_{\Nis}) \ar[d]^{\simeq} \\
E(\A^1_\eta) \ar[r]^{\simeq}_{\eqref{eq4;sheafication}}&  E_{\Nis}(\A^1_\eta) \\
}\]
where $\alpha$ is an isomorphism by Lemma \ref{lem;CIRec}.
Since $h_0(\P^1_\eta,m\infty)_{\Nis}$ is a quotient in $\NST$ of $\Ztr(\A^1_\eta)$, $\beta$ is injective. Hence the diagram implies that $\beta$ is an isomorphism and we get an isomorphism
\begin{equation}\label{eq7;sheafication}
\tau_!F(\A^1_\eta) \simeq E(\A^1_\eta) \simeq
\colim_{m\in \Z_{>0}} \Hom_{\NST}(h_0(\P^1_\eta,m\infty)_{\Nis},E_{\Nis}).
\end{equation}
By \eqref{eq5;sheafication}, \eqref{eq6;sheafication} and \eqref{eq7;sheafication}, $\Ker(\phi)$ injects into 
\[\big(\colim_{m\in \Z_{>0}} \Hom_{\NST}(h_0(\P^1_\eta,m\infty)_{\Nis},E_{\Nis})\big)\cap
\big(\colim_{n\in \Z_{>0}} \Hom_{\NST}(h_0(\P^1_\eta,\infty+n 0)_{\Nis},E_{\Nis})\big),\]
where the intersection is taken in
\[\colim_{m,n\in \Z_{>0}}  \Hom_{\NST}(h_0(\P^1_\eta,m\infty+n 0)_{\Nis},E_{\Nis})\;
\subset E_{\Nis}(\A^1_\eta-0)=\tau_!F(\A^1_\eta-0).\]
We now claim that there is an exact sequence in $\NST$:
\begin{multline}\label{eq8;sheafication}
 h_0(\P^1_\eta,m\infty+n 0)_{\Nis}\to 
h_0(\P^1_\eta,\infty+n0 )_{\Nis} \oplus h_0(\P^1_\eta,m\infty)_{\Nis} \\
\to h_0(\P^1_\eta,\infty)_{\Nis} \to 0.
\end{multline}
Indeed \cite[Th.1.1]{ry} implies an isomorphism of Nisnevich sheaves on $\Sm$:
\begin{equation}\label{eq7.5;sheafication}
 h_0(\P^1_\eta,\frak m)_{\Nis} \simeq \underline{\Pic}(\P^1_\eta,\frak m),
\end{equation}
where $\fm$ is an effective divisor on $\P^1$ and $\underline{\Pic}(\P^1_\eta,\frak m)$ is the sheaf associated to the presheaf  
$U\to \Pic(\P^1_\eta\times U, \frak m\times U\times\eta)$ ($U\in \Sm$).
There is an exact sequence
\begin{equation}\label{eq7.6;sheafication}
 0\to \ulsO^\times \to \ulsO^\times_{\fm} \to \underline{\Pic}(\P^1_\eta,\frak m) \to \Z\to 0,
\end{equation}
where $\ulsO^\times(U)=\sO(U\times\eta)^\times$ and 
$\ulsO^\times_{\fm}(U)=\sO(U\times\fm\times\eta)^\times$ for $U\in\Sm$.
The claim follows easily from this.
By \eqref{eq7.5;sheafication} and \eqref{eq7.6;sheafication} we have isomorphisms
\[h_0(\P^1_\eta,\infty)_{\Nis}\simeq \underline{\Pic}(\P^1_\eta,\infty)  \simeq \Ztr(\eta).\]
Thus \eqref{eq8;sheafication} implies 
\[\Ker(\phi)\hookrightarrow \Hom_{\NST}(\Ztr(\eta),E_{\Nis})\simeq E(\eta)\simeq F(\eta),\]
which proves Lemma \ref{lem1;sheafication}.
\end{proof}

\medbreak\noindent
{\it Proof of Theorem \ref{thm;sheafication}:}\;
By Lemmas \ref{lem1;Mrec}(3) and \ref{lem1;semipure}(4), $\ulaNis F$ has $M$-reciprocity and semipure. Hence it suffices to show $\ulaNis F$ is $\cube$-invariant.
We start with the following.

\begin{claim}\label{claim;sheafication}
For the generic point $\eta$ of $S\in \Sm$,
the natural map $F(\eta) =F_{\Nis}(\eta)\to F_{\Nis}(\bcube\times\eta)$ is an isomorphism (see Definition \ref{def:Nissheafication} for $F_\Nis$). 
\end{claim}
\begin{proof}
Consider a commutative diagram
\[\xymatrix{
0\ar[r]& F(\P_\eta^1,\infty) \ar[r]\ar[d] & F(\P_\eta^1-0,\infty) \ar[r]\ar[d]^{\simeq}&
F(\sO^h_{\P_\eta^1,0}-0)/F(\sO^h_{\P_\eta^1,0})\ar[d]^{\simeq}\\
0\ar[r]& F_{\Nis}(\P_\eta^1,\infty) \ar[r]& F_{\Nis}(\P_\eta^1-0,\infty) \ar[r] & F_{\Nis}(\sO^h_{\P_\eta^1,0}-0)/F_{\Nis}(\sO^h_{\P_\eta^1,0})\\
}\]
The upper sequence is exact by Lemma \ref{lem1;sheafication} thanks to the assumption that $F$ is \lssemipure.
The lower sequence is exact: The injectivity of the first map follows from 
Theorem \ref{thm;localinjectivity} and the semipurity of $F_{\Nis}$ by Lemma \ref{lem1;semipure}(3). The exactness at the middle term comes from the sheaf property.
The right (resp. middle) vertical map is an isomorphism by an obvious reason 
(resp. Theorem \ref{thm-P1}(1)).
Thus we get $F_{\Nis}(\P_\eta^1,\infty)\simeq F(\P_\eta^1,\infty)\simeq F(\eta)$ as desired.
\end{proof}

We now follow the argument of \cite[22.1]{mvw}.
For $\fX=(\Xb,\Xinf)\in \ulMCor$ let
$i^*_\fX: \ulaNis F(\fX\otimes\bcube) \to \ulaNis F(\fX)$
be the pullback along $i_0:\Spec(k) \to \bcube$.
It suffices to show the injectivity of $i^*_\fX$.
Letting $\eta$ be the generic point of $\Xb$, it suffices to show the injectivity of the composite map 
\begin{equation}\label{eq;ifX}
 \ulaNis F(\fX\otimes\bcube) \to \ulaNis F(\fX)\to \ulaNis F(\eta) =
F_\Nis(\eta)=F(\eta).
\end{equation} 
Recall (cf. \ref{eq;ulaNisformular})
\[ \ulaNis F(\fX\otimes\bcube) =\colim_{\fY\in \Sigmafin\downarrow \fX\otimes\bcube}
F_\Nis(\fY).\]
Writing $\fY=(\Yb,\Yinf)\in \Sigmafin\downarrow \fX\otimes\bcube$, we have 
\[\Yb-|\Yinf| =  X\times \A^1 \qaq \fY\times_{\Xb} \eta=(\P^1_\eta,\infty_\eta),\]
where the second equality comes from the fact that any proper birational map
$W\to \P^1_\eta$ with $W$ normal, is an isomorphism.
Hence we have a commutative diagram
\[\xymatrix{
F_\Nis(\fY) \ar[r]^{\alpha}\ar[d]^{\gamma} & F_\Nis(X\times \A^1)\ar[d]^{\beta} \\
F_\Nis(\P^1_\eta,\infty_\eta) \ar[r] & F_\Nis(\A^1_\eta)\;.\\
}\]
The map $\alpha$ is injective since $F_\Nis$ is \lssemipure by Lemma \ref{lem1;semipure}(3) and $\beta$ is injective by Theorem \ref{thm;localinjectivity}(2). Hence $\gamma$ is injective. 
Note that the composite of \eqref{eq;ifX} and $F_\Nis(\fY) \to \ulaNis F(\fX\otimes\bcube)$ coincides with the composite of $\gamma$ and 
$F_\Nis(\P^1_\eta,\emptyset_\eta) \rmapo{i_\eta^*} F_\Nis(\eta)$, which is injective by Claim \ref{claim;sheafication}. This proves the desired injectivity of $\i_{\fX}^*$ and completes the proof of Theorem \ref{thm;sheafication}. .  


\section{Implications on reciprocity sheaves}\label{RSCproof}

In this section we prove Theorems \ref{thm0;rec} and \ref{thm;purityrec}.
We also deduce Voevodsky's theorem \ref{thm;Voev} from Theorem \ref{thm2;purityM}. We need a preliminary.

\medbreak

For $F\in \RSC$ we have isomorphisms
\begin{equation}\label{eq2.5;omegaNis}
 \ulomega_!\ulaNis \tau_!\omegaCI F \overset{\eqref{eq1;omegaNis}}{\simeq}
\aVNis \ulomega_! \tau_!\omegaCI F \overset{(*1)}{\simeq}
\aVNis  \omega_!\omegaCI F \overset{(*2)}{\simeq} \aVNis F,
\end{equation}
where $(*1)$ (resp. $(*2)$) follows from Lemma \ref{lem;MPST}(1) (resp. Lemma \ref{lem2;CIRec}). 

\begin{lemma}\label{lem;omegaCItau} 
For $F\in \PST$, we have $\tau_!\omegaCI F\in \ulMPST$ is $\bcube$-invariant and semipure with $M$-reciprocity.
\end{lemma}
\begin{proof}
By Lemmas \ref{lem;CItau} and Definition \ref{def:m-rec}(2), it suffices to show only the semipurity of $\tau_!\omegaCI F$, namely the injectivity of the unit map 
$u:\tau_!\omegaCI F \to \ulomega^*\ulomega_!\tau_!\omegaCI F$.
By Lemma \ref{lem;hMM} we have $\omegaCI F=h^0_{\bcube}\omega^*F\subset \omega^* F$.
Hence we have a commutative diagram
\[\xymatrix{
\tau_!\omegaCI F \ar[r]^{\hskip -16pt u}\ar[d]^{\hookrightarrow} & \ulomega^*\ulomega_!\tau_!\omegaCI F\ar[d]^{\hookrightarrow}\\
\tau_!\omega^* F \ar[r]^{\hskip -10pt u'} & \ulomega^*\ulomega_!\tau_!\omega^* F\\
}\]
where the vertical maps are injective thanks to the exactness of 
$\omega^*$, $\omega_!$ and $\tau_!$.
By Lemma \ref{lem;MPST} we have $\tau_!\omega^*\simeq \ulomega^*\simeq \ulomega^*\ulomega_!\tau_!\omega^*$ and
$u'$ is identified with the identity through these isomorphisms.
Hence $u$ is injective as desired.
\end{proof}
\medbreak

Now take $F\in \RSC$ and put $G=\tau_!\omegaCI F$.
By \eqref{eq3;omegaCI} and \eqref{eq2.5;omegaNis} we have natural isomorphisms
(note $\omega_!=\ulomega_!\tau_!$ by Lemma \ref{lem;MPST}) 
\begin{equation}\label{eq1;RSCproof}
F \simeq \ulomega_!G \qaq F_\Nis:=\aVNis F \simeq \ulomega_! \ulaNis G.
\end{equation}
By Lemma \ref{lem;omegaCItau} and Theorem \ref{thm;sheafication}, 
\begin{equation}\label{eq2;RSCproof}
\text{$G$ and $\ulaNis G$ are $\bcube$-invariant and semipure with $M$-reciprocity.}
\end{equation}
Hence Theorem \ref{thm0;rec} follows from Lemma \ref{lem;CIRec}.
By \eqref{eq2;omegaNis} we have 
\begin{equation}\label{eq1a;RSCproof}
  H^i(X_{\Nis},(F_{\Nis})_X)\simeq H^i(X_{\Nis},(\ulaNis G)_{(X,\emptyset)})
\end{equation}
and similarly for cohomology with support.
Hence Theorem \ref{thm;purityrec} follows from Corollary \ref{cor;vanishing} and Lemma \ref{lem3;contraction}.

\medbreak

Finally we deduce Theorem \ref{thm;Voev} from Theorem \ref{thm2;purityM}. 
Take $F\in \HI\cap \NST$. We claim 
\begin{equation}\label{eq6;RSCproof}
\ulomega^* F\in \CIt\cap \ulMNST.
\end{equation}
Indeed, by Lemma \ref{lem;MPST}(1) we have $\ulomega^* F =\tau_!\omega^* F$ 
and $\omega^* F\in \CI$ by Lemmas \ref{lem;omegaHICI}. Moreover $F\in \ulMNST$ by Lemma \ref{lem;ulMNSTomega}.
By definition \eqref{eq;adjunction}, for $X\in \Sm$, we have
$(\ulomega^* F)_{(X,\emptyset)} =  F_X$ as sheaves on $X_\Nis$, and hence we get 
a commutative diagram
\[\xymatrix{
 H^i(X_\Nis,(\ulomega^* F)_{(X,\emptyset)}) \ar[r]^{\sim}\ar[d]^{\simeq} &  H^i(X_\Nis,F_X)\ar[d] \\
  H^i((X\times\P^1)_\Nis,(\ulomega^* F)_{(X,\emptyset)\otimes\bcube}) \ar[r]^{\sim} &  
H^i((X\times\A^1)_\Nis,F_{X\times\A^1}), \\
}\]
where the horizontal maps are isomorphisms.
The left vertical map is an isomorphism by \eqref{eq6;RSCproof} and Theorem \ref{thm2;purityM}. This proves Theorem \ref{thm;Voev}.


\section{Appendix}

In this section we collect some technical lemmas used in this paper.

\begin{lemma}\label{lemA1}
Let $A$ be a local ring and $S=\Spec A$ with the closed point $s\in S$.
Let $\pb:\Xb\to S$ be a proper morphism.
Let $Z\subset \Xb$ be closed subschemes.
\begin{itemize}
\item[(1)]
$Z=\emptyset$ if and only if $Z\cap \pb^{-1}(s)=\emptyset$.
\item[(2)]
$Z$ is finite over $S$ if and only if $Z\cap \pb^{-1}(s)$ is finite.
\item[(3)]
Assume further that $S$ is henselian.
If $Z$ is irreducible and $Z\cap \pb^{-1}(s)$ is finite, then it consists of one closed point of $X$.
\end{itemize}
\end{lemma}
\begin{proof}
As for (1) the only-if part is obvious.
Assume $Z\not=\emptyset$. Then $\pb(Z)$ is not empty and closed in $S$ 
since $\pb$ is proper. Hence it contains $s$, which implies 
$Z\cap \pb^{-1}(s)\not=\emptyset$.

As for (2) the only-if part is obvious.
Assume $Z\cap \pb^{-1}(s)$ is finite. Put 
\[F=\{x\in Z|\; \dim_x (p^{-1}(p(x))\cap Z)\geq 1\}.\]
By Chevalley's theorem (see \cite[Th.2.1.1]{sv}),
$F$ is closed in $Z$ so that $\pb(F)$ is closed in $S$ by 
the properness of $p$.
Since $Z\cap \pb^{-1}(s)$ is finite, $s\not\in \pb(F)$, which is absurd since $S$ is local.
(3) follows from (2) and the fact that any finite scheme over
$S$ is the product of henselian local schemes.   
\end{proof}

\begin{lemma}\label{lemA3}
Let the assumption be as Lemmas \ref{lemA1} and assume $S$ is henselian.
Assume further $\dim(\pb^{-1}(s))=1$. Let $X\subset \Xb$ be an open subset such that
$\pb^{-1}(s)\cap X$ is dense in $\pb^{-1}(s)$. If $Z\subset X$ is closed and irreducible
and $Z\cap \pb^{-1}(s)$ is finite and non-empty, then $Z$ is finite over $S$ and
$Z\cap \pb^{-1}(s)$ consists of one closed point of $X$.
\end{lemma}
\begin{proof}
Let $\Zb\subset\Xb$ be the closure of $Z$. By the assumption $\Zb\cap \pb^{-1}(s)$ is finite so that $\Zb$ is finite over $S$ and $\Zb\cap \pb^{-1}(s)$ consists of one closed point of $X$ by Lemma \ref{lemA1}(2) and (3). 
Since $Z\cap \pb^{-1}(s)$ is non-empty, we must have 
$(\Zb-Z)\cap \pb^{-1}(s)=\emptyset$. Hence $\Zb-Z=\emptyset$ by Lemma \ref{lemA1}(1).
This completes the proof.
\end{proof}

\begin{lemma}\label{lem;Levine}
Let $S$ be either the spectrum of an infinite field or a henselian local ring 
with infinite residue field. Let $s\in S$ be the closed point.
Let $p:X\to S$ be smooth of relative dimension one.
Let $x\in X$ be a point such that $p(x)=s$ and $k(s) \simeq k(x)$.
Then there exists a Nisnevich neighbourhood $(X',x)$ of $(X,x)$ and 
a closed immersion $X'\hookrightarrow \A^N_S$ over $S$ such that
letting $\Xb'\subset \P^N_S$ be the closure of $X'$, $\Xb'-X'$ is finite over $S$.
\end{lemma}
\begin{proof}
This follows from \cite[Th.10.0.1]{lev}.
\end{proof}

\end{document}